\documentclass[10pt]{article}

\usepackage{amssymb,amsmath,latexsym}
\voffset
-1.5cm

\usepackage{a4,amscd,graphicx,epsfig}
\usepackage[all]{xy}
\usepackage[francais,english]{babel}
\usepackage[latin1]{inputenc}
\usepackage[T1]{fontenc}
\everymath{\displaystyle}

\newtheorem{teo}{Theorem}[section]
\newtheorem{teo1}[teo]{Théorème}
\newtheorem{lem}[teo]{Lemme}
\newtheorem{cor}[teo]{Corollaire}

\newtheorem{prop}[teo]{Proposition} 
\newtheorem{defi}[teo]{Définition} 
\newtheorem{defis}[teo]{Définitions}

\newtheorem{remark}[teo]{Remarque}
\newtheorem{remarks}[teo]{Remarques}
\newtheorem{conv}[teo]{Convention}

\newcommand{\mr}{\mathbb{R}}
\newcommand{\mc}{\mathbb{C}}
\newcommand{\mz}{\mathbb{Z}}

\title{Dilogarithme Quantique et $6$j-Symboles Cycliques}

\author {St\'ephane Baseilhac}

\date {}

\begin{document}

\maketitle

\noindent Laboratoire E. Picard, CNRS UMR 5580, UFR MIG, Universit\'e Paul Sabatier, 118 route de Narbonne, F-31062 TOULOUSE.

\medskip

\noindent Dipartimento di Matematica, Universit\`a di Pisa, Via
F. Buonarroti, 2, I-56127 PISA. Email: baseilha@mail.dm.unipi.it. 

\vspace{1cm}

\selectlanguage{francais}

\begin{abstract}

\medskip

\noindent Soit $\mathcal{W}_N$ une sous-algèbre de Borel quantifi\'ee de $U_q(sl(2,\mc))$, sp\'ecialis\'ee en une racine de l'unit\'e $\omega = \exp(2i\pi/N)$ d'ordre impair $N >1$. On montre que les $6j$-symboles des repr\'esentations cycliques de $\mathcal{W}_N$ sont des repr\'esentations de l'\'el\'ement canonique d'une extension convenable du double de Heisenberg de $\mathcal{W}_N$. Cet \'el\'ement canonique est un $q$-dilogarithme ``tordu''. En particulier, on donne des formules explicites pour ces $6j$-symboles, et on en construit des symm\'etrisations t\'etrah\'edrales partielles, les c-$6j$-symboles. Les c-$6j$-symboles sont à la base de la construction des invariants hyperboliques quantiques de $3$-vari\'et\'es.
\end{abstract}

\selectlanguage{english}

\vspace{-0.2cm}

\begin{abstract}

\noindent Let $\mathcal{W}_N$ be a quantized Borel subalgebra of $U_q(sl(2,\mc))$, specialized at a primitive root of unity $\omega = \exp(2i\pi/N)$ of odd order $N >1$. One shows that the $6j$-symbols of cyclic representations of $\mathcal{W}_N$ are representations of the canonical element of a certain extension of the Heisenberg double of $\mathcal{W}_N$. This canonical element is a twisted $q$-dilogarithm. In particular, one gives explicit formulas for these $6j$-symbols, and one constructs partial symmetrizations of them, the c-$6j$-symboles. The latters are at the basis of the construction of the quantum hyperbolic invariants of $3$-manifolds.
\end{abstract}

\selectlanguage{francais}

\section{Introduction}

\noindent R\'ecemment, l'auteur a construit en collaboration avec R. Benedetti une nouvelle famille d'invariants à valeurs complexes, les ``invariants hyperboliques quantiques'' (IHQ), pour des triplets $(M,L,\rho)$, o\`u $M$ est une $3$-vari\'et\'e ferm\'ee orient\'ee, $L \subset M$ est un entrelacs plong\'e, et $\rho$ est un fibr\'e principal et plat sur $M$, de groupe structural un sous-groupe de Borel de $SL(2,\mc)$ \cite{BB1}. La partie alg\'ebrique de cette construction repose sur les propri\'et\'es des $6j$-symboles des repr\'esentations cycliques d'une sous-algèbre de Borel quantifi\'ee $\mathcal{W}_N$ de $U_q(sl(2,\mc))$, sp\'ecialis\'ee en une racine de l'unit\'e $\omega = \exp(2i\pi/N)$ d'ordre impair $N >1$. (Dor\'enavant, on nommera ces $6j$-symboles les ``$6j$-symboles cycliques''). Ces propri\'et\'es des $6j$-symboles cycliques, dont la forme particulière de leur \'equation de pentagone, avaient \'et\'e d\'ecouvertes par L. Faddeev et R. Kashaev, et annonc\'ees dans \cite{K1,K3,K4}. 

\smallskip

\noindent Dans cet article on donne des preuves de ces propri\'et\'es, qu'on d\'eduit du r\'esultat suivant (implicitement sugg\'er\'e dans \cite{K3}): on montre que les $6j$-symboles cycliques sont des repr\'esentations de l'\'el\'ement canonique d'une extension convenable du double de Heisenberg de $\mathcal{W}_N$. Les doubles de Heisenberg (ddH) de bigèbres associatives et coassociatives ont \'et\'e introduits en théories des algèbres d'opérateurs \cite{BS,Nov} et des systèmes intégrables \cite{AM,STS,Lu}. Le ddH d'un groupe quantique $U_q(\mathfrak{g})$ (o\`u $\mathfrak{g}$ est une algèbre de Lie complexe et de dimension finie) est une algèbre associative qui peut être interprétée comme une déformation du fibré cotangeant $T^*(G)$ du groupe de Lie $G$ associé à $U_q(\mathfrak{g})$ \cite{STS}. Il est connu que les ddH gouvernent la structure de certaines catégories monoïdales \cite{CR,St,Mil,Dav}. Notre r\'esultat montre que ceci reste vrai pour les repr\'esentations cycliques de $\mathcal{W}_N$. Deux faits sont remarquables ici: on obtient des solutions \emph{à paramètres} de l'\'equation du pentagone, et l'\'el\'ement canonique du ddH consid\'er\'e est, à un facteur près, une version tronqu\'ee d'un $q$-dilogarithme (que nous appellerons un ``dilogarithme cyclique''). Par souci de clart\'e, on prouve tous les r\'esultats directement utilis\'es. Le comportement asymptotique des $6j$-symboles sera considéré dans un futur article.  

\medskip

\noindent  La structure de l'article est la suivante. Les deux premières sections contiennent des r\'esultats pr\'eliminaires. Au \S \ref{repcyc}, on \'etablit un paramétrage $\Psi$ des représentations irréductibles et cycliques de $\mathcal{W}_N$ dans un sous-groupe de Borel de $SL(2,\mc)$ (Proposition \ref{prop1a}); au \S \ref{asy}, on d\'ecrit quelques propri\'et\'es du $q$-dilogarithme.

\noindent Au \S \ref{double}, on construit les doubles de Heisenberg $\mathcal{H}_h(b(2,\mc))$ et $\mathcal{H}_{\omega}$ des algèbres $U_h(b(2,\mc))$ et $\mathcal{W}_N$ respectivement (Proposition \ref{defdoubleh} et D\'efinition \ref{defdoubleq}). On montre que l'élément canonique $R_h$ de $\mathcal{H}_h(b(2,\mc))$ d\'etermine un op\'erateur $R_\omega$ sur une extension convenable de $\mathcal{H}_{\omega}$. Au \S \ref{6jdilocyc} on d\'ecrit les opérateurs de Clebsch-Gordan et les $6j$-symboles des suites régulières de représentations cycliques de $\mathcal{W}_N$ (i.e. les $6j$-symboles cycliques) comme représentations cycliques de $R_\omega$ (Proposition \ref{CGdilocyc} et Corollaire \ref{6j1}). On en donne aussi des formules explicites, n\'ecessaires dans la section suivante (Propositions \ref{CG1} et \ref{6jdesc}).

\noindent Au \S \ref{etendu}, on construit des symmétrisations partielles des $6j$-symboles, les c-$6j$-symboles. Ces derniers sont des op\'erateurs d\'efinis sur des t\'etraèdres ``d\'ecor\'es'', i.e. munis de certaines structures suppl\'ementaires, et sont à la base de la construction des IHQ \cite{BB1}. On calcule les relations de symmétrie t\'etraèdrales des c-$6j$-symboles et on prouve quelques unes de leurs \'equations fonctionnelles remarquables, dont l'équation du pentagone étendue.

\bigskip

\noindent {\bf Remerciements.} Ce travail doit beaucoup aux conseils et critiques de Francis Bonahon et Christian Kassel. 

\section{Les représentations cycliques irr\'eductibles de $\mathcal{W}_N$}\label{repcyc}

\noindent Soient $N \geq 3$ un entier positif impair et $P\in \mathbb{N}$ défini par:
\begin{eqnarray}\label{PN}
N= 2P+1.
\end{eqnarray}
Dans la suite, on note $\omega= \exp(2i\pi/N)$ et on fixe la détermination 
$$\omega^{1/2} = \omega^{P+1}= -\exp(i\pi/N)$$
de sa racine carrée. On écrira souvent $1/2 := P + 1 \pmod N$. Pour les définitions et propriétés élémentaires concernant les algèbres de Hopf, on renvoie le lecteur à \cite{Kas} ou \cite{CP}.

\begin{defi}\label{Weyl}
L'algèbre de Weyl $\mathcal{W}_N$ est la $\mc$-algèbre associative et unitaire engendrée par les éléments $E,\ E^{-1}$ et $D$ vérifiant les relations suivantes : 
$$EE^{-1} = 1, \quad E^{-1}E = 1\ , \quad ED = \omega DE\ .$$
\end{defi}

\noindent L'algèbre de Weyl est isomorphe à une sous-algèbre de Borel quantifiée $U_{\omega}(b(2,\mathbb{C}))$ de la forme intégrale non-restreinte et simplement connexe de $U_q(sl(2,\mathbb{C}))$, o\`u le paramètre de quantification est spécialisé en la racine de l'unité $q=\omega$ \cite[\S 9.1.A]{CP}. On peut munir $\mathcal{W}_N$ de la comultiplication, de l'antipode et de la counité suivantes : 
$$\begin{array}{l}\Delta(E)  = E \otimes E\ ,\ \Delta(D) =  
E \otimes D + D \otimes 1\ ,\\ \epsilon(E)=1\ ,\ \epsilon(D)=0\ ,\ S(E)  = E^{-1}\ ,\ S(D) = -E^{-1}D\ .\end{array}$$
\noindent (Dans certains ouvrages de référence \cite[Ch. 6]{Kas}, \cite[\S 9.1]{CP}, la comultiplication $\Delta$ est d\'efinie en permutant les produits tensoriels ci-dessus). Le centre $\mathcal{Z}$ de $\mathcal{W}_N$ est engendré par $E^N,E^{-N}$ et $D^N$, et $\mathcal{W}_N$ est un $\mathcal{Z}$-module libre ayant pour base l'ensemble des mon\^omes $\{ E^sD^t,\ 0 \leq s,t < N \}$ \cite[\S 9.2.6]{CP}.

\medskip

\noindent Dans la suite, tous les espaces vectoriels consid\'er\'es sont complexes. On notera $V_{\rho}$ l'espace vectoriel support d'une représentation $\rho : \mathcal{W}_N \rightarrow {\rm End}(V_{\rho})$ de $\mathcal{W}_N$. Rappelons que $\rho$ est \emph{irréductible} si $V_{\rho}$ est un $\mathcal{W}_N$-module \emph{simple}, i.e. s'il n'admet pas de sous-modules propres autres que $\{0\}$. On dit que $\rho$ est \emph{cyclique} si $\rho(E), \rho(D) \in {\rm GL}(V_{\rho})$.

\begin{prop}\label{dimfin}
Les représentations irréductibles de $\mathcal{W}_N$ sont de dimension finie, et le centre $\mathcal{Z}$ de $\mathcal{W}_N$ y agit par multiplication scalaire. Celles qui sont irr\'eductibles et cycliques sont de dimension maximale égale à $N$.
\end{prop}

\noindent {\it Démonstration.} La seconde affirmation est une conséquence de la première. En effet, comme $\mc$ est algébriquement clos, si $\rho$ est une représentation irréductible de $\mathcal{W}_N$ de dimension finie, alors il existe un vecteur propre non nul $v \in V_\rho$ pour $E$ : $Ev = \lambda v,\ \lambda\in \mc^*$ ($\lambda$ est non nul car $E$ est inversible). La relation de commutation $ED=\omega DE$ implique donc que $E$ admet $j$ valeurs propres distinctes $\lambda\omega^i,i=0,\ldots,j$, o\`u $j\leq N$ est le plus grand entier tel que $D^jv \ne 0$. Une suite de vecteurs propres associés forme donc une base d'un sous-espace vectoriel de $V_\rho$ laissé stable par $\mathcal{W}_N$. Comme $V_\rho$ est simple, on a forcément ${\rm dim}(V_\rho)=j\leq N$, et lorsque $\rho$ est une représentation cyclique on a l'égalité.

\noindent Montrons donc la première affirmation. Soit $V$ un $\mathcal{W}_N$-module simple, et $\mathcal{Z}(V)$ l'algèbre des opérateurs sur $V$ induite par l'action du centre $\mathcal{Z}$ de $\mathcal{W}_N$. Comme $\mathcal{W}_N$ est un module de type fini sur $\mathcal{Z}$ (voir ci-dessus) et que $V$ est simple sur $\mathcal{W}_N$, on obtient que $V$ est de type fini sur $\mathcal{Z}(V)$. Pour montrer que $V$ est de dimension finie, il suffit donc de prouver que $\mathcal{Z}(V)$ est un corps. Considérons le radical de Jacobson $J$ de $\mathcal{Z}(V)$ : rappelons que c'est l'ensemble des éléments $r \in \mathcal{Z}(V)$ tels que pour tout $s \in \mathcal{Z}(V)$, $1 -rs$ est inversible. Par le lemme de Nakayama \cite[p. 103]{We}, on sait que $JV \ne V$. D'autre part, si $f \in \mathcal{Z}(V)$ est non nul, le module $f(V)$ est nécessairement égal à $V$ puisque $V$ est simple. Alors $J = \{ 0 \}$, et il existe $g \in \mathcal{Z}(V)$ tel que pour tout $v \in V$ on a $(1 -fg)(v)=0$, i.e. $f$ est inversible. Donc $\mathcal{Z}(V)$ est un corps. 

\noindent On observe que $\mathcal{Z}(V)$ agit par multiplication scalaire. En effet, comme plus haut, tout $f \in \mathcal{Z}(V)$ admet une valeur propre $\alpha$. De plus l'application $f - \alpha . id$ est $\mathcal{W}_N$-linéaire, car $f$ est central; donc son noyau est un $\mathcal{W}_N$-sous-module de $V$. Par simplicité de $V$, on en déduit que ce noyau est $V$ tout entier, c'est-à-dire que $f = \alpha . id$. Ceci achève la preuve de la proposition. \hfill $\Box$    
\medskip

\noindent La structure d'algèbre de Hopf de $\mathcal{W}_N$ permet de munir la catégorie de ses représentations d'un produit tensoriel associatif, en posant:
$$\forall \ a \in \mathcal{W}_N,\quad (\rho \otimes \mu) (a) = \rho \otimes \mu(\Delta(a))\ .$$
\noindent Soient $\delta_{i,j}$ le symbole de Kronecker mod($N$), et $X$ et $Z$ les matrices carrées de taille $N$ et de composantes $X_{ij} = \delta_{i,j+1}$ et $Z_{ij} = \omega^i \delta_{i,j}$ dans la base canonique de $\mc^N$. 

\begin{defis}\label{repstandard}
i) Une suite $\rho_1,\ldots,\rho_n$ de repr\'esentations cycliques de $\mathcal{W}_N$ est \emph{régulière} si tout produit tensoriel $\rho_i \otimes \ldots \otimes \rho_{i+j},\ 1 \leq i \leq n,\  1 \leq j \leq n -i$ est cyclique.

\smallskip

\noindent ii) Une représentation $\rho$ de $\mathcal{W}_N$ est \emph{standard} si $\rho(E) = a_\rho^{2}Z$ et $\rho(D) = a_\rho y_\rho X$, o\`u $a_\rho$, $y_\rho \in \mathbb{C}^*$.

\smallskip

\noindent iii) Si $\rho$ est standard, la repr\'esentation conjugu\'ee $\rho^*$ et la repr\'esentation inverse $\bar{\rho}$ sont définies par :
$$\begin{array}{l}
a_{\bar{\rho}} = 1/a_\rho\ ,\quad  y_{\bar{\rho}} = -y_\rho \\ 
a_{\rho^*} = (a_\rho)^*, \quad y_{\rho^*} = (y_\rho)^* \ .   
\end{array}$$
\end{defis} 

\vspace{1mm}

\noindent Pour deux représentations standards $\rho, \mu$ on v\'erifie facilement que
$$\begin{array}{lll} 
    \rho \otimes \mu(E^N) & = & \rho(E^N) \otimes \mu(E^N) = (a_\rho a_\mu)^{2N} \ id \otimes id \ , \\ \\
    \rho \otimes \mu(D^N) & = & \rho(E^N) \otimes \mu(D^N) + \rho(D^N) \otimes \mu(1) \\ \\
& = & (a_\rho a_\mu)^N \left(a_\rho^Ny_\mu^N + \frac{y_\rho^N}{a_\mu^N}\right) \ id \otimes id\ , 
\end{array}$$
\noindent o\`u $1$ est l'unité de $\mathcal{W}_N$ et $id$ est la matrice identité de dimension $N$. En particulier, pour toute représentation standard $\rho$ la représentation $\rho \otimes \bar{\rho}$ n'est pas cyclique, puisque $D$ y est nilpotent. Rappelons que deux représentations $\rho, \mu$ de $\mathcal{W}_N$ sont \emph{équivalentes}, ce que l'on notera $\rho \sim \mu$, s'il existe un isomorphisme  d'espaces vectoriels entre $V_\rho$ et $V_\mu$ qui commute avec l'action de $\mathcal{W}_N$. 

\smallskip

\noindent Notons $\mathcal{C}$ l'ensemble des repr\'esentations irr\'eductibles et cycliques de $\mathcal{W}_N$.

\begin{prop} \label{prop1a} i) Les représentations standards sont irréductibles et cycliques, et deux représentations standards $\rho$ et $\mu$ sont équivalentes si et seulement si :
$$a_\rho^{2N} = a_\mu^{2N},\quad a_\rho^Ny_\rho^N = a_\mu^Ny_\mu^N\ .$$
ii) Toute représentation $\rho \in \mathcal{C}$ est équivalente à une représentation standard. \hfill\break 
iii) Fixons une détermination de la racine $N$-i\`eme de l'unité. Si $(\rho,\mu) \in \mathcal{C}^2$ est régulière, alors $\rho \otimes \mu: \mathcal{W}_N \rightarrow {\rm End}(V_\rho) \otimes {\rm End}(V_\mu)$ se scinde en la somme directe de $N$ repr\'esentations équivalentes à la représentation standard $\rho\mu$ définie par : 
\begin{eqnarray}
a_{\rho\mu} = a_\rho a_\mu,\ y_{\rho\mu} = \left(a_\rho^Ny_\mu^N + \frac{y_\rho^N}{a_\mu^N}\right)^{1/N}. \label{multiplication}
\end{eqnarray}
\end{prop}

\noindent {\it Démonstration.}  Notons $b_0 = \{e_i,\ i=0,\ldots,\ N-1\}$ la base canonique de $\mc^N$. Une représentation standard $\rho$ est cyclique, puisque $\rho(E^N)$ et $\rho(D^N)$ sont des multiples scalaires de la matrice identité $id$, et elle est irréductible, car la relation de commutation $ED=\omega DE$ implique que l'action de $D$ sur un quelconque vecteur $e_i \in b_0$ engendre $b_0$. Si deux représentations standards déterminent la même paire $(a_\rho^{2N},\ a_\rho^Ny_\rho^N)$, elles sont nécessairement obtenues l'une de l'autre par une permutation des vecteurs de $b_0$ et sont donc équivalentes. Réciproquement, deux représentations standards équivalentes déterminent des polynômes minimaux identiques pour les endomorphismes définis par $E$ et $D$, et donc la même paire $(a_\rho^{2N},\ a_\rho^Ny_\rho^N)$.

\medskip

\noindent Montrons le second point. Soit $\rho \in \mathcal{C}$. Fixons un vecteur propre $v \in V_\rho$ de $\rho(E)$, de valeur propre $\lambda \in \mathbb{C}^*$; il en existe car $\mc$ est algébriquement clos. La relation $ED = \omega DE$ implique que les valeurs propres de $\rho(E)$ sont toutes distinctes et de la forme $\{\lambda \omega^i\},\ i=0,\ldots,\ N-1$. Soit $b = \{v_i,\ i=0,\ldots,\ N-1\},\ v_0 = v,$ une base de vecteurs propres pour $\rho(E)$ associés dans le même ordre à ces valeurs propres. Alors la base $\hat{b}$ de vecteurs obtenus par transformée de Fourier des vecteurs de $b$ est une base de vecteurs propres pour $\rho(D)$. En effet, on a :
\begin{eqnarray}
 \rho(D)\hat{v}_i = \frac{1}{N} \sum_{l=0}^{N-1} \omega^{-il} \rho(D)v_l = \frac{1}{N} \sum_{l=0}^{N-1} \omega^{-il} v_{l+1} = \omega^i\hat{v}_i\ . \label{Fourier}
\end{eqnarray} 
\noindent Posons alors $\lambda = a_\rho^2$, et soit $y_\rho\in \mc^*$ tel que $\rho(D^N) = a_\rho^N y_\rho^N id_V$. La relation (\ref{Fourier}) montre que l'isomorphisme de $V_\rho$ vers $\mc^N$ défini en envoyant $v_i$ sur $e_i$ commute avec $\rho$ et l'action de $\mathcal{W}_N$ donnée dans la définition \ref{repstandard} ii). Donc $\rho$ est équivalente à cette dernière représentation standard.

\medskip

\noindent Finalement, montrons $iii)$. Par $ii)$ on peut supposer que $\rho$ et $\mu$ sont des représentations standards. Soit $V_\alpha$ l'espace caractéristique de $\rho \otimes \mu(E) \in {\rm End}(\mc^N \otimes \mc^N)$ associé à la valeur propre $\omega^\alpha a_{\rho\mu}^2$ :
$$V_\alpha = {\rm Ker}\left(\rho \otimes \mu(E) - \omega^\alpha a_{\rho\mu}^2 \ id \otimes id\right)\ .$$
\noindent De façon classique, on a un isomorphisme d'espaces vectoriels :
$$\displaystyle \mc^N \otimes \mc^N \cong \oplus_{\alpha = 0}^{\alpha = N-1} \ V_{\alpha}\ .$$
\noindent Comme $(\rho,\mu)$  est une paire régulière et $\rho \otimes \mu(D)V_\alpha \subset V_{\alpha+1}$ à cause de la relation $ED = \omega DE$, l'action de $\rho \otimes \mu(D)$ sur cette décomposition est cyclique, et tous les $V_\alpha$ ont la même dimension $N$. Considérons la base suivante de $V_\alpha$ :
$$\{u_{i,\alpha} := v_{i+\alpha} \otimes v_{N-i},\ i = 0,\ldots,N-1\}\ .$$
\noindent Pour chaque $i,\alpha$ fixés, l'espace vectoriel engendré par l'action de $\mathcal{W}_N$ sur $u_{i,\alpha}$ est une sous-représentation de $V_\rho \otimes V_\mu$ de dimension $N$, qui est cyclique et irréductible. Par conséquent, elle est équivalente à une représentation standard, d\'efinie par les racines des polynômes minimaux de $\rho\otimes \mu(E)$ et $\rho \otimes \mu(D)$. Les formules de $\rho \otimes \mu(E^N)$ et $\rho \otimes \mu(D^N)$ montrent que ces racines sont égales respectivement aux scalaires $a_{\rho\mu}^2$ et $a_{\rho\mu}y_{\rho\mu}$ décrits en (\ref{multiplication}), modulo une racine $N$-ième de l'unité.\hfill $\Box$

\medskip

\noindent Soit $B$ le sous-goupe de Borel de $SL(2,\mc)$ des matrices triangulaires sup\'erieures. La proposition \ref{prop1a} permet de d\'efinir une application injective $\Phi: \mathcal{C}/\sim \ \rightarrow B$, o\`u Im($\Phi$) est l'ensemble des matrices non diagonales de $B$ et 
\begin{equation} \label{param}
 \Psi([\rho]) = \left( \begin{array}{ll} 
                a_\rho^N & y_\rho^N \\
                0     & a_\rho^{-N}
                \end{array} \right)\ .
\end{equation}
\noindent De plus, pour toute paire régulière $(\rho,\mu) \in \mathcal{C}^2$ on a $\Psi([\rho\mu]) = \Psi([\rho]) \cdot \Psi([\mu])$. 

\smallskip

\begin{remark} \label{coadjointe} {\rm Bien que les sous-algèbres de Borel de $sl(2,\mc)$ ne soient pas semi-simples, l'existence de $\Psi$ est aussi une conséquence de la théorie de l'action quantique coadjointe de De Concini, Kac et Procesi \cite{DC-P},\ \cite[\S 9.2,11.1]{CP}. Pour toute algèbre de Lie complexe, semi-simple et de dimension finie $\mathfrak{g}$ et toute racine de l'unité $\epsilon$ d'ordre impair $l$ assez grand, ces derniers construisent un paramétrage $\Theta$ des représentations irréductibles de $U_{\epsilon}(\mathfrak{g})$, à valeurs dans un groupe de Lie. Comme pour $\Psi$, le paramétrage $\Theta$ se fait par l'intermédiaire de l'ensemble $Spec(\mathcal{Z}_{\epsilon}(\mathfrak{g}))$ des caractères centraux de ces représentations, qui sont les homomorphismes $\chi : \mathcal{Z}_{\epsilon}(\mathfrak{g}) \rightarrow \mc$ définis par l'action du centre $\mathcal{Z}_{\epsilon}(\mathfrak{g})$ de $U_{\epsilon}(\mathfrak{g})$. Ainsi, les représentations irréductibles et cycliques de $U_{\epsilon}(sl(2,\mc))$ peuvent toutes être obtenues à partir de produits tensoriels tordus (``pull-back'') de celles de $\mathcal{W}_N$ \cite{DJMM,CP}.}
\end{remark}

\section{Le $q$-dilogarithme}\label{asy}

\noindent Soit $\mc(q)$ le corps des fractions rationnelles en une indéterminée $q$. Le $q$-dilogarithme est la série formelle $(x;q)_{\infty} \in \mc(q)[[x]]$ définie par :
\begin{equation}
(x;q)_{\infty} = \prod_{n=0}^{\infty} (1-xq^n) = \sum_{n=0}^{\infty} \frac{q^{\frac{n(n-1)}{2}}(-x)^n}{(q)_n}\ ,\quad \mathrm{avec} \ (q)_n = \frac{(q;q)_{\infty}}{(q^{n+1};q)_{\infty}}\ \label{dieq0}.
\end{equation}
\noindent On s'int\'eresse au $q$-dilogarithme d'abord comme op\'erateur sur l'algèbre de Weyl. Soit $\mathcal{A}_q$ la $\mc(q)$-algèbre engendrée par $u,\ v$ et $w$ tels que :
\begin{equation} \label{extcent}
uv -qvu = w\ ,\quad [u,w] = [v,w] = 0\ .
\end{equation}
\noindent L'algèbre $\mathcal{A}_q$ est une extension centrale de l'algèbre de Weyl $\mathcal{W}_q$, que l'on obtient en posant $w=0$. Considérons la complétion $q$-adique $\widehat{\mathcal{A}}_q$ (resp. $\widehat{\mathcal{W}}_q$), qui est l'espace vectoriel des séries formelles en $u,\ v$ et $w$ (resp. $u$ et $v$) et à coefficients dans $\mc(q)$. L'ensemble $\{u^kv^lw^m \vert\ k,l,m \in \mathbb{N}\}$ (resp. $\{u^kv^l \vert \ k,l \in \mathbb{N}\}$) est une base topologique du $\mc(q)$-module $\widehat{\mathcal{A}}_q$ (resp. $\widehat{\mathcal{W}}_q$).

\begin{lem} \label{qdilogen} On a la relation suivante dans $\widehat{\mathcal{A}}_q$ :
\begin{equation} \label{qdilocentral}
(u;q)_{\infty}\ \left( \frac{[u,v]}{(1-q)};q\right)_{\infty}\ (v;q)_{\infty} = (v;q)_{\infty}\ (u;q)_{\infty}\ ;\end{equation}
\noindent cette relation devient \emph{l'équation du $q$-dilogarithme} dans $\widehat{\mathcal{W}}_q$ :
\begin{equation} \label{qdiloeq}
(u;q)_{\infty}\ \left(-vu;q\right)_{\infty}\ (v;q)_{\infty} = (v;q)_{\infty}\ (u;q)_{\infty}\ .
\end{equation}
\end{lem}

\medskip

\noindent  La relation (\ref{qdiloeq}) est due à L. Faddeev et A. Volkov \cite{Kir}, et la relation (\ref{qdilocentral}) a été annoncée par Kashaev dans \cite{K3} comme due à Volkov. 

\medskip

\noindent {\it Démonstration}. On va montrer que les deux membres de la relation (\ref{qdilocentral}) vérifient la même équation linéaire algébrique aux $q$-différences, lorsqu'on les considère comme fonctions de $u$. On en d\'eduit le r\'esultat par unicité des solutions formelles de telles équations, et le fait que les deux membres de (\ref{qdilocentral}) sont \'egaux lorsque $u=0$. On a 
$$\begin{array}{lll}
(v;q)_{\infty}\ (u;q)_{\infty}\ (qu;q)_{\infty}^{-1}\ (v;q)_{\infty}^{-1} & = & (v;q)_{\infty}\ (1 - u)\ (v;q)_{\infty}^{-1} \\ \\
                         & = & 1 - (v;q)_{\infty}\ u \ (v;q)_{\infty}^{-1}\ .
\end{array}$$
\noindent Par réccurence, on v\'erifie avec (\ref{extcent}) que 
$$v^n u = q^{-n}uv^n - (q^{-1} + \ldots + q^{-n}) \ wv^{n-1} = q^{-n}uv^n - q^{-n}\ \frac{1-q^n}{1-q}\ wv^{n-1}$$
\noindent pour tout entier $n \geq 1$. Supposons pour simplifier que $u$ et $v$ sont inversibles, d'inverses $u^{-1}$ et $v^{-1}$ (le résultat ne dépend pas de cette hypothèse, car $u^{-1}$ et $v^{-1}$ n'apparaissent pas dans l'équation aux $q$-différences vérifiée par les deux membres de (\ref{qdilocentral})). Alors
$$\begin{array}{lll}
(v;q)_{\infty}\ u & = & \sum_{n=0}^{\infty} \frac{q^{\frac{n(n-1)}{2}}(-v)^n}{(q)_n}\ u\\
& & \\                  
& = & \sum_{n=0}^{\infty} \frac{q^{\frac{n(n-1)}{2}}(-1)^n}{(q)_n} \ (q^{-n}uv^n - q^{-n}\ \frac{1-q^n}{1-q}\ wv^{n-1}) \\
& & \\                  
& = & u\ (q^{-1}v;q)_{\infty} + \frac{wv^{-1}}{1-q}\ \left( (v;q)_{\infty} - (q^{-1}v;q)_{\infty} \right) \ .
\end{array}$$

\noindent On en déduit que 

$$\begin{array}{lll} 
(v;q)_{\infty}\ (u;q)_{\infty}\ (qu;q)_{\infty}^{-1}\ (v;q)_{\infty}^{-1} & = & 1 - \bigl( u(1-q^{-1}v) + \frac{wv^{-1}}{1-q}\ ( 1- (1-q^{-1}v)) \bigr) \\
& & \\                  
& = & 1 - u + q^{-1}uv - \frac{q^{-1}w}{1-q} \\
& & \\
& = & 1 - u + vu - \frac{w}{1-q} \\
& & \\
& = & 1 - u -  \frac{[u,v]}{1-q}\ .
\end{array}$$
\noindent D'autre part, on a 
$$\begin{array}{l}
(u;q)_{\infty}\ \left(\frac{[u,v]}{1-q};q\right)_{\infty}\ (v;q)_{\infty}\ (v;q)_{\infty}^{-1}\ \left(\frac{q[u,v]}{1-q};q\right)_{\infty}^{-1}\ (qu;q)_{\infty}^{-1} \hspace{3cm} \\
 \\
\hspace{7cm} = (u;q)_{\infty}\ \left(1 - \frac{[u,v]}{1-q}\right)\ (qu;q)_{\infty}^{-1} \\
\\  
\hspace{7cm} =  1 - u - (u;q)_{\infty}\ \frac{[u,v]}{1-q} \ (qu;q)_{\infty}^{-1}\ .
\end{array}$$
\noindent Comme plus haut, on vérifie par réccurence que
$$u^n v = q^{n}vu^n + \frac{1-q^n}{1-q}\ wu^{n-1}$$
\noindent pour tout entier $n \geq 1$. Alors 
$$(u;q)_{\infty}\ v = v\ (qu;q)_{\infty} + \frac{wu^{-1}}{1-q}\ \left( (u;q)_{\infty} - (qu;q)_{\infty}\right)\ , $$
\noindent d'où l'on obtient
$$(u;q)_{\infty}\ \frac{[u,v]}{1-q}\ (qu;q)_{\infty}^{-1} = \frac{[u,v]}{1-q}\ .$$
\noindent Par conséquent, on a 
$$1 - u - (u;q)_{\infty}\ \frac{[u,v]}{1-q} \ (qu;q)_{\infty}^{-1}  =  1 - u -  \frac{[u,v]}{1-q}\ .$$
\noindent Donc les deux membres de la relation (\ref{qdilocentral}) vérifient la même équation aux $q$-différences en $u$ ($v$ étant considéré comme paramètre) : 
$$f_v(u) - \left( 1 - u -  \frac{[u,v]}{1-q} \right) \ f_v(qu) = 0\ ,$$
\noindent ce qui achève la preuve.\hfill $\Box$

\medskip

\noindent Le $q$-dilogarithme (\ref{dieq0}) est solution de l'équation linéaire aux $q$-différences 
$$(1-x)f(qx) - f(x) = 0\ .$$
Lorsque $\vert q \vert < 1$, $(x;q)_{\infty}$ converge normalement sur tout compact de $\mc$ et définit une fonction entière $E_q(x)$ que l'on peut écrire sous la forme du produit infini $\textstyle (x ;q)_{\infty} = \prod_{n=0}^{\infty} (1-xq^n)$ \cite[Prop. 5.1]{R}. Les zéros de $E_q(x)$ sont simples et parcourent $\{q^{-n} ,\ n \geq 0\}$. Lorsque $\vert q \vert > 1$, le rayon de convergence de $(x;q)_{\infty}$ est égal à $\vert q \vert$, et sa somme définit dans le disque ouvert de convergence une fonction holomorphe $e_q(x)$. De plus, $e_q(x)$ admet un prolongement méromorphe à $\mc$ que l'on peut écrire $e_q(x) = E_{q^{-1}}(xq^{-1})^{-1}$. La fonction $e_q(x)$ n'a pas de z\'eros; ses pôles sont simples et parcourent $\{q^{n}\ ,\ n \geq 1\}$. On va montrer maintenant que $(x;\ q)_{\infty}$ a des singularités essentielles en les racines de l'unité (Proposition \ref{devqdilo} ci-dessous). Dans la suite on n'utilisera pas ce r\'esultat. 

\medskip

\begin{lem} {\rm \cite{Kir}}
Soit ${\rm Li}_2(x;q) \in \mc(q)[[x]]$ la série formelle définie par : 
$${\rm Li}_2(x;q) = \log((x;q)_{\infty})\ ,$$ 
\noindent o\`u $\log$ signifie son développement en série au voisinage de $1$. On a :
\begin{equation} \label{li2}
{\rm Li}_2(x;q)=-\sum_{n=1}^{\infty} \frac{x^n}{n(1-q^n)}\ .
\end{equation}
\end{lem}

\noindent {\it Démonstration}. La relation $(x;q)_{\infty}=(1-x)(qx;q)_{\infty}$ implique
$${\rm Li}_2(x;q) = {\rm Li}_2(qx;q) + \log(1-x)\ .$$
\noindent En écrivant $\textstyle {\rm Li}_2(x;q)= \sum_{n=1}^{\infty} a_nx^n$ avec $a_n \in \mathbb{Q}[q,\ q^{-1}]$, on déduit que $a_n=q^na_n-\frac{1}{n}$, d'o\`u le résultat. \hfill $\Box$ 

\medskip

\noindent Soit ${\rm Li}_2$ le dilogarithme d'Euler, défini pour $x \in [0,1]$ par 
$${\rm Li}_2(x) := \sum_{n=1}^{\infty} \frac{x^n}{n^2} = -\int_0^x \frac{\log(1-t)}{t}\ dt\ .$$
\noindent La représentation intégrale à droite permet de prolonger ${\rm Li}_2$ analytiquement à $\mc \setminus ]1;\infty[$. Avec (\ref{li2}) on voit que $(q-1){\rm Li}_2(x;q)$ converge formellement pour $q=e^h$ et $h \rightarrow 0,$ vers le développement en série de ${\rm Li}_2(x)$ dans le disque unité ouvert de $\mc$. Précisons ce résultat. Rappelons que si $f$ est une fonction à valeurs complexe de classe $\mathcal{C}^{2n}$ sur un intervalle $[M,N]$, où $M$ et $N$ sont des entiers distincts, on a la formule d'Euler-Mac Laurin
$$\begin{array}{lll}
\sum_{i=M}^N f(i) & = & \int_M^N f(x)dx + \frac{1}{2}(f(M) + f(N)) +\\
& = & \sum_{k=1}^{n} \frac{B_{2k}}{(2k)!} ( f^{(2k-1)}(N) - f^{(2k-1)}(M) )- \int_{M}^N \frac{\overline{B}_{2n}(t)}{(2n)!} f^{(2n)}(t)dt\ ,
\end{array}$$
\noindent où $B_k(t)$ est le $k$-ième polynôme de Bernoulli, défini par
$$ \frac{ze^{tz}}{e^z - 1} = \sum_{k \geq 0}\frac{B_k(t)}{k!} z^k\ ,$$
\noindent et $B_k := B_k(0)$ et $\overline{B}_n(t)=B_n(\{ t \})$, où $\{ t \} = t$ (mod $\mz$). Posons $q=e^{-\epsilon/N^2}\zeta$, où $N \geq 2$ est un entier, $\zeta$ est une racine primitive $N$-ième de l'unité, et $\epsilon \in \mathbb{C}^{\ast}$ vérifie ${\rm Re}(\epsilon) > 0$. Consid\'erons 
$$S(x)=S(x,\epsilon)=(1-x)^{1/2} \exp(-\frac{{\rm Li}_2(x)}{\epsilon})$$
$$\Phi(x)=\Phi(x,\zeta,N))=(1-x^N)^{(N-1)/2N} \prod_{k=1}^{N-1} (1-\zeta^kx)^{-k/N}\ ,$$
\noindent où les racines doivent être comprises comme leurs développements en série au voisinage de $x=0$. La s\'erie formelle $\Phi(x)$ définit une fonction analytique sur le disque unité ouvert $D = \{ x \in \mc \vert \ \vert x \vert < 1 \}$, qui pour tout $\theta \in \mr$ se prolonge analytiquement à 
$$\mathfrak{D}_{\theta} = \mathbb{C} \setminus \{ \vert x \vert \geq 1,\ \arg(x)= \theta + 2\pi k/N,\ k=0, \ldots ,\ N-1 \}\ .$$ 
\begin{prop} \label{devqdilo} {\rm \cite{BR}}
Soit $q=\exp(-\epsilon/N^2)\zeta$, o\`u $\zeta$ est une racine primitive $N$-ième de l'unité et ${\rm Re}(\epsilon) > 0$. Pour tout nombre complexe $x$ tel que $\vert x \vert < 1$, le $q$-dilogarithme $(x;q)_{\infty}$ a pour $\epsilon \rightarrow 0$ le développement asymptotique
\begin{equation} \label{dev}
(x;q)_{\infty}= (1-x^N)^{(1-N)/2N} S(x^N) \Phi(x) (1 + \mathcal{O}(\epsilon))\ .
\end{equation}
\end{prop} 

\noindent {\it Démonstration}. On remarque d'abord que $(x;q)_{\infty}$ se décompose en $N$ facteurs dont les $q$-paramètres convergent vers $1$ :
\begin{equation} \label{qdilsplit1}
(x;q)_{\infty}= \prod_{k=0}^{N-1} \prod_{j=0}^{\infty}(1-q^{Nj}q^kx) \stackrel{\mathrm{def}}{=} \ \prod_{k=0}^{N-1} (xq^k;\ q^N)\ .
\end{equation}
\noindent D'autre part, on a 
$$\begin{array}{lll}
{\rm Li}_2(e^{\alpha \epsilon} x) & = & - \int_0^{xe^{\alpha \epsilon}} \frac{\log(1 - t)}{t} \ dt \\
 & &  \\                                
& = & - \int_0^x \frac{\log(1 - e^{\alpha \epsilon}u)}{u} \ du \ .
\end{array}$$
\noindent Supposons que $\alpha \in \mathbb{R},\ \alpha < 0$, et que $x \in D = \{ x \in \mc \ \vert \ \vert x \vert < 1\}$. Dans la dernière \'egalit\'e, intégrons le long d'un chemin $\gamma$ de $\left( \mc \ \setminus \ ]1;\infty[ \right) \ \cap \ D$ allant de $0$ à $x$. Alors $u \in D$ et 
$$\begin{array}{lll}
\log(1 - e^{\alpha \epsilon}u) & = & \log(1 - (1 + \alpha \epsilon + \mathcal{O}(\epsilon^2))u) \\
                               & = & \log(1 - u - \epsilon (\alpha u + \mathcal{O}(\epsilon)u))\ .
\end{array}$$
\noindent Posons $h(\epsilon,u) = \epsilon (\alpha u + \mathcal{O}(\epsilon)u)$ . On a
$$\log(1 - e^{\alpha \epsilon}u) = \log(1 - u) - \frac{h(\epsilon,u)}{1-u} + \mathcal{O}(h(\epsilon,u)^2)\ .$$
\noindent Comme $h(\epsilon,u) = \mathcal{O}(\epsilon)$, en utilisant la définition explicite de $h$ et en intégrant le long du chemin $\gamma$ (chaque intégrale étant définie) on trouve pour tous $x \in D$ et $\alpha < 0$:
$$\begin{array}{lll}
{\rm Li}_2(e^{\alpha \epsilon} x) & = & {\rm Li}_2(x) + \alpha \epsilon \int_0^x \frac{du}{1-u} + \mathcal{O}(\epsilon^2) \\ \\
& = & {\rm Li}_2(x) - \alpha \epsilon \log(1-x) + \mathcal{O}(\epsilon^2)\ ,
\end{array}$$
\begin{equation} \label{temp1}
\exp(\frac{-{\rm Li}_2(e^{\alpha \epsilon} x)}{\epsilon})=
(1-x)^{\alpha} \exp(\frac{-{\rm Li}_2(x)}{\epsilon})(1 + \mathcal{O}(\epsilon))\ .
\end{equation}
\noindent Appliquons la formule d'Euler-Mac Laurin au logarithme de $(x;e^{-\epsilon})_{\infty}$, avec $M=0$ et $N \rightarrow \infty$. Comme ${\rm Re}(\epsilon) > 0$, pour tous $x \in D$ et $t \in [0;\infty[$, on a $\vert xe^{-\epsilon t} \vert < 1$. Donc
$$\log\left( (x;e^{-\epsilon})_{\infty}\right) = \sum_{n=0}^{\infty} \log(1-xe^{-\epsilon n}) = \int_0^{\infty} \log(1-xe^{-\epsilon t}) \ dt + \frac{1}{2} \log(1 - x) + \mathcal{O}(\epsilon)\ .$$
\noindent On remarque que
$$\frac{d\log(1-xe^{-\epsilon t})}{dt} = \epsilon \ \frac{x e^{-\epsilon t}}{1 - xe^{-\epsilon t}} = \mathcal{O}(\epsilon)$$
\noindent et de même pour les dérivées successives de $\log(1-xe^{-\epsilon t})$. Comme 
$$\int_0^{\infty} \log(1-xe^{-\epsilon t}) \ dt = -\frac{1}{\epsilon} \int_x^0 \frac{\log(1-u)}{u} \ du\ ,$$
\noindent on a le développement asymptotique
$$\forall \ x \in D,\quad (x;e^{-\epsilon})_{\infty} = S(x)(1 + \mathcal{O}(\epsilon))\ ,\quad \epsilon \rightarrow 0\ .$$
\noindent On obtient donc pour $\epsilon \rightarrow 0$ et $x \in D$:
$$\begin{array}{lll}
(xq^k;q^N)_{\infty} & = & S(e^{-k\epsilon / N^2}\zeta^k x) (1 + \mathcal{O}(\epsilon) )\ ,\\ \\
& = & (1-x\zeta^k)^{1/2-k/N} \exp(\frac{-N{\rm Li}_2(x\zeta^k)}{\epsilon})
(1 + \mathcal{O}(\epsilon))\ .
\end{array}$$
\noindent Dans la seconde \'egalit\'e, on utilise (\ref{temp1}) avec $\alpha = -k/N$. Alors pour $\epsilon \rightarrow 0$ et $x \in D$, (\ref{qdilsplit1}) donne
$$(x;q)_{\infty}= (1-x^N)^{1/2} \prod_{k=1}^{N-1} (1-\zeta^kx)^{-k/N} \exp(-N
\sum_{k=0}^{N-1} \frac{{\rm Li}_2(\zeta^kx)}{\epsilon})(1 + \mathcal{O}(\epsilon))\ .$$
\noindent Lorsque $\vert x \vert < 1$, on a la c\'elèbre relation
$$N \sum_{k=0}^{N-1} {\rm Li}_2(\zeta^kx) = N \sum_{k=0}^{N-1} \sum_{n \geq 1} \frac{x^n \zeta^{kn}}{n^2} = N \sum_{n \geq 1} \frac{x^n}{n^2} \sum_{k=0}^{N-1} \zeta^{kn} = \sum_{j \geq 1} \frac{ x^{Nj}}{j^2} = {\rm Li}_2(x^N)\ ,$$
\noindent où l'on remarque que $\textstyle \sum_{k=0}^{N-1} \zeta^{kn}$ est non nul si et seulement si $N$ divise $n$. Ceci achève la preuve de la proposition. \hfill $\Box$

\medskip

\noindent Considérons la spécialisation en $q=\exp(-\epsilon/N^2)\zeta$ de $\widehat{\mathcal{W}}_q$. De façon classique, on peut d\'efinir $\widehat{\mathcal{W}}_q$ comme quantification par déformation de $\widehat{\mathcal{W}}_\zeta$ \cite[\S 6]{CP}: le crochet de Poisson du centre $\mathcal{Z}(\widehat{\mathcal{W}}_\zeta)$ de $\widehat{\mathcal{W}}_\zeta$ et son action de Poisson sur $\widehat{\mathcal{W}}_\zeta$ sont définis par
$$\{\bar{a},\bar{b}\} = \lim_{\epsilon \rightarrow 0}\ \frac{ab -ba}{\epsilon}\ ,$$
\noindent o\`u $a,\ b \in \widehat{\mathcal{W}}_q$, $\bar{a} \equiv a$ mod($\epsilon$) et $\bar{b} \equiv b$ mod($\epsilon$), et $\bar{a}$ ou $\bar{b}$ est dans $\mathcal{Z}(\widehat{\mathcal{W}}_\zeta)$. La proposition \ref{devqdilo} implique que $\Phi$ est la composante non ab\'elienne de l'action de $(x;q)_{\infty}$ pour $q \rightarrow \zeta$. On a l'analogue cyclique suivant de la d\'efinition de $(x;q)_{\infty}$ comme produit infini. 

\begin{lem}\label{gdilocyc}
Posons $r (x) = (1 - x^N)^{1/N}$. Pour tout $x \in \mathfrak{D}_{0},\ x \ne \omega^j$ avec $j=1,\ldots,\ N$, et $k \in \{0,\ldots,\ N-1 \}$, on a l'identité :
$$\Phi(x\zeta^{k}) = \Phi(x)\prod_{j=0}^{k - 1} \frac{ \zeta r(x)}{1-x\zeta^j}\ .$$
\end{lem}

\noindent {\it Démonstration du lemme.} Dans la d\'efinition de $\Phi(x)$, le facteur $(1-x^N)^{(N-1)/2N}$ ne pose pas de problème pour prouver l'identit\'e. Posons donc
$$\Phi'(x)=\prod_{j=1}^{N-1} (1-x\zeta^j)^{-j/N}\ .$$
\noindent On a
$$\begin{array}{lll}
(\Phi'(x\zeta^k))^N & = & \prod_{j=1}^{N-1}(1-x\zeta^{k}\zeta^j)^{-j} \\ \\
& = & \prod_{j=1}^{N-1}(1-x\zeta^{k+j})^{k}\ \prod_{j=N-k+1}^{N-1}(1-x\zeta^{k+j})^{-(k+j)}  \\ & & \hspace{5cm} \times \  \prod_{j=1}^{N-k}(1-x\zeta^{k+j})^{-(k+j)}
\end{array}$$
$$\begin{array}{lll}
 & = & \left( \frac{1-x^N}{1-x\zeta^{k}}\right)^{k} \ \left(\ \prod_{l=1}^{k-1}(1-x\zeta^{l})^{-l}\  \prod_{l=1}^{k-1}(1-x\zeta^{l})^{-N} \right)\ \prod_{l=k+1}^{N}(1-x\zeta^{l})^{-l} \\ \\
 & = & (\Phi'(x))^N \ (1-x)^{-N}\ (1-x^N)^{k} \ \prod_{l=1}^{k-1}(1-x\zeta^{l})^{-N}\\ \\
 & = & (\Phi'(x))^N \ \prod_{j=0}^{k - 1} \frac{(1-x^N)}{(1-x\zeta^j)^N} \ .
\end{array}$$

\noindent D'autre part, lorsqu'on passe dans le sens positif à travers la coupure de $\mathfrak{D}_{0}$ partant du point $\exp(2i \pi k/N) = \zeta^{k}$, on change de détermination de la racine $N$-ième et la valeur de $\Phi'$ est multipliée par $\zeta^{k}$. D'o\`u le r\'esultat.  \hfill $\Box$

\medskip

\noindent  Dans \cite{BR}, V. Bazhanov et N. Reshetikhin ont montr\'e que lorsque $q \rightarrow \zeta$, le terme dominant en $1/\epsilon$ de l'équation du $q$-dilogarithme (\ref{qdiloeq}) se scinde en deux \'equations fonctionnelles plus \'el\'ementaires.  La première est l'identité de Rogers \cite{Kir}, obtenue sur $\mathcal{Z}(\widehat{\mathcal{W}}_\zeta)$ à l'aide de $S(x)$. La seconde est d\'ecrite dans le th\'eorème suivant; on en prouvera une version \'equivalente au \S \ref{descdc}.

\noindent Soient $\alpha$ et $\beta$ deux nombres complexes non nuls, et $H_\zeta(\alpha,\beta)$ le quotient de $\widehat{\mathcal{W}}_\zeta$ par l'idéal engendré par $(u^N -\alpha^N)$ et $(v^N - \beta^N)$; on a ${\rm dim}_{\mc}(H_\zeta(\alpha,\beta)) = N^2$.

\begin{teo1} \label{eqdiloasy} {\rm \cite[Th. 3.3]{BR}} La limite semi-classique pour $\epsilon \rightarrow 0$ de l'équation du $q$-dilogarithme sur $\widehat{\mathcal{W}}_q$, o\`u $q=\exp(-\epsilon/N^2)\zeta$, ${\rm Re}(\epsilon) > 0$, et $\zeta$ est une racine primitive $N$-ième de l'unité, induit \emph{l'\'equation du dilogarithme cyclique} sur $H_\zeta(\alpha,\beta)$ :
$$\Phi(v)\ \Phi(u) = \Phi\left(\frac{u}{(1-\beta^N)^\frac{1}{N}}\right) \ \Phi\left(\frac{-vu}{(1-\alpha^N-\beta^N)^\frac{1}{N}}\right)\ \Phi\left(\frac{v}{(1-\alpha^N)^\frac{1}{N}}\right)\ .$$ \hfill $\Box$
\end{teo1}

\section{Le double de Heisenberg de $\mathcal{W}_N$}\label{double}

\noindent Cette section est organis\'ee comme suit. Au \S \ref{can}, on construit le double de Heisenberg $\mathcal{H}_h(b(2,\mc))$ de la QUE-algèbre $U_h(b(2,\mc))$, et on d\'ecrit son \'el\'ement canonique $R_h$ à l'aide du $q$-dilogarithme (Proposition \ref{defdoubleh}). Au \S \ref{rationnalisation}, on considère une \emph{forme intégrale} $\mathcal{H}_q$ de $\mathcal{H}_h(b(2,\mc))$, et lorsque $q$ est une racine primitive $N$-ième de l'unit\'e, on d\'ecrit l'op\'erateur $R_\omega$ induit par $R_h$ sur (une extension convenable de) $\mathcal{H}_q$. On d\'ecrit explicitement la partie ``$q$-dilogarithmique'' de $R_\omega$ au \S \ref{descdc}.

\subsection{$\mathcal{H}_h(b(2,\mc))$ et son élément canonique} \label{can}

\noindent On commence par rappeler les propri\'et\'es \'el\'ementaires des doubles de Heisenberg \cite{BS,K3,STS}.

\medskip

\noindent Soit $A = (1,\epsilon,m,\Delta,S)$ une algèbre de Hopf unitaire sur un anneau $k$, o\`u l'on note respectivement $1$, $\epsilon$, $m$, $\Delta$ et $S$ l'unité, la counité, la multiplication, la comultiplication et l'antipode de $A$. Si $A$ est de dimension infinie, on suppose que $k=\mc[[h]]$ est l'anneau des séries formelles sur $\mc$ d'indéterminée $h$, et que $A = U_h(\mathfrak{g})$ est l'algèbre enveloppante quantique (QUE) d'une algèbre de Lie complexe $\mathfrak{g}$ de dimension finie \cite[\S 16-17]{Kas},\ \cite[\S 6-8]{CP}.  Tous les produits tensoriels écrits ci-dessous sont définis sur $k$, et lorsque $A=U_h(\mathfrak{g})$, ils sont implicitement \emph{complétés} dans la topologie $h$-adique.

\medskip

\noindent De façon classique, si $A$ est de dimension finie on peut munir le $k$-module dual $A^* = {\rm Hom}(A,k)$ d'une structure de bigèbre \emph{duale} $A^* = (\epsilon^*,1^*,\Delta^*,m^*,S^*)$ avec
$$\begin{array}{l}
\langle x,m(a \otimes b)\rangle = \langle m^*(x),a \otimes b\rangle\ ,\quad \langle x \otimes y,\Delta(a)\rangle = \langle \Delta^*(x \otimes y),a \rangle\\ \\
\langle x,1\rangle = 1^*(x)\ ,\quad \langle \epsilon^*,a\rangle = \epsilon(a)\ ,\quad \langle x,S(a)\rangle = \langle S^*(x),a\rangle\ ,
\end{array}$$
\noindent o\`u $\langle \ ,\ \rangle : A^* \otimes A \rightarrow k$ désigne une application $k$-bilinéaire non dégénérée, et $a$, $b \in A$ et $x$, $y \in A^*$. Lorsque $A=U_h(\mathfrak{g})$ il existe une notion de QUE-\emph{dual} que l'on notera aussi $A^* = U_h^*(\mathfrak{g})$ \cite[\S 6.3.D-8.3]{CP}. Le QUE-dual $U_h^*(\mathfrak{g})$ est une QUE-algèbre isomorphe à $U_h(\mathfrak{g})$ en tant que $\mc[[h]]$-module, et sa structure d'algèbre de Hopf est duale de celle de $U_h(\mathfrak{g})$ dans le même sens que ci-dessus. (La limite classique de $U_h^*(\mathfrak{g})$ au sens des quantifications d'algèbres de co-Poisson Hopf est une bigèbre de Lie duale de $\mathfrak{g}$). Nous décrirons plus loin en détail le QUE-dual d'une algèbre de Borel quantifiée $U_h(b(2,\mc))$. 

\medskip

\noindent Fixons une fois pour toute un crochet de dualité $\langle \ ,\ \rangle$ entre $A$ et $A^*$. Soient $\{e_\alpha\}_{\alpha}$ et $\{e^{\beta}\}_{\beta}$ des bases (topologiques) duales de $A$ et $A^*$ : $\langle e^\beta,e_\alpha \rangle = \delta_{\alpha}^{\beta}$ ($\delta$ désigne le symbole de Kronecker). Pour tout $a \in A$, notons ${\rm ev}_a : A^* \rightarrow k$ le morphisme d'évaluation en $a$ associé au crochet $\langle \ ,\ \rangle$ : 
$$\forall \ x \in A^*,\quad {\rm ev}_a(x) = \langle x,a \rangle\ .$$
\noindent Soient $\pi_A: A \rightarrow {\rm End}_k(A^*)$ l'homomorphisme d\'efini par
$$\forall \ a \in A, \quad \pi_A(a) = (id \otimes {\rm ev}_a) \ m^* \ ,$$
\noindent et $\pi_{A^*}: A^* \rightarrow {\rm End}_k(A^*)$ l'action par multiplication à gauche.
\begin{defi} \label{doubleA}
Le \emph{double de Heisenberg} $\mathcal{H}(A)$ est la $k$-sous-algèbre de ${\rm End}_k(A^*)$ (topologiquement) engendrée par l'image des homomorphismes $\pi_A$ et $\pi_{A^*}$.
\end{defi}
\noindent En utilisant $\langle \ ,\ \rangle$, les constantes de structures (i.e. les coefficients dans les expressions suivantes) de $A^*$ et $A$ sont
\begin{eqnarray} \label{relations}
m(e_{\alpha} \otimes e_{\beta}) = \sum_{\gamma} m_{\alpha,\beta}^{\gamma}\ e_{\gamma}\ , & \ \Delta(e_{\alpha}) = \sum_{\beta , \gamma} \Delta_\alpha^{\beta,\gamma}\ e_\beta \otimes e_\gamma\ , \nonumber \\
\Delta^*(e^{\alpha} \otimes e^{\beta}) = \sum_{\gamma} \Delta_\gamma^{\alpha,\beta}\ e^{\gamma}\ , & \ m^*(e^{\alpha}) = \sum_{\beta , \gamma} m_{\beta,\gamma}^{\alpha}\ e^\beta \otimes e^{\gamma}\ .\label{relations}
\end{eqnarray}
\begin{lem}
La $k$-algèbre associative (topologiquement) engendrée par  les éléments $\{e_{\alpha}\}_\alpha$ et $\{ e^{\beta} \}_\beta$ vérifiant les relations
\begin{eqnarray} 
e_{\alpha}e_{\beta} = \sum_{\gamma} m_{\alpha,\beta}^{\gamma}\ e_{\gamma}\ ,\quad e^{\alpha}e^{\beta} = \sum_{\gamma} \Delta_\gamma^{\alpha,\beta}\ e^{\gamma}\ ,\nonumber \\ 
e_{\alpha}e^{\beta} = \sum_{\rho,\gamma,\sigma} m_{\rho,\gamma}^\beta \ \Delta_\alpha^{\gamma,\sigma} e^{\rho}e_{\sigma} \ \hspace{1cm}\label{relations2}
\end{eqnarray}
est isomorphe à $\mathcal{H}(A)$.
\end{lem}
\noindent {\it D\'emonstration.} Les deux premières relations sont \'evidentes, elles sont obtenues en omettant de noter $m$, $\Delta^*$. Les actions de $A$ et $A^*$ sur le $k$-module $A^*$, associées respectivement aux homomorphismes $\pi_A$ et $\pi_{A^*}$, sont décrites dans les bases $\{e_{\alpha}\}_{\alpha}$ et $\{e^{\alpha}\}_{\alpha}$ par
\begin{equation} \label{repadj}
e_{\alpha}(e^\beta) = (id \otimes {\rm ev}_{e_{\alpha}})\ m^*(e^\beta)\ ,\quad e^\alpha (e^\beta) = \Delta^*(e^\alpha \otimes e^\beta)\ ,
\end{equation}
\noindent ce qui s'écrit matriciellement
$$(e_{\alpha})_{i,j} = m_{i,\alpha}^j\ , \quad (e^{\alpha})_{i,j} = \Delta_i^{\alpha,j}\ .$$
\noindent Maintenant, remarquons que $\textstyle (e_{\alpha}e^{\beta})_{i,j} = \sum_k m_{i,\alpha}^k \ \Delta_k^{\beta,j}$. Par définition, la comultiplication $\Delta$ est un morphisme d'algèbre : $\Delta \circ m = (m \otimes m) \circ p_{23} \circ (\Delta \circ \Delta)$ (o\`u $p_{23}$ est la permutation des deuxièmes et troisièmes facteurs tensoriels). En termes de constantes de structures, cela donne 
$$(e_{\alpha}e^{\beta})_{i,j} = \sum_k \Delta_k^{\beta,j} \ m_{i,\alpha}^k= \sum_{\rho,\gamma,\sigma,k} m_{\rho,\gamma}^\beta \ m_{k,\sigma}^j \ \Delta_i^{\rho,k}\ \Delta_\alpha^{\gamma,\sigma}  = \sum_{\rho,\gamma,\sigma} m_{\rho,\gamma}^\beta \ \Delta_\alpha^{\gamma,\sigma} \ (e^{\rho}e_{\sigma})_{i,j}\ .$$
On en d\'eduit la troisième relation. \hfill $\Box$

\medskip 

\noindent En utilisant cette présentation de $\mathcal{H}(A)$, on v\'erifie imm\'ediatement que son \emph{élément canonique} 
$$R = \sum_{\gamma} e_\gamma \otimes e^\gamma \quad \in \mathcal{H}(A) \otimes \mathcal{H}(A)$$
\noindent satisfait pour tout $\alpha$ les égalités :
\begin{eqnarray} \label{pentagone}
(1 \otimes e_{\alpha}) \ R = R \ \Delta(e_\alpha) \ ,\quad R \ (e^{\alpha} \otimes 1) = m^*(e^\alpha) \ R \ , \nonumber \\
R_{12} \ R_{13}\ R_{23} = R_{23}\ R_{12}\ ,\hspace{1.8cm}
\end{eqnarray}
\noindent o\`u $R_{12} = R \otimes 1,\ R_{23} = 1 \otimes R$, etc... La dernière identité est appelée \emph{l'équation du pentagone}. Les deux autres permettent de définir des comultiplications sur les images de représentations de $A$ et $A^*$, et donc des opérateurs de Clebsch-Gordan (cf. Remarque \ref{Tannaka}). Notons que les relations (\ref{repadj}) impliquent 
$$(\pi_A \otimes \pi_{A^*})(R) = (\epsilon^* \otimes \Delta^*)\ (m^* \otimes \epsilon^*)\ \in {\rm End}_k((A^*)^{\otimes 2})\ .$$
\noindent (Rappelons que $\epsilon^*$ est l'unité de $A^*$). L'antipode $S$ de $A$ permet d'inverser $(\pi_A \otimes \pi_{A^*})(R)$ :
$$(\pi_A \otimes \pi_{A^*})(R)^{-1} = (\epsilon^* \otimes \Delta^*)\ (\epsilon^* \otimes S^* \otimes \epsilon^*) \ (m^* \otimes \epsilon^*)\ ,$$
\noindent o\`u $S^*$ est l'antipode de $A^*$. Enfin, on peut définir pour les doubles de Heisenberg un analogue de la méthode FRT, qui permet de construire des bigèbres cotressées à partir de solutions inversibles de l'équation de Yang-Baxter quantique \cite[\S 8]{Kas}. On aboutit alors au théorème suivant. Pour sa preuve, on renvoie le lecteur à \cite[Th. 5.7]{Dav}, qui donne un énoncé plus général pour des solutions de l'équation du pentagone sur des modules de Hopf quelconques.

\begin{teo1} \label{FRT} {\rm \cite[Th. 5.7]{Dav}} Soit $(A,1,\epsilon,m,\Delta)$ comme ci-dessus. L'image $(\pi_A \otimes \pi_{A^*})(R)$ de l'élément canonique de $\mathcal{H}(A)$ dans ${\rm End}_k((A^*)^{\otimes 2})$ est l'unique solution de l'équation du pentagone. Les relations {\rm (\ref{pentagone})} vérifiées par $(\pi_A \otimes \pi_{A^*})(R)$ définissent la structure d'algèbre de $\mathcal{H}(A)$.
\end{teo1}

\begin{conv} \label{brousouf} {\rm Dans la suite, pour des raisons techniques (essentiellement pour pouvoir inverser $h$ dans la proposition \ref{defdoubleh}), on a besoin d'introduire le corps des fractions $\mc((h))$  de $\mc[[h]]$, qui est un anneau intègre. Une algèbre topologique $A$ sur $\mc[[h]]$ sera donc considérée comme $\mc((h))$-module $\mc((h)) \otimes_{\mc[[h]]} A$ .}
\end{conv}

\noindent On va maintenant construire explicitement le double de Heisenberg $\mathcal{H}_h(b(2,\mc))$ de $A = U_h(b(2,\mc))$. La sous-algèbre de Borel quantifiée \emph{positive} $U_h(b(2,\mc))$ de $U_h(sl(2,\mc))$ est la $\mc((h))$-QUE-algèbre de Hopf topologiquement engendrée sur $\mc[[h]]$ par $H$ et $D$ vérifiant la relation $HD - DH = D$, et munie des comultiplications et antipodes définies par \cite[\S 6.4]{CP}, \cite[\S 17]{Kas} :
$$\begin{array}{l}
\Delta(H) = H \otimes 1 + 1 \otimes H\ , \quad \Delta(D) = 1 \otimes D + D \otimes e^{hH}\ ,\\
S(H) = -H \ , \quad S(D) = -De^{-hH}\ ,\quad \epsilon(H) = \epsilon(D) = 0\ , \quad \epsilon(1) = 1\ .
\end{array}$$
\noindent (On note $1$ l'unité de $U_h(b(2,\mc))$ et de $\mc((h))$). La $\mc((h))$-algèbre QUE-duale $U_h^*(b(2,\mc))$ de $U_h(b(2,\mc))$ est isomorphe en tant qu'algèbre de Hopf topologique sur $\mc[[h]]$ à la sous-algèbre de Borel quantifiée \emph{négative} de $U_h(sl(2,\mc))$, munie de la comultiplication opposée \cite[Prop. 8.3.2]{CP}. Donc $U_h^*(b(2,\mc))$ est topologiquement engendrée sur $\mc[[h]]$ par $\bar{H}$ et $\bar{D}$ vérifiant la relation $\bar{H}\bar{D} - \bar{D}\bar{H} = -h\bar{D}$, avec :
$$\begin{array}{ll}
\Delta(\bar{H}) = \bar{H} \otimes 1 + 1 \otimes \bar{H}\ , \quad \Delta(\bar{D}) = 1 \otimes \bar{D} + \bar{D} \otimes e^{-\bar{H}}\ ,\\
S(\bar{H}) = -\bar{H} \ , \quad S(\bar{D}) = -\bar{D}e^{\bar{H}}\ ,\quad \epsilon(\bar{H}) = \epsilon(\bar{D}) = 0\ , \quad \epsilon(1) = 1\ .
\end{array}$$
\noindent Clairement, l'application : $\bar{H} \rightarrow -hH$, $\bar{D} \rightarrow D$ est un isomorphisme des $\mc((h))$-algèbres topologiques $U_h^*(b(2,\mc))$ et $U_h(b(2,\mc))$.

\begin{prop} \label{defdoubleh} Le double de Heisenberg $\mathcal{H}_h(b(2,\mc))$ de $U_h(b(2,\mc))$ est isomorphe à la $\mc((h))$-algèbre topologiquement engendrée sur $\mc[[h]]$ par $H,\ D,\ \bar{H}$ et $\bar{D}$ satisfaisant les relations :
$$\begin{array}{ll}
HD - DH = D\ , \quad \bar{H} \bar{D} - \bar{D} \bar{H} = -h\bar{D}\ , \\
H\bar{H} - \bar{H} H = 1\ , \quad \ D\bar{H} = \bar{H} D\ ,\\
H\bar{D} - \bar{D} H = -\bar{D}\ ,\quad  D\bar{D} - \bar{D} D = (1 - q) \ e^{hH} \ \ \mathrm{avec} \ q = e^{-h}\ .
\end{array}$$
\noindent Pour cette pr\'esentation, son élément canonique s'écrit
$$R_h = e^{H \otimes \bar{H}}\ (D \otimes \bar{D};q)_{\infty}^{-1}\ .$$
\end{prop}

\bigskip

\noindent {\it Démonstration}. Soit  $e_{m,n} =H^mD^n/m!(q)_n$ une base linéaire de $U_h(b(2,\mc))$, et considérons la base duale $e^{m,n} = \bar{H}^m\bar{D}^n$ pour le crochet de dualité $\langle e_{m,n},e^{k,l} \rangle = \delta_{m,k}\ \delta_{n,l}$. En utilisant la définition de $\mathcal{H}_h(b(2,\mc))$ par les relations (\ref{relations2}), et notamment la troisième d'entre elles, on trouve le résultat suivant.  

\noindent Les seules constantes de structure du type $\Delta_D^{x,y}$ et $m_{z,x}^{H}$ qui sont simultanément non nulles sont $\Delta_D^{1,D} = 1$ et $m_{H,1}^{H} = 1$. Comme $\bar{H}$ est dual à $H$, on en déduit que $D\bar{H} = \bar{H} D$. 

\noindent De même, les paires $(\Delta_H^{x,y},m_{z,x}^{D})$ simultanément non nulles sont $(\Delta_H^{1,H},m_{D,1}^{D}) = (1,1)$ et $(\Delta_H^{H,1},m_{D,H}^{D}) = (1,-1)$. Comme on a évidemment $m_{D,1}^{D} = m_{D(1-q)^{-1},1}^{D(1-q)^{-1}}$, que la transformation  $D \rightarrow D(1-q)^{-1}$ n'affecte ni la valeur de $\Delta_H^{1,H}$ ni celle de $\Delta_H^{H,1}$, et que $\bar{D}$ est dual à $D(1-q)^{-1}$, on trouve donc $H\bar{D} =\bar{D} H - \bar{D}.1 =\bar{D} H - \bar{D}$. 

\noindent Pour la relation $H\bar{H} - \bar{H} H = 1$, on procède de même. Par contre, la dernière relation montre l'importance du choix des bases duales pour obtenir des relations agréables : les paires $(\Delta_D^{x,y},m_{z,x}^{D})$ simultanément non nulles sont $(\Delta_D^{1,D},m_{D,1}^{D}) = (1,1)$ et $(\Delta_D^{D,e^{hH}},m_{1,D}^{D}) = (1,1)$. Mais ici on a $\Delta_D^{D,e^{hH}} = (1-q)\Delta_D^{D(1-q)^{-1},e^{hH}}$, et comme $\bar{D}$ est dual à $D(1-q)^{-1}$ on obtient $D\bar{D} = \bar{D}D - (1-q)\ e^{hH}$. 

\smallskip

\noindent Montrons que la formule proposée pour $R_h$ vérifie l'équation du pentagone (\ref{pentagone}). Par le théorème \ref{FRT}, on en déduira que $R_h$, considéré comme endomorphisme sur le $k$-module $(A^*)^{\otimes 2} = \left(U_h^*(b(2,\mc)) \right)^{\otimes 2}$, est l'élément canonique de $\mathcal{H}_h(b(2,\mc))$. Dans la suite, pour simplifier, on note $H_i,\ \bar{H}_i,\ D_i,\ \bar{D}_i$ l'action de $H,\ \bar{H},\ D$ et $\bar{D}$ sur le i-ème facteur tensoriel. Considérons l'inverse de la relation du pentagone, pour $R_h^{-1}$. Elle est équivalente à 
\begin{eqnarray}
(D_2\bar{D}_3;q)_{\infty} e^{-H_2\bar{H}_3}\ (D_1\bar{D}_3;q)_{\infty} e^{-H_1\bar{H}_3}\ (D_1\bar{D}_2;q)_{\infty} e^{-H_1\bar{H}_2} = \hspace{3cm} \nonumber \\ (D_1\bar{D}_2;q)_{\infty} e^{-H_1\bar{H}_2}\ (D_2\bar{D}_3;q)_{\infty} e^{-H_2\bar{H}_3}\ .\label{qdilopent} \nonumber
\end{eqnarray}
\noindent  La relation $D\bar{H} = \bar{H} D$ implique que le membre de droite de est égal à
$$(D_1\bar{D}_2;q)_{\infty}(D_2\bar{D}_3;q)_{\infty}\ e^{-H_1\bar{H}_2}e^{-H_2\bar{H}_3}\ .$$
\noindent Pour tous $m,\ n$ entiers on prouve par récurrence que $\bar{H}^m\bar{D}^n =\bar{D}^n\ (\bar{H} - nh1)^m$ . Donc 
$$(H_2\bar{H}_3)^m \ (D_1\bar{D}_3)^n = D_1^n H_2^m \bar{H}_3^m\bar{D}_3^n = D_1^n H_2^m \bar{D}_3^n\ (\bar{H}_3 - nh1_3)^m\ ,$$ 
\noindent et on a :
$$\begin{array}{lll}
e^{-H_2\bar{H}_3}\ (D_1\bar{D}_3;q)_{\infty} & = & \sum_{m \geq 0} \frac{(-H_2\bar{H}_3)^m}{m!} \ \cdot \ \sum_{n \geq 0} \frac{q^{\frac{n(n-1)}{2}}(-D_1\bar{D}_3)^n}{(q)_n} \\ \\
& = & \sum_{m,n \geq 0} \frac{q^{\frac{n(n-1)}{2}}}{m!(q)_n}\ (-D_1 \bar{D}_3)^n\ \left( -H_2(\bar{H}_3 - nh1_3) \right)^m  \\ \\
& = & \sum_{n \geq 0} \frac{q^{\frac{n(n-1)}{2}}}{(q)_n}\ (-D_1 \bar{D}_3)^n \exp( -H_2(\bar{H}_3 - nh1_3)) \\ \\
& = & \sum_{n \geq 0} \frac{q^{\frac{n(n-1)}{2}}}{(q)_n}\ (- e^{hH_2}D_1 \bar{D}_3)^n \exp( -H_2\bar{H}_3) \\ \\
& = & (e^{hH_2}D_1 \bar{D}_3;q)_{\infty}\ \exp( -H_2\bar{H}_3)\ . 
\end{array}$$ 
\noindent De la même façon, on a $H^mD^n = D^n\ (H + n1)^m$ et $H^m\bar{D}^n = \bar{D}^n\ (H - n1)^m$, d'o\`u l'on tire respectivement les identités :
$$\begin{array}{lll}
e^{-H_1\bar{H}_3}\ (D_1\bar{D}_2;q)_{\infty} & = & \sum_{m,n \geq 0} \frac{(-H_1\bar{H}_3)^m}{m!} \ \cdot \ \frac{q^{\frac{n(n-1)}{2}}(-D_1\bar{D}_2)^n}{(q)_n} \\ \\
& = & \sum_{m,n \geq 0} \frac{q^{\frac{n(n-1)}{2}}}{m!(q)_n}\ (-D_1 \bar{D}_2)^n\ \left( -\bar{H}_3(H_1 + n1_1) \right)^m  \\ \\
& = & \sum_{n \geq 0} \frac{q^{\frac{n(n-1)}{2}}}{(q)_n}\ (-D_1 \bar{D}_2)^n \exp(-\bar{H}_3(H_1 + n1_1))\ , 
\end{array}$$  
$$\begin{array}{lll}
e^{-H_2\bar{H}_3}\ e^{-H_1\bar{H}_3}\ (D_1\bar{D}_2;q)_{\infty} & = & \sum_{m,n \geq 0} \frac{q^{\frac{n(n-1)}{2}}}{m!(q)_n}\ (-H_2\bar{H}_3)^m  (-D_1 \bar{D}_2)^n \\ & & \hspace{4.5cm} \times \ \ e^{-n\bar{H}_3}e^{-H_1\bar{H}_3} \\
& = & \sum_{m,n \geq 0} \frac{q^{\frac{n(n-1)}{2}}}{m!(q)_n}\ (-D_1 \bar{D}_2)^n (-\bar{H}_3(H_2 - n1_2))^m \\  & & \hspace{4.5cm} \times  \ \ e^{-n\bar{H}_3}e^{-H_1\bar{H}_3} \\
& = & \sum_{n \geq 0} \frac{q^{\frac{n(n-1)}{2}}}{(q)_n}\ (-D_1 \bar{D}_2)^n \ e^{-H_2\bar{H}_3}e^{-H_1\bar{H}_3} \\ \\
& = & (D_1 \bar{D}_2;q)_{\infty} \ e^{-H_2\bar{H}_3}e^{-H_1\bar{H}_3}\ .
\end{array}$$
\noindent Avec ces trois relations de commutation, on peut transformer le membre de gauche en
$$(D_2 \bar{D}_3;q)_{\infty}(e^{hH_2}D_1 \bar{D}_3;q)_{\infty}(D_1 \bar{D}_2;q)_{\infty}\ e^{-H_1\bar{H}_3}e^{-H_2\bar{H}_3}e^{-H_1\bar{H}_2}\ .$$
\noindent Pour montrer que l'équation du pentagone pour $R$ est vérifiée, il suffit donc de montrer les deux identités suivantes :
\begin{eqnarray}
(D_2 \bar{D}_3;q)_{\infty}(e^{hH_2}D_1 \bar{D}_3;q)_{\infty}(D_1 \bar{D}_2;q)_{\infty} = (D_1 \bar{D}_2;q)_{\infty}(D_2 \bar{D}_3;q)_{\infty} \label{egal1} \ ,\\ \nonumber \\
e^{-H_1\bar{H}_3}e^{-H_2\bar{H}_3}e^{-H_1\bar{H}_2} = e^{-H_1\bar{H}_2} \ e^{-H_2\bar{H}_3}\label{egal2} \ .\hspace{4.3cm}
\end{eqnarray}
\noindent L'équation (\ref{egal2}) est une conséquence de la formule classique de Baker-Campbell-Haussdorff appliquée à la sous-algèbre de Lie complexe de $\mathcal{H}_h(b(2,\mc))$ engendrée par $H$ et $\bar{H}$ (\emph{l'algèbre de Heisenberg}) :
$$\begin{array}{lll}
\exp(-H_1\bar{H}_2) \ \exp(-H_2\bar{H}_3) & = & \exp(-H_1\bar{H}_2 - H_2\bar{H}_3 + \frac{1}{2} H_1\ [\bar{H_2},H_2]\ \bar{H}_3) \\ \\
                                         & = & \exp(-H_1\bar{H}_2 - H_2\bar{H}_3 - \frac{1}{2} H_1\bar{H}_3) \\ \\
                                         & = & \exp(-H_1\bar{H}_3)\exp(-H_2\bar{H}_3)\exp(-H_1\bar{H}_2)\ .
\end{array}$$ 
\noindent Posons maintenant $U = D_2\bar{D}_3$ et $V = D_1\bar{D}_2$ . On a 
$$[U,V] = UV-VU = D_1D_2\bar{D}_2\bar{D}_3 - D_1\bar{D}_2D_2\bar{D}_3 = [D_2,\bar{D}_2] D_1\bar{D}_3 = (1-q)e^{hH_2}D_1\bar{D}_3\ ,$$
\noindent et on vérifie sans difficultés que $W = UV-qVU$ commute avec $U$ et $V$ :
$$[W,U] = [W,V] = 0\ .$$
\noindent Alors l'équation (\ref{egal1}) est une conséquence du lemme \ref{qdilogen}, car on peut l'écrire :
\begin{equation}\label{bobof}
(U;q)_{\infty} \left(\frac{[U,V]}{(1-q)};q\right)_{\infty}(V;q)_{\infty} = (V;q)_{\infty}(U;q)_{\infty}\ .
\end{equation}
\noindent On a donc achevé la preuve de la proposition. \hfill $\Box$

\subsection{Rationalisation et dilogarithme cyclique}\label{rationnalisation}

\noindent Considérons les éléments $E = e^{hH},\ \bar{E} = e^{-\bar{H}} \in \mathcal{H}_h(b(2,\mc))$. En utilisant les mêmes relations de commutation que dans la preuve de la proposition \ref{defdoubleh}, on trouve :
$$\begin{array}{l}
ED = e^{hH}D = \sum_{n = 0}^{\infty} \frac{1}{n!}\ (hH)^nD = \sum_{n = 0}^{\infty} \frac{1}{n!}\ D(h(H + 1))^n  = D e^{h + hH} = q^{-1} DE\ ,\\ \\
\bar{E}\bar{D} = e^{-\bar{H}}\bar{D} = \sum_{n = 0}^{\infty} \frac{1}{n!}\ (-\bar{H})^n\bar{D} = \sum_{n = 0}^{\infty} \frac{1}{n!}\ \bar{D}(-\bar{H} +h ))^n  = \bar{D} e^{h - \bar{H}} = q^{-1} \bar{D}\bar{E}\ ,\\ \\
\left. \begin{array}{l} E\bar{E} = e^{hH}e^{-\bar{H}} = e^{hH - \bar{H} - \frac{h}{2}\ [H,\bar{H}]} \\
                 \bar{E}E = e^{hH - \bar{H} + \frac{h}{2}\ [H,\bar{H}]} 
        \end{array} \right\rbrace \ \Longrightarrow \bar{E}E = e^h E\bar{E} = q^{-1} E\bar{E},\ \\ \\
E\bar{D} = \sum_{n = 0}^{\infty} \frac{1}{n!}\ (hH)^n\bar{D} = \sum_{n= 0}^{\infty} \frac{1}{n!}\ \bar{D}(h(H - 1))^n  = \bar{D} e^{hH - h} = q \bar{D}E\ ,
\end{array}$$
\noindent La relation $D\bar{E} = \bar{E}D$ est une conséquence immédiate de la relation $D\bar{H} = \bar{H}D$. De façon parallèle à la construction des formes intégrales de groupes quantiques \cite[\S 9]{CP}, on aboutit donc à la définition suivante :  

\begin{defi} \label{defdoubleq} La forme intégrale de $\mathcal{H}_h(b(2,\mc))$ est la $\mc[q,q^{-1}]$-algèbre $\mathcal{H}_q = \mathcal{H}_q(b(2,\mc))$ engendrée par $E,\ E^{-1},\ \bar{E},\ \bar{E}^{-1},\ D$ et $\bar{D}$ satisfaisant les relations :
$$\begin{array}{ll}
EE^{-1} = E^{-1}E = 1\ ,\\
DE = qED\ , \quad \bar{D} \bar{E} = q\bar{E} \bar{D}\ , \\
E\bar{E} = q \bar{E}E\ ,\quad D\bar{E} = \bar{E} D\ ,\\
E\bar{D} = q \bar{D}E\ ,\quad D\bar{D} - \bar{D} D = (1 - q) \ E \ .
\end{array}$$
\noindent Soit $\varphi : \mc[q,q^{-1}] \rightarrow \mc$ l'homomorphisme qui envoie $q$ sur $\epsilon$. Pour tout nombre complexe $\epsilon$ non nul, on dit que l'algèbre $\mathcal{H}_{\epsilon} = \mathcal{H}_q \otimes_{\varphi} \mc$ est la spécialisation de $\mathcal{H}_q$ en $q=\epsilon$.
\end{defi}
\noindent Clairement, l'élément canonique $R_h \in \mathcal{H}_h(b(2,\mc))^{\otimes 2}$ ne définit pas un élément de $\mathcal{H}_q^{\otimes 2}$. D'abord parce que $(D \otimes \bar{D};q)_{\infty}$ est une série \emph{infinie} en les générateurs $D$ et $\bar{D}$, ensuite et de façon plus grave parce que \emph{$e^{H \otimes \bar{H}}$ ne peut pas s'écrire en termes des générateurs de $\mathcal{H}_q$} (ou même d'une complétion de $\mathcal{H}_q$). On va donc construire un ``analogue cyclique'' de $R_h$, qui impl\'emente l'action de $R_h$ sur les représentations cycliques d'un quotient de $\mathcal{H}_{\omega^{-1}}$ (ce qui nous suffira pour la suite). Un problème similaire, résolu par Tanisaki \cite{Ta},\ \cite[\S 10.1.D]{CP}, a lieu pour définir des solutions de l'équation de Yang-Baxter quantique à partir des formes intégrales de QUE-algèbres. 

\subsubsection{L'analogue cyclique de $e^{H \otimes \bar{H}}$} \label{cycexp}

\noindent Considérons la sous-algèbre de Hopf de $U_h(b(2,\mc))$ topologiquement engendrée sur $\mc[[h]]$ par $H$. On a
$$\Delta(H) = 1 \otimes H + H \otimes 1\ .$$
\noindent En reprenant la preuve de la proposition \ref{defdoubleh}, on voit que la QUE-algèbre duale est topologiquement engendrée par $\bar{H}$, avec le même coproduit. Son double de Heisenberg $\mathcal{H}_h^0$ est la $\mc((h))$-algèbre topologiquement engendrée sur $\mc[[h]]$ par $H$ et $\bar{H}$ et vérifiant la relation de commutation de Heisenberg $H\bar{H} - \bar{H}H = 1$. Son élément canonique est $e^{H \otimes \bar{H}}$. La forme intégrale de $\mathcal{H}_h^0$ est la $\mc[q,q^{-1}]$-sous-algèbre $\mathcal{W}_q^0$ de $\mathcal{H}_q$ engendrée par $E$ et $\bar{E}$ et tels que $E\bar{E} = q \bar{E}E$. D'autre part, on peut naturellement identifier la $\mc[q,q^{-1}]$-algèbre du groupe $\mz$ avec l'algèbre du sous-groupe de $\mathcal{H}_q$ engendré par $E$, en prolongeant lin\'eairement l'isomorphisme $i \in \mz \rightarrow E^i$. Lorsque $q$ est une racine primitive $N$-ième de l'unité, il existe un isomorphisme analogue : 

\begin{lem} Le double de Heisenberg $\mathcal{H}(\mc[\mz/N\mz])$ est isomorphe au quotient de $\mathcal{W}_{\omega^{-1}}^0$ par l'idéal engendré par les relations $E^N = 1$ et $\bar{E}^N = 1$. Son \'el\'ement canonique peut s'\'ecrire
$$S_N = \sum_{i \in\mz/N\mz} e^i \otimes \bar{e}^i =\frac{1}{N}\ \sum_{i,j \in \mz} \omega^{ij}\ e^i\otimes \hat{\bar{e}}^j\ ,$$
o\`u $\textstyle \hat{\bar{e}}^j = \sum_{i \in \mz/N\mz} \omega^{-ij}\ \bar{e}^i$ est la transformée de Fourier de $\bar{e}^i$.
\end{lem}

\noindent On a sp\'ecialis\'e $\mathcal{W}_q^0$ en $q=\omega^{-1}$ pour obtenir exactement cette forme pour $S_N$.

\medskip

\noindent {\it D\'emonstration.} De fa\c{c}on g\'en\'erale, si $G$ un groupe fini et $\mc[G]$ l'algèbre de $G$. on peut munir $\mc[G]$ d'une structure d'algèbre de Hopf, de comultiplication $\Delta(g) = g \otimes g$, d'antipode $S(g) = g^{-1}$ et de counité $\epsilon(g) = 1$ pour tout $g \in G$ \cite[\S 3]{Kas}. Avec (\ref{relations2}), on voit que le double de Heisenberg $\mathcal{H}(\mc[G])$ de $\mc[G]$ est défini par les relations :
\begin{equation}\label{relfini}
e_ge_h = e_{gh}\ ,\quad e^ge^h = \delta_{g,h}e^h,\quad e^ge_h = e_he^{gh}\ ,
\end{equation}
\noindent o\`u $e_g$ et $e^g$ sont des bases linéaires de $\mc[G]$ et $\mc[G]^*$ respectivement, et $\delta_{g,h}$ est égal à $1$ si $g=h$ et $0$ sinon. Considérons le cas $G=\mz/N\mz$. Soient $e$ un générateur de $\mz/N\mz$ de dual $\bar{e}$ dans $\mc[G]^*$, et $\textstyle S_N = \sum_{i \in\mz/N\mz} e^i \otimes \bar{e}^i$ l' élément canonique de $\mathcal{H}(\mc[\mz/N\mz])$ . On a 
$$\bar{e}^i = \frac{1}{N} \sum_{j \in \mz/N\mz} \omega^{ij} \ \hat{\bar{e}}^j \quad \Longrightarrow \quad \hat{\bar{e}}^j = \sum_{i \in \mz/N\mz} \omega^{-ij} \ \bar{e}^i$$
$$e\hat{\bar{e}} = \sum_{i \in \mz/N\mz} \omega^{-i} \ e \bar{e}^i = \sum_{i \in \mz/N\mz} \omega^{-i} \ \bar{e}^{i-1}e = \omega^{-1} \sum_{i \in \mz/N\mz} \omega^{-(i-1)} \ \bar{e}^{i-1}e = \omega^{-1} \hat{\bar{e}}e\ .$$
$$\hat{\bar{e}}^{i} \hat{\bar{e}}^{j} = \sum_{k,l \in \mz/N\mz}\omega^{-ik}\omega^{-jl}\ \bar{e}^k \bar{e}^l = \sum_{l \in \mz/N\mz} \omega^{-(i+j)l} \ \bar{e}^l = \hat{\bar{e}}^{i+j}\ .$$
\noindent On utilise (\ref{relfini}) dans la seconde \'egalit\'e de la deuxième et de la troisième ligne. Donc l'application d\'efinie par $e \rightarrow E$ et $\hat{\bar{e}} \rightarrow \bar{E}$ induit un isomorphisme de $\mathcal{H}(\mc[\mz/N\mz])$ vers $\mathcal{W}_q^0/(E^N = 1)(\bar{E}^N = 1)$.  \hfill $\Box$

\begin{cor} \label{élémzn} L'image de $S_N$ par une représentation irréductible et cyclique est de la forme :
$$\sum_{i \in \mz/N\mz} Z_1^{-i} \widehat{Y}_2^i = \frac{1}{N}\ \sum_{i,j \in \mz} \omega^{ij}\ Z_1^{-i}\ Y_2^j\ ,$$ 
\noindent o\`u $Y$ et $Z$ sont des matrices unipotentes de taille $N$ et d'ordre $N$ telles que $ZY = \omega YZ$, et $\widehat{Y}^i$ désigne la transformée de Fourier inverse normalisée de $Y^i$.
\end{cor}

\noindent {\it Démonstration.} L'image de $e$ et $\hat{\bar{e}}$ par une représentation cyclique de $\mathcal{H}(\mc[\mz/N\mz])$ sont nécessairement des matrices unipotentes d'ordre au plus $N$, qui sont de taille et d'ordre exactement \'egal à $N$ si la représentation est de plus irréductible. L'application définie par $e^i \rightarrow Z^{-i}$ et $\bar{e}^i \rightarrow \widehat{Y}^i$ est clairement une représentation de $\mathcal{H}(\mc[\mz/N\mz]) \cong \mathcal{W}_q^0/(E^N = 1)(\bar{E}^N = 1)$.\hfill $\Box$

\medskip

\noindent On peut d\'efinir l'image de $S_N$ dans une extension de $(\mathcal{W}_{\omega^{-1}}^0)^{\otimes 2}$ de la fa\c{c}on suivante. L'algèbre $\mathcal{W}_{\omega^{-1}}^0$ est sans diviseur de zéros : par exemple, c'est une extension de Ore d'un anneau de polynômes \cite[\S 4]{Kas}. Donc son centre $\mathcal{Z}(\mathcal{W}_{\omega^{-1}}^0)$ est un anneau intègre et on peut considérer son corps quotient $Q(\mathcal{Z}(\mathcal{W}_{\omega^{-1}}^0))$ . Posons 
\begin{equation}\label{corps}
Q(\mathcal{W}_{\omega^{-1}}^0) := \mathcal{W}_{\omega^{-1}}^0 \otimes_{\mathcal{Z}(\mathcal{W}_{\omega^{-1}}^0)} Q(\mathcal{Z}(\mathcal{W}_{\omega^{-1}}^0)) \ ,
\end{equation}
\noindent et soient $E' = c_{E}^{-1}E,\ \bar{E} ' = c_{\bar{E}}^{-1}\bar{E} \in Q(\mathcal{W}_{\omega^{-1}}^0)$, o\`u $c_{E},\ c_{\bar{E}} \in \mathcal{Z}(\mathcal{W}_{\omega^{-1}}^0)$ sont définis respectivement par $c_{E}^N = E^N$ et $c_{\bar{E}}^N = \bar{E}^N$. Les éléments $c_{E}$ et $c_{\bar{E}}$ existent car $\mathcal{Z}(\mathcal{W}_{\omega^{-1}}^0)$ est intégralement clos \cite[Prop. 11.1.2]{CP}. Consid\'erons l'\'el\'ement $\Upsilon_{\omega} \in Q(\mathcal{W}_{\omega^{-1}}^0)^{\otimes 2}$, image de $S_N$, défini par :

\begin{equation} \label{Upsilon}
\Upsilon_{\omega} = \frac{1}{N}\ \sum_{i,j \in \mz} \omega^{ij}\ (E')^{i} \otimes (\bar{E}')^j \ .
\end{equation} 
\noindent En utilisant les relations du d\'ebut du \S \ref{rationnalisation}, on v\'erifie facilement que l'action par conjugaison de $e^{H \otimes \bar{H}}$ sur $\mathcal{H}_q^{\otimes 2}$ d\'efinit un automorphisme $\mathcal{E}_q$. Lorsque $q=\omega^{-1}$, l'action par conjugaison de $\Upsilon_{\omega}$ sur $\mathcal{H}_{\omega^{-1}}^{\otimes 2}$ s'identifie à $\mathcal{E}_{\omega^{-1}}$. Donc $\Upsilon_{\omega}$ est l'analogue cyclique de $e^{H \otimes \bar{H}}$. On trouve la formule du lemme \ref{élémzn} pour $\Upsilon$ en consid\'erant une repr\'esentation irr\'eductible $\rho$ de $\mathcal{W}_{\omega^{-1}}^0$ avec $\rho(E')=Z^{-1}$ et $\rho(\bar{E }')=Y$.

\subsubsection{Les analogues cycliques de $(D \otimes \bar{D};q)_{\infty}$ et $R_h$}\label{introdc}

\noindent Soient $\bar{\mathcal{H}}_q$ le quotient de $\mathcal{H}_q$ par l'idéal engendré par la relation $D\bar{D} = q\bar{D}D$, et $\widehat{\bar{\mathcal{H}}_q}$ l'espace vectoriel des séries formelles en les générateurs de $\bar{\mathcal{H}}_q$, et à coefficients dans $\mc[q,q^{-1}]$. Alors $(D \otimes \bar{D};q)_{\infty}$ définit un élément de $\widehat{\bar{\mathcal{H}}_q}^{\otimes 2}$. De plus, la relation (\ref{qdilocentral}) devient l'équation du $q$-dilogarithme (\ref{qdiloeq}) pour la sous-algèbre topologique de $\widehat{\bar{\mathcal{H}}_q}^{\otimes 3}$ engendrée par $U = D_2\bar{D}_3$ et $V =D_1\bar{D}_2$, car $UV = qVU$. 

\smallskip

\noindent Considérons la spécialisation de $\widehat{\bar{\mathcal{H}}_q}$ en $q=\exp(-\epsilon/N^2)\omega^{-1}$. On a vu au \S \ref{asy} que l'analogue cyclique de $(x;q)_{\infty}$ sur $\widehat{\mathcal{W}}_{\omega^{-1}}$ est l'opérateur inversible $\Phi(x) \in {\rm End}( \widehat{\mathcal{W}}_{\omega^{-1}})$. Donc l'analogue cyclique de $(D \otimes \bar{D};q)_{\infty}$ est $\Phi\left( D \otimes \bar{D} \right) \in {\rm End}\left(\widehat{\bar{\mathcal{H}}_{\omega^{-1}}}^{\otimes 2}\right)$. Pour écrire son équation du dilogarithme cyclique (Th\'eorème \ref{eqdiloasy}), il faut pouvoir inverser certains éléments de $\widehat{\bar{\mathcal{H}}_q}$. On le fait de la façon suivante. 

\noindent Notons $\widehat{\mathcal{W}}_q^+$ la sous-algèbre de $\widehat{\bar{\mathcal{H}}_q}$ topologiquement engendrée par $D$ et $\bar{D}$. On pose
$$Q(\widehat{\bar{\mathcal{H}}_{\omega^{-1}}}) := \widehat{\bar{\mathcal{H}}_{\omega^{-1}}} \otimes_{\mathcal{Z}(\widehat{\mathcal{W}}_{\omega^{-1}}^+)} Q(\mathcal{Z}(\widehat{\mathcal{W}}_{\omega^{-1}}^+)) \otimes_{\mathcal{Z}(\widehat{\mathcal{W}}_{\omega^{-1}}^0)} Q(\mathcal{Z}(\widehat{\mathcal{W}}_{\omega^{-1}}^0))\ ,$$
\noindent o\`u $Q(\mathcal{Z}(\widehat{\mathcal{W}}_{\omega^{-1}}^0))$ et $Q(\mathcal{Z}(\widehat{\mathcal{W}}_{\omega^{-1}}^+))$ sont définis comme dans (\ref{corps}) et $\mathcal{W}_{\omega^{-1}}^0$ est défini au \S \ref{cycexp}. Soit $\Psi_{\omega}\left( D \otimes \bar{D} \right) := \Phi^{-1}\left( D \otimes \bar{D} \right) \in {\rm End}\left( Q(\widehat{\bar{\mathcal{H}}_{\omega^{-1}}})^{\otimes 2} \right)$ . Consid\'erons
\begin{equation} \label{élémcanomega}
R_\omega = \Upsilon_{\omega} \cdot \Psi_{\omega}\left( D \otimes \bar{D} \right) \in {\rm End}\left( Q(\widehat{\bar{\mathcal{H}}_{\omega^{-1}}})^{\otimes 2} \right)
\end{equation}
\noindent On sait que $\Upsilon_{\omega}$ v\'erifie l'\'equation du pentagone. On montrera à la proposition \ref{6jdesc} que l'image de $R_\omega$ sur les représentations cycliques de $\mathcal{H}_{\omega^{-1}}$ (i.e. les $6j$-symboles cycliques) v\'erifie aussi l'équation du pentagone. La même preuve montre que $R_\omega$ satisfait aussi l'équation du pentagone. Donc $R_\omega$ est l'analogue cyclique de $R_h$. 

\subsubsection{Description explicite des dilogarithmes cycliques}\label{descdc}

\noindent Avec le lemme \ref{gdilocyc}, on aboutit naturellement à la d\'efinition suivante, originellement due à L. Faddeev et R. Kashaev \cite{FK} :

\begin{defi} \label{dilocyc}
Soient $E$ un espace vectoriel complexe de dimension finie, $A \in {\rm End}(E)$ un endomorphisme de spectre $\{-\omega^n,\ n=0,\ldots,\ N-1\}$, et des nombres complexes non nuls $a,\ b$ et $c$ tels que $a^N + b^N = c^N$. Un dilogarithme cyclique d'ordre $N$ sur $E$ est un opérateur $\Psi_{a,b,c}(A) \in {\rm End}(E)$ qui commute avec $A$ et dont le spectre est de la forme $\{ h_{a,b,c}\ \omega(a,b,c \vert n)\ ,\ n=0,\ldots,\ N-1\}$, o\`u $h_{a,b,c}$ est une fonction complexe non nulle quelconque et 
$$\omega(a,b,c \vert n) = \prod_{j=1}^n \frac{b}{c-a\omega^j}\ .$$
\end{defi}

\noindent De cette d\'efinition, on tire imm\'ediatement l'identit\'e
\begin{equation} \label{eqfonc}
\Psi_{a,b,c}(\omega^{-1}A)\ \Psi_{a,b,c}(A)^{-1} = \frac{(c - aA)}{b}\ ,
\end{equation}
\noindent qui par identification des coefficients donne
\begin{equation} \label{dilocycbis}
\Psi_{a,b,c}(A) = h_{a,b,c}\ \sum_{n=0}^{N-1} A^n\ \prod_{s=1}^n \ \frac{a}{c - \omega^{-s}b} = h_{a,b,c}\ \sum_{n=0}^{N-1} \left(\frac{a}{c}A\right)^n\ \prod_{s=1}^n \ \frac{1}{1 - \omega^{-s}\ \frac{b}{c}} \ .
\end{equation}
\noindent Le théorème suivant est équivalent au théorème \ref{eqdiloasy}. (On en donne une preuve car il n'y en a pas dans \cite{FK}).

\smallskip

\begin{teo1} \label{FadKas} {\rm \cite{FK}} \ Soient $U,\ V \in {\rm End}(E)$ tels que $U^N = V^N = -1$ et $UV = \omega VU$. Posons $\Psi_i = \Psi_{a_i,b_i,c_i}$, o\`u $x_i = a_i/b_i$ et $y_i = c_i/b_i$ sont tels que
\begin{eqnarray}\label{paramètres}
y_0y_2 = y_1y_4\ ,\quad y_1 = y_2y_3\ ,\quad x_3 = x_0x_1\ ,\nonumber \\
x_2 = x_1y_4\ ,\quad x_4 = x_0y_2\ .\hspace{1cm}
\end{eqnarray}
\noindent Alors on a l'identité suivante :
$$\Psi_0(V)\ \Psi_1(U) = \Psi_2(U)\ \Psi_3(-UV)\ \Psi_4(V)\ $$
\noindent lorsque les déterminants des deux membres sont égaux. 
\end{teo1}

\medskip

\noindent {\it Démonstration}. On montre d'abord que $\Psi_2(U)^{-1}\ \Psi_0(V)\ \Psi_1(U)\ \Psi_4(V)^{-1}$ commute avec $-UV$. On a :
$$\begin{array}{l}
\hspace{-2cm} \Psi_2(U)^{-1}\ \Psi_0(V)\ \Psi_1(U)\ \Psi_4(V)^{-1}\ (-UV) \\ \\
= \Psi_2(U)^{-1}\ \Psi_0(V)\ \Psi_1(U)\ (-UV)\ \Psi_4(\omega^{-1}V)^{-1} \\ \\
= \Psi_2(U)^{-1}\ \Psi_0(V)\ (-UV)\ \Psi_1(\omega U)\ \Psi_4(\omega^{-1}V)^{-1}\\ \\
=  \Psi_2(U)^{-1}\ (-UV)\ \Psi_0(\omega^{-1}V)\ \Psi_1(\omega U)\ \Psi_4(\omega^{-1}V)^{-1}\\ \\
=  -UV\ \Psi_2(\omega U)^{-1}\ \Psi_0(\omega^{-1}V)\ \Psi_1(\omega U)\ \Psi_4(\omega^{-1}V)^{-1}\ .
\end{array}$$
\noindent En utilisant (\ref{eqfonc}) on peut transformer le produit des dilogarithmes cycliques dans cette dernière expression en :
$$\begin{array}{l}
\hspace{-0.5cm} \Psi_2(U)^{-1}\ \frac{(c_2 - \omega a_2 U)}{b_2} \ \Psi_0(\omega^{-1}V)\ \Psi_1(\omega U)\ \Psi_4(\omega^{-1}V)^{-1} \\ \\
 = \Psi_2(U)^{-1}\ \left( \Psi_0(\omega^{-1}V)\ \frac{c_2}{b_2} - \frac{\omega a_2}{b_2}\ \Psi_0(V)\ U \right) \ \Psi_1(\omega U)\ \Psi_4(\omega^{-1}V)^{-1}\\ \\
 = \Psi_2(U)^{-1}\ \Psi_0(V)\ \left( \left( \frac{c_0 - a_0 V}{b_0}\right)\ \frac{c_2}{b_2} - \frac{\omega a_2}{b_2}\ U \right) \ \Psi_1(\omega U)\ \Psi_4(\omega^{-1}V)^{-1} \\ \\
\end{array}$$
$$\begin{array}{l}
 = \Psi_2(U)^{-1}\ \Psi_0(V)\ \left( \Psi_1(\omega U)\ \left( \frac{c_0c_2}{b_0b_2} - \frac{\omega a_2}{b_2}\ U \right) - \frac{a_0c_2}{b_0b_2} \ \Psi_1(U)\ V \right)\ \Psi_4(\omega^{-1}V)^{-1} \\ \\
 = \Psi_2(U)^{-1}\ \Psi_0(V)\ \Psi_1(\omega U)\ \left( \left( \frac{c_0c_2}{b_0b_2} - \frac{\omega a_2}{b_2}\ U \right) - \frac{a_0c_2}{b_0b_2}\left( \frac{c_1 - \omega a_1 U}{b_1}\right)\ V \right)\\ \hspace{10cm}\times \ \Psi_4(\omega^{-1}V)^{-1}.
\end{array}$$
\noindent Avec les notations de l'énoncé, on a :
$$\begin{array}{c}
y_1y_4 = y_0y_2 \Rightarrow \frac{c_0c_2}{b_0b_2} = \frac{c_1c_4}{b_1b_4}\ ,\quad x_0y_2 = x_4 \Rightarrow \frac{a_0c_2}{b_0b_2} = \frac{a_4}{b_4}\ ,\\ \\
x_1y_4 = x_2 \Rightarrow \frac{a_1c_4}{b_1b_4} = \frac{a_2}{b_2}\ .
\end{array}$$
\noindent Alors le terme entre parenthèses dans la dernière expression est égal à 
$$\begin{array}{l}
\left( \left( \frac{c_0c_2}{b_0b_2} - \frac{\omega a_2}{b_2}\ U \right) + \left( \frac{(c_1 - \omega a_1 U)}{b_1} \right)\ \left( \frac{(c_4 - a_4 V)}{b_4} \right) - \left( \frac{(c_1 - \omega a_1 U)}{b_1} \right)\ \frac{c_4}{b_4} \right) \\ \\
\hspace{5cm} = \left( \frac{(c_1 - \omega a_1 U)}{b_1} \right)\ \left( \frac{(c_4 - a_4 V)}{b_4} \right)\ .
\end{array}$$
\noindent En utilisant (\ref{eqfonc}) encore deux fois, on trouve donc :
$$\begin{array}{l}
\Psi_2(U)^{-1}\ \Psi_0(V)\ \Psi_1(U)\ \Psi_4(V)^{-1}\ (-UV) = \hspace{6cm}\\
\hspace{5cm}(-UV)\ \Psi_2(U)^{-1}\ \Psi_0(V)\ \Psi_1(U)\ \Psi_4(\omega^{-1}V)^{-1}\ .
\end{array}$$
\noindent Par conséquent
$$P(-UV) = \Psi_2(U)^{-1}\ \Psi_0(V)\ \Psi_1(U)\ \Psi_4(V)^{-1}$$
\noindent est une fonctionnelle de $-UV$. Montrons qu'elle vérifie l'équation :
\begin{equation} \label{eqnouv}
P(-UV) \ \left( \frac{c_3 - a_3 (-UV)}{b_3} \right) = P(-\omega^{-1}UV)\ .
\end{equation}
\noindent Considérons le changement de variable $V \rightarrow \omega^{-1}V$ dans l'expression de $P(-UV)$. On a :
$$\begin{array}{l}
\hspace{-0.5cm} \Psi_2(U)^{-1}\ \Psi_0(\omega^{-1}V)\ \Psi_1(U)\ \Psi_4(\omega^{-1}V)^{-1}  \\ \\
 \hspace{0.3cm} =  \Psi_2(U)^{-1}\ \Psi_0(V)\ \left( \frac{c_0 - a_0 V}{b_0} \right)\ \Psi_1(U)\ \Psi_4(\omega^{-1}V)^{-1}\\ \\
 \hspace{0.3cm} =  \Psi_2(U)^{-1}\ \Psi_0(V)\ \left( \frac{c_0}{b_0}\ \Psi_1(U) - \frac{a_0}{b_0}\ \Psi_1(\omega^{-1}U)\ V \right)\ \Psi_4(\omega^{-1}V)^{-1} 
\\ \\
\hspace{0.3cm} = \Psi_2(U)^{-1}\ \Psi_0(V)\ \Psi_1(U)\ \left( \frac{c_0}{b_0} - \frac{a_0}{b_0}\ \left( \frac{c_1 - a_1 U}{b_1} \right)\ V \right)\ \Psi_4(\omega^{-1}V)^{-1} \\ \\
\end{array}$$
$$\begin{array}{l}
 \hspace{-0.2cm} = \Psi_2(U)^{-1}\ \Psi_0(V)\ \Psi_1(U)\ \left( \frac{c_0}{b_0}\ \Psi_4(\omega^{-1}V)^{-1} - \frac{a_0}{b_0}\ \left( \frac{c_1 - a_1 U}{b_1} \right)\ \Psi_4(\omega^{-1}V)^{-1}\ V \right) \\ \\
  \hspace{-0.2cm} = \Psi_2(U)^{-1}\ \Psi_0(V)\ \Psi_1(U)\ \left(\Psi_4(\omega^{-1}V)^{-1}\ \left( \frac{c_0}{b_0} - \frac{a_0c_1}{b_0b_1}\ V \right) + \frac{a_0a_1}{b_0b_1}\ \Psi_4(V)^{-1}\ UV \right) .
\end{array}$$
\noindent Or
$$\begin{array}{c}
\left. \begin{array}{l} y_0y_2 = y_1y_4 \\
                     y_1 = y_2y_3 \end{array} \right\rbrace \Rightarrow \frac{c_4c_3}{b_4b_3} = \frac{c_0}{b_0}\ ,\quad \left. \begin{array}{l}  y_1 = y_2y_3 \\
                         x_4 = x_0y_2 \end{array} \right\rbrace \Rightarrow \frac{a_4c_3}{b_4b_3} = \frac{a_0c_1}{b_0b_1} \\ \\
x_3 = x_0x_1 \Rightarrow \frac{a_0a_1}{b_0b_1} = \frac{a_3}{b_3}
\ .
\end{array}$$
\noindent On en déduit que la dernière expression est égale à :
$$\begin{array}{l}
\Psi_2(U)^{-1}\ \Psi_0(V)\ \Psi_1(U)\ \left(\Psi_4(\omega^{-1}V)^{-1}\ \left( \frac{c_4 - a_4 V}{b_4}\right)\ \frac{c_3}{b_3} + 
\frac{a_3}{b_3}\ \Psi_4(V)^{-1}\ UV \right)\\ \\
\hspace{2cm}  = \Psi_2(U)^{-1}\ \Psi_0(V)\ \Psi_1(U)\ \Psi_4(V)^{-1}\ \left(\frac{c_3 - a_3 (-UV)}{b_3}\right)\ ,
\end{array}$$
ce qui prouve (\ref{eqnouv}). Comme cette équation fonctionnelle détermine $P(-UV)$ a une constante près, il suffit d'imposer que les déterminants de part et d'autre de l'équation du dilogarithme cyclique soient égaux pour que cette constante soit égale à $1$. Ceci achève la preuve du théorème.\hfill $\Box$

\section{$6j$-symboles et dilogarithme cyclique}\label{6jdilocyc}

\subsection{Description explicite des opérateurs de Clebsch-Gordan} \label{symboles}

\noindent Le \emph{module de muliplicit\'e} de repr\'esentations $\rho$, $\mu$ de $\mathcal{W}_N$ est l'ensemble
$$M_{\rho,\mu} = {\rm End}_{\mathcal{W}_N}(V_\rho,V_\mu) = \{ U : V_\rho \rightarrow  V_\mu \ \vert \ U\rho(a) = \mu(a)U,\ \forall \ a \in \mathcal{W}_N \}.$$
\noindent En tant qu'espace de morphismes entre $\mathcal{W}_N$-modules, $M_{\rho,\mu}$ est naturellement muni d'une structure de $\mathcal{W}_N$-module. La proposition \ref{dimfin} montre que la dimension de $M_{\rho,\mu}$ est finie. La proposition \ref{prop1a} $ii)$-$iii)$ implique que pour des représentations irréductibles et cycliques $\rho,\mu,\nu$ de $\mathcal{W}_N$ telles que $(\rho,\mu)$ est une paire régulière, on a 
$${\rm dim}_{\mc}(M_{\nu,\rho \otimes \mu}) = {\rm dim}_{\mc}(M_{\rho \otimes \mu,\nu}) = \left\lbrace \begin{array}{l} 
                                                       N \ \mathrm{si} \ \nu \sim \rho\mu \\
                                                       0 \ \mathrm{sinon}\ .\end{array} \right.$$
\noindent Les éléments de $M_{\rho\mu,\rho \otimes \mu}$ sont des injections, et sont communément appelés \emph{opérateurs de Clebsch-Gordan} (OCG). Les éléments de $M_{\rho \otimes \mu,\rho\mu}$ sont des projections, et sont appelés les OCG \emph{duaux}. Consid\'erons la courbe
$$\Gamma = \{ [x,y,z] \ \vert\  x^N + y^N = z^N \} \subset \mathbb{CP}^2,$$ 
\noindent o\`u $[x,y,z]$ sont les coordonnées homogènes de $\mathbb{CP}^2$. Rappelons que $\omega = \exp(2i\pi/N)$. Pour tout entier positif $n$, soit $\omega$ (la confusion des notations est volontaire) la fonction rationnelle sur $\Gamma$ d\'efinie par
$$\omega(x,y,z \vert n) = \prod_{j=1}^n \frac{y}{z-x\omega^j}\ , [x,y,z] \in \Gamma \setminus \{ [1,0,\omega^j],j=1,\ldots,n\}\ .$$
On notera $\omega(x,y,z \vert m, n) = \omega(x,y,z \vert m-n)\ \omega^{n^2/2}$. On renvoie à \cite{BBP,KMS} pour des d\'etails sur leurs propri\'et\'es. On a immédiatement :
\begin{equation}
\omega(x,y,z \vert m+n)  =  \omega(x,y,z \vert n) \ \omega(x\omega^n,y,z \vert m)\ .\label{dec}
\end{equation} 
\noindent Pour tout $n \geq 1$ et tout $x \ne \omega^j,\ j=-1,\ldots,-n$, posons :
\begin{equation} \label{ompetit}
\omega(x \vert 0) = 1,\ \ \omega(x \vert n) =  \prod_{j=1}^n \frac{1}{1-x\omega^j}\ .
\end{equation}
\noindent On définit la fonction rationnelle $f(x,y \vert z)$ par :
\begin{equation} \label{sommf}
f(x,y \vert z) = \sum_{\sigma = 0}^{N-1} \frac{\omega(x \vert \sigma)}{\omega(y \vert \sigma)}\  z^{\sigma}\ ,
\end{equation}
\noindent o\`u $z^N(1-y^N) = 1 - x^N$. Notons de plus
$$[x] = N^{-1}\ \frac{1-x^N}{1-x}\ ,  \ \ \mathrm{et} \ \ \delta(n) = \left\lbrace \begin{array}{l} 
                           1 \ \mathrm{si} \ n \equiv 0 \pmod N \\
                           0 \ \mathrm{sinon}\ .
                           \end{array} \right.$$
\noindent La fonction $\delta(n)$ est le symbole de Kronecker \emph{réduit modulo $N$}.

\begin{prop} \label{CG1}
Soit $(\rho,\mu)$ une paire régulière de représentations standards. Pour tout $h_{\rho,\mu} \in \mc^*$, les applications linéaires $K_\alpha(\rho,\mu)$ et $\bar{K}^\alpha(\rho,\mu)$, $\alpha = 0,\ldots,\ N-1$, décrites par des matrices de composantes
$$\begin{array}{l}
K_\alpha(\rho,\mu)_{i,j}^k = h_{\rho,\mu} \ \omega^{\alpha j}\omega(a_\rho y_\mu,\frac{y_\rho}{a_\mu},y_{\rho\mu}\vert i,\alpha) \ \delta(i+j-k) \\
 \\
\bar{K}^\alpha(\rho,\mu)_k^{i,j}= \frac{[\frac{a_\rho y_\mu}{y_{\rho\mu}}]}{h_{\rho,\mu}} \ \frac{\omega^{-\alpha j} \ \delta(i+j-k)}{\omega(\frac{a_\rho y_\mu}{\omega},\frac{y_\rho}{a_\mu},y_{\rho\mu}\vert i,\alpha)}
\end{array}$$
\noindent sont des bases respectives de $M_{\rho\mu, \rho\otimes \mu}$ et $M_{\rho \otimes \mu,\rho\mu}$, qu'on appelle \emph{standard}. Ces bases sont duales, i.e. on a :
$$\bar{K}^\alpha(\rho,\mu)K_\beta(\rho,\mu) = \delta (\alpha - \beta) \ id_{V_{\rho\mu}}\ \ \mathrm{et} \ \ \sum_{\alpha=0}^{N-1} K_\alpha(\rho,\mu)\bar{K}^\alpha(\rho,\mu) = id_{V_\rho} \otimes id_{V_\mu}\ .$$
\end{prop}

\medskip

\noindent Le r\^ole du scalaire $h_{\rho,\mu}$ sera pr\'ecis\'e dans la proposition \ref{6jdesc}.

\medskip

\noindent {\it Démonstration.} \  On cherche une famille de solutions linéairement indépendantes $K_\alpha(\rho,\mu) :  V_{\rho\mu} \rightarrow  V_{\rho \otimes \mu},\ \alpha = 0,\ldots,\ N-1$ du système d'équations 
\begin{eqnarray} \label{sys}
 \ (\rho \otimes \mu) (\Delta (a))K_\alpha(\rho,\mu) =  K_\alpha(\rho,\mu)\rho\mu(a),\ a \in \mathcal{W}_N\ .
\end{eqnarray}
\noindent En écrivant $(\ref{sys})$ dans la base canonique de $\mathbb{C}^N \otimes \mathbb{C}^N$ et successivement pour $a=E$ puis $a=D$, on obtient
\begin{eqnarray} \label{eq-1}
(K_\alpha(\rho,\mu)\rho\mu(E))_{i,j}^k & = & a_{\rho\mu}^2\omega^kK_\alpha(\rho,\mu)_{i,j}^k \ ,\nonumber  \\
(\rho\otimes \mu(\Delta (E))K_\alpha(\rho,\mu))_{i,j}^k & = & a_\rho^2a_\mu^2\omega^{i+j} K_\alpha(\rho,\mu)_{i,j}^k\ . 
\end{eqnarray}
\noindent Donc $K_\alpha(\rho,\mu)_{i,j}^k \ne 0$ seulement si $i+j=k$, et  
\begin{eqnarray}\label{eq0}
(K_\alpha(\rho,\mu)\rho\mu(D))_{i,j}^{k-1} & = & a_{\rho\mu}y_{\rho\mu}K_\alpha(\rho,\mu)_{i,j}^k\ , \nonumber \\
(\rho\otimes \mu(\Delta (D))K_\alpha(\rho,\mu))_{i,j}^{k-1} & = & \left( (a_\rho^2a_\mu y_\mu Z \otimes X + a_\rho y_\rho X \otimes id)K_\alpha(\rho,\mu) \right)_{i,j}^{k-1} \nonumber \\
& = & a_\rho^2a_\mu y_\mu\omega^i K_\alpha(\rho,\mu)_{i,j-1}^{k-1} + a_\rho y_\rho K_\alpha(\rho,\mu)_{i-1,j}^{k-1}
\end{eqnarray}
\noindent Il est clair que $(\ref{sys})$ est sous-déterminé. On va spécifier les solutions cherchées en ajoutant les équations suivantes :
\begin{eqnarray}\label{eq1}
K_{\alpha+1}(\rho,\mu)_{i,j}^k = \lambda_{\alpha}(i,j\vert \rho,\mu) \ K_\alpha(\rho,\mu)_{i,j}^k\ ,\ \alpha = 0,\ldots,\ N-1,
\end{eqnarray} 

\noindent o\`u $\{ \lambda_{\alpha}(i,j\vert \rho,\mu) \}_{\alpha}$ est pour $i,\ j,\ \rho,\ \mu$ donnés une famille de nombres complexes non nuls. Ces équations peuvent être interprétées comme suit : si $U \in {\rm End}( M_{\rho\mu, \rho\otimes \mu})$, alors (\ref{sys}) montre que $U$ commute avec l'action de $\mathcal{W}_N$ sur $V_\rho \otimes V_\mu$, et (\ref{eq1}) peut s'écrire $UK_\alpha(\rho,\mu) = K_{\alpha+1}(\rho,\mu)$, o\`u $U_{i,j}^{k,l} = \lambda_{\alpha}(i,j\vert \rho,\mu) \ \delta(i-k) \ \delta(j-l)$ dans la base canonique de $\mc^N \otimes \mc^N$. 

\noindent D'après (\ref{eq1}), on doit avoir $\textstyle \prod_{\alpha = 0}^{N-1 } \lambda_{\alpha}(i,j\vert \rho,\mu) = 1$. Parmi les fonctions rationnelles en $a_\rho$, $y_\rho$, $a_\mu$ et $y_\mu$ qui vérifient cette égalité, les plus simples sont :
$$\lambda_{\alpha}(i,j\vert \rho,\mu) = \frac{z - x\omega^{g(i,j)-\alpha}}{y} \omega^{f(i,j,\alpha)}\ ,$$
\noindent o\`u $x,\ y$ et $z$ dépendent de $a_\rho$, $y_\rho$, $a_\mu$ et $y_\mu$ et v\'erifient $x^N + y^N = z^N$, et $g(i,j)$ et $f(i,j,\alpha)$ sont des fonctions affines que l'on va d\'eterminer. Notons qu'on doit avoir $\textstyle \sum_{l=0}^{N-1} f(i,j,l) = \frac{Nk}{2}$ avec $k \in \mathbb{Z}$. On peut donc considérer des solutions de la forme : 
$$K_\alpha(\rho,\mu)_{i,j}^k = \nu_{\rho,\mu} \ \omega^{\sum_{l=0}^{\alpha - 1}f(i,j,l)} \ \delta(i+j-k) \prod_{l = 1}^{g(i,j) - \alpha}\frac{y}{z - x\omega^l}\ .$$  
\noindent Le terme $\delta(i+j-k)$ est d\^u à (\ref{eq-1}). Quand $g(i,j) = i,\ x = a_\rho a_\mu y_\mu,\ y = y_\rho$, et $ z = a_\mu y_{\rho\mu}$, ces composantes matricielles sont solutions de (\ref{eq0}) si on impose aussi la condition : 
$$K_\alpha(\rho,\mu)_{i,j-1}^{k-1} = \omega^{-\alpha}K_\alpha(\rho,\mu)_{i,j}^k\ . \nonumber$$
\noindent Pour qu'elle soit vérifiée, on doit avoir
$$\sum_{l=0}^{\alpha - 1}f(i,j,l) = \alpha + \sum_{l=0}^{\alpha - 1}f(i,j-1,l)\ .$$
\noindent Négligeant l'indice $i$, supposons que $\textstyle \sum_{l=0}^{\alpha - 1}f(i,j,l) = j\alpha + f''(\alpha)$ pour une certaine fonction $f ''$. En écrivant $f(i,j,l) = l + 1/2 + j$, on trouve $f''(\alpha) = \alpha^2 /2$, d'o\`u les solutions de $(\ref{sys})$ décrites dans l'énoncé. 

\medskip

\noindent Consid\'erons maintenant les OCG duaux. On cherche une famille d'applications linéaires $\{\bar{K}^\alpha(\rho,\mu), \ \alpha = 0,\ldots,\ N-1\}$ v\'erifiant 
$$\sum_{k,\alpha = 0}^{N-1} h_{\rho,\mu} \ \omega^{\alpha j}\ \omega(a_\rho y_\mu,\frac{y_\rho}{a_\mu},y_{\rho\mu} \vert i,\alpha) \ \delta(i + j - k) \ \bar{K}^\alpha(\rho,\mu) _k^{i',j '} = \delta(i - i')\ \delta(j-j')\ .$$
\noindent Remplacons dans cette \'egalit\'e $\bar{K}^\alpha(\rho,\mu) _k^{i',j '}$ par les solutions de l'\'enonc\'e. On trouve
$$\begin{array}{l}
[\frac{a_\rho y_\mu}{y_{\rho\mu}}] \ \sum_{\alpha = 0}^{N-1} \omega^{\alpha (j-j')} \ \frac{\omega(a_\rho y_\mu,\frac{y_\rho}{a_\mu},y_{\rho\mu} \vert i - \alpha)}{\omega(\frac{a_\rho y_\mu}{\omega},\frac{y_\rho}{a_\mu},y_{\rho\mu} \vert i' - \alpha)} \ \delta(i + j - k)\ \delta(i' + j' - k)   \hspace{1cm}\\ \\
= [\frac{a_\rho y_\mu}{y_{\rho\mu}}] \  \delta(i + j - i' - j') \ \omega^{-i'(j' - j)} \ \omega(a_\rho y_\mu,\frac{y_\rho}{a_\mu},y_{\rho\mu} \vert i - i') \ \times \\
\hspace{4cm} \sum_{\alpha = 0}^{N-1} \omega^{(i' - \alpha)(j' - j)} \ \frac{\omega(a_\rho y_\mu \omega^{i-i'},\frac{y_\rho}{a_\mu},y_{\rho\mu} \vert i' - \alpha)}{\omega(\frac{a_\rho y_\mu}{\omega},\frac{y_\rho}{a_\mu},y_{\rho\mu} \vert i' - \alpha)} \\ \\
= [\frac{a_\rho y_\mu}{y_{\rho\mu}}] \ \delta(i + j - i' - j') \ \omega^{-i'(j' - j)} \ \omega(a_\rho y_\mu,\frac{y_\rho}{a_\mu},y_{\rho\mu} \vert i - i') \ \times\\
\hspace{7cm} f(\frac{a_\rho y_\mu\omega^{i - i'}}{y_{\rho\mu}},\frac{a_\rho y_\mu}{\omega y_{\rho\mu}} \vert \omega^{j ' - j})\ , 
\end{array}$$
\noindent o\`u l'on utilise (\ref{dec}) dans la seconde égalité. D'apr\'es la relation (A.14) de \cite[App. A]{KMS}, on a 
$$f(\frac{a_\rho y_\mu\omega^{i - i'}}{y_{\rho\mu}},\frac{a_\rho y_\mu}{y_{\rho\mu\omega}} \vert \omega^{j ' - j}) = \delta (i - i') \ [ \frac{y_{\rho\mu}}{a_\rho y_\mu} ]^{-1} \ (\frac{a_\rho y_\mu}{y_{\rho\mu}})^{-[j'-j-1]_N}\ ,$$
\noindent o\`u $[a]_N$ désigne le reste de la division euclidienne de $a \in \mathbb{N}$ par $N$. Or on vérifie immédiatement que 
$$[ \frac{a_\rho y_\mu}{y_{\rho\mu}} ] = (\frac{a_\rho y_\mu}{y_{\rho\mu}})^{N-1} \ [ \frac{y_{\rho\mu}}{a_\rho y_\mu} ] \ .$$ 
\noindent D'o\`u le r\'esultat.\hfill $\Box$

\bigskip

\begin{remark} {\rm On a vu que ${\rm dim}_{\mc}(M_{r,\rho \otimes \mu}) > 1$.  Ce fait est remarquable pour la construction d'invariants quantiques de noeuds ou de $3$-vari\'et\'es. En effet, le produit tensoriel $\rho \otimes \mu$ de deux représentations irréductibles de dimension finie de $U_h(sl(2,\mc))$, ou de $U_q(sl(2,\mc))$ lorsque $q$ est générique, se d\'ecompose de fa\c{c}on unique en la somme directe d'irréductibles de multiplicité $1$ chacun \cite[\S 7.7]{Kas} (``formule de Clebsch-Gordan''). Ces repr\'esentations sont utilis\'ees pour d\'efinir, par exemple, le polyn\^ome de Jones d'entrelacs dans $S^3$. De m\^eme, lorsque $q=\epsilon$ est une racine de l'unité, les représentations de dimension finie de $U_{\epsilon}^{res}(sl(2,\mc))$ (cf. e.g. \cite[\S 9.3]{CP}) ont une structure tensorielle très sophistiqu\'ee, mais les seules sous-représentations d'un produit tensoriel d'irréductibles qui interviennent dans la formule des invariants de Reshetikhin-Turaev de $3$-vari\'et\'es sont là aussi de multiplicité $1$ \cite[\S 11.2-3,15.3]{CP},\ \cite{KM}.}
\end{remark}

\medskip

\noindent On conclut cette section en montrant, à l'aide de la proposition \ref{CG1}, que les OCG sont des représentations de $R_{\omega} \in {\rm End}\left( Q(\widehat{\bar{\mathcal{H}}_{\omega^{-1}}})^{\otimes 2} \right)$ (cf. \S \ref{rationnalisation}). Soient $Z,\ X$ et $Y = \omega^{\frac{1}{2}}XZ$ les matrices de composantes
$$Z_{i,j} = \omega^i\ \delta (i-j),\ \ X_{i,j} = \delta (i-j-1),\ \ \mathrm{et} \ \ Y_{i,j} = \omega^{\frac{1}{2} +j}\ \delta (i-j-1)$$

\noindent dans la base canonique de $\mc^N$ ($\delta$ est le symbole de Kronecker mod($N$)). On v\'erifie facilement qu'on peut d\'efinir une  repr\'esentation $\rho$ de $\bar{\mathcal{H}}_{\omega^{-1}}$ par
$$\begin{array}{ll}
\rho(E) = a_\rho^{-2}\ Z^{-1},\ & \ \rho(D) = -\frac{y_\rho}{a_\rho}\ Y^{-1},\\
\rho(\bar{E}) = a_\rho^{-2}\ Y,\ & \rho(\bar{D}) = \frac{1}{a_\rho y_\rho}\ Z^{-1}Y = \frac{\omega^{-\frac{1}{2}}}{a_\rho y_\rho}\ X\ ,
\end{array}$$
\noindent o\`u $a_\rho$ et $y_\rho \in \mc^*$. (En particulier, $\rho$ respecte la relation $D\bar{D} - \bar{D}D = (1-\omega^{-1})E$ de $\mathcal{H}_{\omega^{-1}}$ - cf. Définition \ref{defdoubleq}). Clairement, $\rho$ est cyclique et irréductible; on dira que $\rho$ est une repr\'esentation \emph{normale}.  

\smallskip

\noindent Comme dans la définition \ref{repstandard} i) et la proposition \ref{prop1a} iii), on dit qu'une paire $(\rho,\mu)$ de représentations normales de $\bar{\mathcal{H}}_{\omega^{-1}}$ est \emph{régulière} si tous les $\mathcal{W}_N$-sous-modules simples de $\rho \otimes \mu$ sont cycliques. Soit $(\rho,\mu)$ une paire régulière de représentations normales de $\bar{\mathcal{H}}_{\omega^{-1}}$. Notons $\Upsilon,\ \Psi(\rho,\mu) \in {\rm End}(V_\rho \otimes V_\mu)$ les images par $\rho \otimes \mu$ de $\Upsilon_{\omega}$, $\Psi_{\omega}\in {\rm End}\left( Q(\widehat{\bar{\mathcal{H}}_{\omega^{-1}}})^{\otimes 2} \right)$, définis dans (\ref{Upsilon}) et (\ref{élémcanomega}). Le lemme suivant est une conséquence immédiate de (\ref{Upsilon}), (\ref{dilocycbis}) et de la définition des représentations normales de  $\bar{\mathcal{H}}_{\omega^{-1}}$. 

\begin{lem} \label{lemrep} Les endomorphismes $\Upsilon$ et $\Psi(\rho,\mu)$ ont les r\'ealisations matricielles suivantes dans les bases canoniques de $V_\rho$ et $V_\mu$ :
$$\begin{array}{l}
\Upsilon = \frac{1}{N}\ \sum_{i,j = 0}^{N-1} \omega^{ij}\ Z_1^{-i}\ Y_2^j = \sum_{i}^{N-1} Z_1^{-i}\ \widehat{Y}_2^i\ ,\\
\Psi(\rho,\mu) = h_{\rho,\mu} \sum_{t=0}^{N-1} \left( -\frac{y_\rho}{a_\rho a_\mu y_\mu}\ Y_1^{-1}Z_2^{-1}Y_2 \right)^t \ \prod_{s=1}^t \ \frac{1}{1 - \omega^{-s} \frac{y_{\rho\mu}}{a_\rho y_\mu}}\ .
\end{array}$$
\noindent o\`u $Z_1=Z \otimes id$, etc, et $\textstyle \widehat{Y}^i = \frac{1}{N}\sum_{k = 0}^{N-1}\omega^{ij} \ Y^j$.
\end{lem}

\medskip

\noindent On pr\'ecisera la valeur de $h_{\rho,\mu}$ dans le corollaire \ref{6j1}. Soit $(\rho,\mu)$ une paire régulière de représentations standards de $\mathcal {W}_N$. Consid\'erons les applications linéaires $K(\rho,\mu) :  M_{\rho\mu,\rho \otimes \mu} \otimes V_{\rho\mu} \rightarrow V_\rho \otimes V_\mu $ et $R(\rho,\mu) : V_\rho \otimes V_\mu \rightarrow M_{\rho\mu,\rho \otimes \mu} \otimes V_{\rho\mu}$ définies par $K(\rho,\mu)_{i,j}^{\alpha,k} = K_{\alpha}(\rho,\mu)_{i,j}^k$ et $R(\rho,\mu):= K^T(\rho,\mu)\ $ ($^T$ est la transposition).

\begin{prop}\label{CGdilocyc}  Les OCG sont des représentations de $R_{\omega} \in  {\rm End}\left( Q(\widehat{\bar{\mathcal{H}}_{\omega^{-1}}})^{\otimes 2} \right)$.
\end{prop}

\noindent {\it Démonstration} On va montrer que $R(\rho,\mu) = \Upsilon \cdot \Psi(\rho,\mu)$. Notons $\{ v_n \}_n$ la base canonique de $\mc^N$. On a $\Upsilon (v_k \otimes v_l) = v_k \otimes Y^k v_l$ et $\Upsilon^{-1} = v_k \otimes Y^{-k}v_l$. Or 
$$(Y^k)_{i,j} = \omega^{\frac{k^2}{2} + kj}\ \delta (i-j-k) \ \ \mathrm{et} \ \  (Y^{-k})_{i,j} = \omega^{-\frac{k^2}{2} - ki}\ \delta (i+k-j)\ .$$
\noindent Alors $\left( \Upsilon^{-1} \right)_{i,j}^{k,l} = \omega^{-\frac{k^2}{2}-kj}\ \delta (i-k) \delta (j+k-l)$, et
$$\begin{array}{lll}
\left( \Upsilon^{-1} \cdot R(\rho,\mu) \right)_{k,l}^{i,j} & = & h_{\rho,\mu} \sum_{m,n =0}^{N-1} \omega^{mj}\ \omega(a_\rho y_\mu,\frac{y_\rho}{a_\mu},y_{\rho\mu} \vert i,m) \ \delta (i+j-n) \ \times \\
& & \hspace{3cm}  \omega^{-\frac{m^2}{2}-ml}\ \delta (k-m) \ \delta (l+m-n) \\ \\
 & = & h_{\rho,\mu}\ \omega^{-k(i-k)} \  \omega(a_\rho y_\mu,\frac{y_\rho}{a_\mu},y_{\rho\mu} \vert i-k) \ \delta (i+j-l-k)\ .
\end{array}$$
\noindent D'autre part $-Y_1^{-1}Z_2^{-1}Y_2 = -\omega^{-1} (XZ)_1^{-1}X_2$ \ et $(XZ)_{i,j}^{-1} = \omega^{-i}\ \delta (i+1-j)$ donnent 
$$\left(  -\omega^{-1} (XZ)_1^{-1}X_2 \right)^t\ _{k,l}^{i,j} = (-1)^t \omega^{-kt-\frac{t(t+1)}{2}}\ \delta (l-j-t)\ \delta (k+t-i)\ .$$
\noindent Donc
$$\begin{array}{lll}
\Psi(\rho,\mu)_{k,l}^{i,j} & = & h_{\rho,\mu} \sum_{t=0}^{N-1} \left( -\frac{y_\rho}{a_\rho a_\mu y_\mu}\ Y_1^{-1}Z_2^{-1}Y_2 \right)^t \ \prod_{s=1}^t \ \frac{1}{1- \omega^{-s} \frac{y_{\rho\mu}}{a_\rho y_\mu}}\\ \\
 & = & h_{\rho,\mu} \sum_{t=0}^{N-1} \omega^{-kt-\frac{t(t+1)}{2}}\ \delta (l-j-t)\ \delta (k+t-i)\ \times \\
 & & \hspace{3.5cm} \left(-\frac{y_\rho}{a_\rho a_\mu y_\mu}\right)^t\ \ \prod_{s=1}^t \ \frac{1}{1 - \omega^{-s} \frac{y_{\rho\mu}}{a_\rho y_\mu}}\\ \\
 & = & h_{\rho,\mu} \ \omega^{-k(i-k)}\ \delta (l+k-i-j)\ \prod_{s=1}^{i-k} \ \frac{-\frac{y_\rho}{a_\mu}\omega^{-s}}{a_\rho y_\mu - \omega^{-s} y_{\rho\mu}}
\end{array}$$
$$\begin{array}{lll}
\hspace{1cm}  & = & h_{\rho,\mu} \ \omega^{-k(i-k)}\ \delta (l+k-i-j)\ \prod_{s=1}^{i-k} \ \frac{\frac{y_\rho}{a_\mu}}{ y_{\rho\mu} - a_\rho y_\mu \omega^{s}}\\ \\
 & = & h_{\rho,\mu} \ \omega^{-k(i-k)}\ \omega(a_\rho y_\mu,\frac{y_\rho}{a_\mu},y_{\rho\mu} \vert i-k)\ \delta (l+k-i-j)\\ \\
 & = & \left( \Upsilon^{-1} \cdot R(\rho,\mu) \right)_{k,l}^{i,j} \ .
\end{array}$$\hfill $\Box$

\subsection{Les $6j$-symboles et l'équation du pentagone} \label{symbvol}

\noindent Dans cette section on calcule les $6j$-symboles cycliques pour la base standard d' OCG. On en déduit que les $6j$-symboles cycliques sont des représentations de $R_{\omega} \in {\rm End}\left( Q(\widehat{\bar{\mathcal{H}}_{\omega^{-1}}})^{\otimes 2} \right)$.

\medskip 

\noindent Rappelons que $\mathcal{C}$ est l'ensemble des repr\'esentations cycliques irr\'eductibles de $\mathcal{W_N}$. Pour tout triplet r\'egulier $(\rho,\mu,\nu) \in \mathcal{C}^3$, notons 
$$\begin{array}{l}
M_{\rho,(\mu,\nu )} = {\rm End}_{\mathcal{W}_N} \left( V_{\rho \mu \nu},\left( V_\rho \otimes \left( V_\mu \otimes V_\nu \right)\right) \right) \\ \\ 
M_{(\rho,\mu),\nu} = {\rm End}_{\mathcal{W}_N} \left( V_{\rho \mu \nu},\left( \left( V_{\rho} \otimes V_\mu \right) \otimes V_\nu \right) \right)\ .
\end{array}$$
\noindent On a 
$$\begin{array}{l}
M_{\rho , (\mu ,\nu)} \cong  M_{\rho \mu \nu , \rho \otimes \mu \nu}\otimes  M_{\mu \nu , \mu \otimes \nu}\\ \\
M_{(\rho ,\mu), \nu} \cong   M_{\rho \mu  , \rho \otimes \mu} \otimes M_{\rho \mu \nu , \rho \mu \otimes \nu}\ .
\end{array}$$
\noindent L'isomorphisme de $\mathcal{W}_N$-modules $\left( V_\rho \otimes \left( V_\mu \otimes V_\nu \right)\right) \cong \left( \left( V_\rho \otimes V_\mu \right) \otimes V_\nu \right)$ induit un isomorphisme
\begin{equation} \label{iso6j}
R(\rho,\mu,\nu) : M_{\rho , (\mu , \nu)} \xrightarrow{\ \cong\ } M_{(\rho ,\mu ), \nu}
\end{equation}
\noindent qu'on appellera un $6j$-\emph{symbol}. Considérons le diagramme
$$\xymatrix{V_{\rho\mu \nu} \ar[r] \ar[d] & V_{\rho} \otimes V_{\mu \nu} \ar[d] \\ V_{\rho\mu} \otimes V_{\nu} \ar[r]  & V_{\rho} \otimes V_{\mu} \otimes V_{\nu}}$$
\noindent o\`u $(\rho,\mu,\nu) \in \mathcal{C}^3$ est régulière, et les flèches indiquent des plongements de représentations, i.e. des OCG. Par (\ref{iso6j}), ce diagramme est commutatif. Notons $\widetilde{K}(\rho,\mu) \in {\rm End}_{\mc}(M_{\rho\mu,\rho\otimes \mu})$ l'endomorphisme identit\'e, défini dans une  base $\{ K_{\alpha}(\rho,\mu)\}_{\alpha}$ de $M_{\rho\mu,\rho\otimes \mu}$ par 
$$\widetilde{K}(\rho,\mu)(\alpha) = K_{\alpha}(\rho,\mu)\ .$$
\noindent Posons $\widetilde{K}_1(\rho,\mu) = \widetilde{K}(\rho,\mu) \otimes id$ et $\widetilde{K}_2(\rho,\mu) = id \otimes \widetilde{K}(\rho,\mu)$. Alors la commutativit\'e du diagramme ci-dessus s'\'ecrit
$$\widetilde{K}_1(\rho,\mu) \ \widetilde{K}_2(\rho\mu,\nu) = R_{12}(\rho,\mu,\nu)\  \widetilde{K}_2(\mu,\nu) \ \widetilde{K}_1(\rho,\mu \nu)\ ,$$
\noindent ou encore
\begin{eqnarray}\label{6jdef1}
K_\alpha(\rho,\mu) \ K_\beta(\rho\mu,\nu) = \sum_{\delta,\gamma = 0}^{N-1} R(\rho,\mu,\nu)_{\alpha,\beta}^{\gamma,\delta} \ K_\delta(\mu,\nu) \ K_\gamma(\rho,\mu \nu)\ .
\end{eqnarray}
\noindent Ici, $R(\rho,\mu,\nu)$ appara\^it comme la matrice de changement de base entre les $N^2$ sous-modules simples de $V_\rho \otimes V_\mu \otimes V_\nu$. Ils sont tous isomorphes à $V_{\rho\mu \nu}$, mais plong\'es diff\'eremment par les OCG $K_\alpha(\rho,\mu)$.
\begin{lem} Pour tout suite r\'egulière $(\rho,\mu,\nu,\upsilon) \in \mathcal{C}^4$ on a \emph{l'\'equation du pentagone}
\begin{equation} \label{P1}
R_{12}(\rho,\mu,\nu)\ R_{13}(\rho,\mu \nu,\upsilon)\ R_{23}(\mu,\nu\upsilon,) = R_{23}(\rho\mu,\nu,\upsilon)\ R_{12}(\rho,\mu,\nu\upsilon)\ .
\end{equation}
\end{lem}

\noindent  La relation (\ref{P1}) traduit l'associativité du produit tensoriel pour \emph{quatre} représentations cycliques de $\mathcal{W}_N$, i.e. la coh\'erence de l'isomorphisme (\ref{iso6j}) \cite[\S 12]{Kas}, \cite[\S 5]{CP}.

\medskip

\noindent {\it D\'emonstration.} On a 
$$\begin{array}{l}
R_{12}(\rho,\mu,\mu)\ R_{13}(\rho,\mu \mu,\upsilon)\ R_{23}(\mu,\mu,\upsilon)\ \widetilde{K}_3(\mu,\upsilon) \ \widetilde{K}_2(\mu,\mu\upsilon)\ \widetilde{K}_1(\rho,\mu \mu\upsilon) \hspace{5cm}\\
 \hspace{3cm} = R_{12}(\rho,\mu,\mu)\ R_{13}(\rho,\mu \mu,\upsilon)\  \widetilde{K}_2(\mu,\mu)\ \widetilde{K}_3(\mu \mu,\upsilon)\ \widetilde{K}_1(\rho,\mu \mu\upsilon) \\
 \hspace{3cm}  =  R_{12}(\rho,\mu,\mu)\ \widetilde{K}_2(\mu,\mu)\ \widetilde{K}_1(\rho,\mu \mu)\ \widetilde{K}_3(\rho\mu \mu,\upsilon)\\
 \hspace{3cm} =  \widetilde{K}_1(\rho,\mu)\ \widetilde{K}_2(\rho\mu,\mu)\ \widetilde{K}_3(\rho\mu \mu,\upsilon)\ .  
\end{array}$$
\noindent Ensuite
$$\begin{array}{l}
\widetilde{K}_1(\rho,\mu)\ \widetilde{K}_2(\rho\mu,\mu)\ \widetilde{K}_3(\rho\mu \mu,\upsilon) \hspace{5cm} \\
\hspace{3cm} = \widetilde{K}_1(\rho,\mu)\ R_{23}(\rho\mu,\mu,\upsilon)\ \widetilde{K}_3(\mu,\upsilon) \ \widetilde{K}_2(\rho\mu,\mu\upsilon) \\
\hspace{3cm} =  R_{23}(\rho\mu,\mu,\upsilon)\ R_{12}(\rho,\mu,\mu\upsilon)\ \widetilde{K}_3(\mu,\upsilon)\ \widetilde{K}_2(\mu,\mu\upsilon)\ \widetilde{K}_1(\rho,\mu \mu\upsilon)\ .
\end{array}$$
\noindent Les morphismes $\widetilde{K}$ \'etant inversibles, on en d\'eduit le r\'esultat.\hfill $\Box$

\medskip

\noindent On va d\'ecrire explicitement les $6j$-symboles. Soient $\theta \in \mathbb{R}$ et 
$$\mathfrak{D}_{\theta} =\mc \setminus \{ \vert x \vert \geq 1, \ {\rm arg}(x) = \theta + 2k\pi/N,\ k=0,\ldots,\ N-1\}\ .$$
Posons
\begin{equation}\label{ghdel}
r(x) := (1 - x^N)^{1/N}\ ,\quad g(x) := \prod_{j=1}^{N-1} (1 - x\omega^j)^{j/N} \quad \mathrm{et} \ \ h(x) := x^{-P} \ \frac{g(x)}{g(1)}\ ,
\end{equation}
\noindent  o\`u $x \in \mc^*$ et $P = N-1/2$ (cf. \ref{PN}). Les fonctions $r$ et $g$ sont comprises comme le prolongement analytique à $\mathfrak{D}_{\theta}$ de leur développement en série pour $x=0$ (qui converge sur le disque unité ouvert de $\mc$). Montrons que $\vert g(1) \vert = N^{\frac{1}{2}}$, et donc que $h$ est bien définie. En effet :
$$g^N(1) = \prod_{j=1}^{N-1} (1 - \omega^{-j})^j = \prod_{j=1}^{N-1} (1 - \omega^{N-j})^j =  \prod_{l=1}^{N-1} (1 - \omega^{l})^{N-l}\ ,$$
\noindent Donc
$$\vert g^N(1) \vert^2= \prod_{k=1}^{N-1} (1 - \omega^{k})^{k+N-k} = \prod_{k=1}^{N-1} (1 - \omega^{k})^N = \left( \frac{1-x^N}{1-x} \right)_{\vert x=1} ^N = N^N\ .$$

\begin{conv} \label{conv1} {\rm Dans la suite, on suppose toujours implicitement que $\theta$ est tel que $g$ est \'evalu\'e \emph{hors} des coupures de $\mathfrak{D}_{\theta}$. De plus, on suppose toujours que les suites régulières de représentations standards $(\rho,\mu,\nu)$ sont telles que $-y_\rho y_\nu/y_{\rho\mu \nu} y_\mu = r\left( y_{\rho\mu} y_{\mu \nu}/y_{\rho\mu \nu} y_\mu \right)$ (cf. Proposition \ref{prop1a} i)).} 
\end{conv}

\noindent  Pour tout triplet r\'egulier $(\rho,\mu,\nu)$ de repr\'esentations standard, posons 
$$h_{\rho,\mu} = h\left( \frac{y_{\rho\mu}}{a_\rho y_\mu} \right)\ ,\quad h_{\rho,\mu,\nu} = h \left(\frac{y_{\rho\mu}y_{\mu \nu}}{y_{\rho\mu \nu}y_\mu}\right)\ .$$

\begin{prop} \label{6jdesc} Dans les bases standard d'OCG, les $6j$-symboles cycliques et leurs inverses ont les composantes matricielles suivantes : 
$$R(\rho,\mu,\nu)_{\alpha,\beta}^{\gamma,\delta} = h_{\rho,\mu,\nu} \ \omega^{\alpha\delta}\ \omega(y_{\rho\mu \nu}y_\mu,y_\rho y_\nu,y_{\rho\mu}y_{\mu \nu}\vert \gamma,\alpha) \ \delta(\gamma + \delta - \beta)$$
$$\bar{R}(\rho,\mu,\nu)_{\gamma,\delta}^{\alpha,\beta} =  \frac{\mathbb{[}\frac{y_{\rho\mu \nu}y_{\mu}}{y_{\rho\mu}y_{\mu \nu}}\mathbb{]}}{h_{\rho,\mu,\nu}}\ \omega^{-\alpha\delta}\ \frac{\delta(\gamma + \delta - \beta)}{\omega(\frac{y_{\rho\mu \nu}y_\mu}{\omega},y_\rho y_\nu,y_{\rho\mu}y_{\mu \nu}\vert \gamma,\alpha)}\ .$$
\end{prop}

\noindent La preuve montre en m\^eme temps que :

\begin{cor} \label{6j1} Les $6j$-symboles cycliques sont des représentations cycliques de $R_{\omega} \in  {\rm End}\left( Q(\widehat{\bar{\mathcal{H}}_{\omega^{-1}}})^{\otimes 2} \right)$.
\end{cor}

\noindent {\it Démonstration.} Réécrivons les OCG dans la relation (\ref{6jdef1}) à l'aide du corollaire \ref{CGdilocyc}. On rappelle que $R(\rho,\mu):= K^T(\rho,\mu)\ $, et on remarque que cette notation se prolonge naturellement aux $6j$-symboles $R(\rho,\mu,\nu)$. On obtient :
$$R_{23}(\rho\mu,\nu)\ R_{12}(\rho,\mu) = R_{12}(\rho,\mu,\nu)\ R_{13}(\rho,\mu \nu)\ R_{23}(\mu,\nu)\ ,$$
\noindent ou encore :
\begin{equation} \label{OKun}
\Upsilon_{23}\Psi_{23}(\rho\mu,\nu)\ \Upsilon_{12}\Psi_{12}(\rho,\mu) = \Upsilon_{12}\Psi_{12}(\rho,\mu,\nu)\ \Upsilon_{13}\Psi_{13}(\rho,\mu \nu)\ \Upsilon_{23}\Psi_{23}(\mu,\nu)\ .
\end{equation}
\noindent On va montrer que cette relation se factorise en l'équation du pentagone pour $\Upsilon$ et l'équation du dilogarithme cyclique pour $\Psi$. Posons 
$$U = -Y_1^{-1}Z_2^{-1}Y_2 = -\omega^{-1} (XZ)_1^{-1}X_2\ ,\quad V = -Y_2^{-1}Z_3^{-1}Y_3\ .$$
\noindent On a 
$$\Psi_{13}(\rho,\mu \nu) = \Psi_{13}(\rho,\mu \nu)(-\omega^{-1} (XZ)_1^{-1}X_3)$$
\noindent et on vérifie facilement que 
$$\begin{array}{l}
\Psi_{13}(\rho,\mu \nu)(-\omega^{-1} (XZ)_1^{-1}X_3)\ \Upsilon_{23} \hspace{8cm}\\ \\

\hspace {2cm} = \Upsilon_{23}\ \Psi_{13}(\rho,\mu \nu)(-\omega^{-1} (XZ)_1^{-1}Z_2^{-1}X_3) \\ \\
\hspace{2cm} = \Upsilon_{23}\ \Psi_{13}(\rho,\mu \nu)\left(-(-\omega^{-1} (XZ)_1^{-1}X_2)(-\omega^{-1} (XZ)_2^{-1}X_3)\right) \\ \\
\hspace {2cm} = \Upsilon_{23}\ \Psi_{13}(\rho,\mu \nu)(-UV)
\end{array}$$
$$\begin{array}{lll}
\Psi_{12}(\rho,\mu,\nu)(-\omega^{-1} (XZ)_1^{-1}X_2)\ \Upsilon_{13} & = & \Upsilon_{13}\ \Psi_{12}(\rho,\mu,\nu)(-\omega^{-1} (XZ)_1^{-1}X_2(XZ)_3)\\ \\
& = & \Upsilon_{13}\ \Psi_{12}(\rho,\mu,\nu)(U(XZ)_3)\ ,\\ \\
\Psi_{12}(\rho,\mu,\nu)(-\omega^{-1} (XZ)_1^{-1}X_2)\ \Upsilon_{23} & = & \Upsilon_{23}\ \Psi_{12}(\rho,\mu,\nu)(-\omega^{-1} (XZ)_1^{-1}X_2(XZ)_3^{-1}) \\ \\
& = & \Upsilon_{23}\ \Psi_{12}(\rho,\mu,\nu)(U(XZ)_3^{-1})\ .
\end{array}$$
\noindent Alors on peut transformer le membre de droite de (\ref{OKun}) comme suit (dans les dilogarithmes cycliques, on omet les arguments matriciels sans contribution significative) :
$$\begin{array}{l}
\Upsilon_{12}\Psi_{12}(\rho,\mu,\nu)\ \Upsilon_{13}\Psi_{13}(\rho,\mu \nu)\ \Upsilon_{23}\Psi_{23}(\mu,\nu)(V) \hspace{2cm} \\ \\
\hspace{2cm} = \Upsilon_{12}\Psi_{12}(\rho,\mu,\nu)\ \Upsilon_{13}\Upsilon_{23}\ \Psi_{13}(\rho,\mu \nu)(-UV)\ \Psi_{23}(\mu,\nu)(V) \\ \\ 
\hspace{2cm} =  \Upsilon_{12}\Upsilon_{13}\ \Psi_{12}(\rho,\mu,\nu)(U(XZ)_3)\Upsilon_{23}\ \Psi_{13}(\rho,\mu \nu)(-UV)\ \Psi_{23}(\mu,\nu)(V) \\ \\ 
\hspace{2cm} =\Upsilon_{12}\Upsilon_{13}\Upsilon_{23}\ \Psi_{12}(\rho,\mu,\nu)(U)\ \Psi_{13}(\rho,\mu \nu)(-UV)\ \Psi_{23}(\mu,\nu)(V)\ .
\end{array}$$
\noindent A gauche on a trivialement
$$\Upsilon_{23}\Psi_{23}(\rho\mu,\nu)\ \Upsilon_{12}\Psi_{12}(\rho,\mu) = \Upsilon_{23}\Upsilon_{12}\ \Psi_{23}(\rho\mu,\nu)\Psi_{12}(\rho,\mu)\ .$$
\noindent Or $\Upsilon$ est une solution de l'équation du pentagone (Lemme \ref{lemrep}). De plus, $\Psi$ est une solution de l'équation du dilogarithme cyclique
\begin{equation} \label{abon}
\Psi_{23}(\rho\mu,\nu)(V)\Psi_{12}(\rho,\mu)(U) = \Psi_{12}(\rho,\mu,\nu)(U)\Psi_{13}(\rho,\mu \nu)(-UV)\Psi_{23}(\mu,\nu)(V)
\end{equation}
\noindent lorsque ses paramètres vérifient (\ref{paramètres}). On en d\'eduit les paramètres de $R(\rho,\mu,\nu)=\Psi_{12}(\rho,\mu,\nu)(U)$. On détermine complètement $R(\rho,\mu,\nu)$ en imposant l'égalité des déterminants des deux membres de (\ref{abon}). Un calcul facile montre que cette égalité est équivalente à
$$\begin{array}{l}
 (\frac{a_\mu y_\nu}{y_{\mu \nu}})^{\frac{N(N-1)}{2}}\ h_{\rho,\mu}^Nh_{\rho\mu,\nu}^N\ \frac{g^{2N}(1)}{g^N(\frac{y_{\rho\mu}}{a_\rho y_\mu})g^N(\frac{y_{\rho\mu \nu}}{a_{\rho\mu}y_\nu})} = \hspace{6cm}\\ \\
\hspace{2.8cm}  (\frac{a_{\rho\mu}y_\nu}{y_{\rho\mu \nu}})^{\frac{N(N-1)}{2}}\ h_{\mu,\nu}^Nh_{\rho,\mu \nu}^Nh_{\rho,\mu,\nu}^N\ \frac{g^{3N}(1)}{g^N(\frac{y_{\mu \nu}}{a_\mu y_\nu})g^N(\frac{y_{\rho\mu \nu}}{a_{\rho}y_{\mu \nu}})g^N(\frac{y_{\rho\mu}y_{\mu \nu}}{y_{\rho\mu \nu}y_{\mu}})}\ .\end{array}$$
On en d\'eduit $h_{\rho,\mu}$ et $h_{\rho,\mu,\nu}$. La formule pour l'inverse des $6j$-symboles est une conséquence immédiate de celle pour les OCG (cf. Corollaire \ref{CG1}). \hfill $\Box$

\begin{remarks} \label{Tannaka}  {\rm 1) On peut interpreter la relation (\ref{6jdef1}) comme une \emph{équation du pentagone modifiée}, en termes de dualité Tannaka-Krein \cite[Prop. 7.1-7.2]{Dav}, \cite{Mil,CR}. Il serait naturel d'étendre cette dualité à notre contexte, o\`u l'on a des solutions de l'équation du pentagone \emph{à paramètres}.

\smallskip

\noindent 2) Si l'on applique $\log \circ \det$, pour une détermination fixée du logarithme, aux deux membres de l'équation (\ref{P1}), on obtient une relation de $3$-cocycle \emph{abélienne} sur un ouvert de Zariski du groupe $B$ (cf. (\ref{param})). D'autre part, les équivalences de quasi-bigèbres préservent la structure monoïdale de leurs catégories de représentations \cite[Th. 15.3.5]{Kas}, et induisent des \emph{transformations de jauge} sur leurs associateurs, i.e. leurs $6j$-symboles. Existe-t-il $C(\rho,\mu,\nu)$, obtenu à partir de $R(\rho,\mu,\nu)$ par une transformation de jauge, tel que $\log \circ \det(C(\rho,\mu,\nu))$ soit un $3$-cocycle sur $B$ ?}
\end{remarks}

\section{Les c-$6j$-symboles, leurs symmétries et l'équation du pentagone étendue} \label{etendu}

\noindent On construit dans cette section les c-$6j$-symboles, gr\^ace auxquels on d\'efinit dans \cite{BB1} les invariants quantiques hyperboliques de $3$-vari\'et\'es introduits.

\subsection{Calcul et unitarité des c-$6j$-symboles}

\noindent Les trois lemmes suivants sont utilisés à la proposition \ref{symmetry1} pour la construction des c-$6j$-symboles. Rappelons que les fonctions $f$, $\omega$ et $[\ ]$ ont été définies au \S \ref{symboles}, et les fonctions $g$ et $h$ dans (\ref{ghdel}).

\begin{lem} \emph{\cite[C.7]{KMS}} \label{sommation} On a la factorisation {\rm (cf. Convention \ref{conv1})} :

$$f(x,y \vert z) = (y\omega )^{P} \frac{g(1) g(y\omega /x) g(x/yz)}{g(1/x) g(y \omega) g(\omega /z)} \ ,$$

\vspace*{1mm}

\noindent o\`u $x$,$y$ et $z$ sont des nombres complexes non nuls vérifiant $z^N(1-y^N) = 1 - x^N$, et $P = (N-1)/2$.
\end{lem}

\begin{lem} \label{inversion} On a la formule d'inversion:
$$\forall \ k,l \in \mathbb{Z}\ ,\ \ \omega(x,y,z \vert k-l) = \frac{\omega^{kl} \ \omega^{\frac{-l^2}{2} - \frac{k^2}{2}}}{\omega(z,-\omega^{1/2}y,\omega x \vert l-k)}\ ,$$

\noindent o\`u $x$,$y$ et $z$ sont des nombres complexes non nuls vérifiant $x^N + y^N = z^N$.\hfill $\Box$
\end{lem}

\noindent {\it Démonstration.} Rappelons qu'on a choisi $\omega^{\frac{1}{2}}$ tel que $\omega^{\frac{N}{2}} = 1$. On a :
$$\prod_{j=1}^{k-l}\frac{y}{z-x\omega^j} \ \prod_{j=1}^{N+l-k}\frac{-\omega^{\frac{1}{2}}y}{\omega x-z\omega^j} = \frac{y^N\ \omega^{\frac{l-k}{2}}}{\prod_{j=1}^{k-l}(z-x\omega^j) \ \prod_{j=1}^{N+l-k} (z\omega^j - \omega x)}\ ,$$
\noindent et
$$\begin{array}{lll}
\prod_{j=1}^{N+l-k} (z\omega^j - \omega x) & = & \omega^{\frac{(N+l-k)(N+l-k+1)}{2}}\ \prod_{j=1}^{N+l-k} (z - x\omega^{1-j}) \\
                                           & = & \omega^{\frac{(l-k)^2}{2} + \frac{(l-k)}{2}} \ \prod_{n=1+k-l}^{N} (z - x\omega^{n})\ .
\end{array}$$
\noindent Donc  
$$\begin{array}{lll}
\prod_{j=1}^{k-l}\frac{y}{z-x\omega^j} \ \prod_{j=1}^{N+l-k}\frac{-\omega^{\frac{1}{2}}y}{\omega x-z\omega^j} & = & \frac{y^N}{\omega^{\frac{(l-k)^2}{2}}\ \prod_{j=1}^{k-l}(z-x\omega^j) \ \prod_{n=1+k-l}^N (z - x\omega^n)} \\ \\
     & = & \omega^{\frac{(l-k)^2}{2}} = \omega^{kl} \ \omega^{\frac{-l^2}{2} - \frac{k^2}{2}}\ .
\end{array}$$ \hfill $\Box$

\begin{lem}\label{lem1}
Pout tout $x \in \mc^*$, on a l'identité :
$$[x] = h(1/x)h(x)x^P \ .$$
\end{lem}

\noindent {\it Démonstration.} Le lemme \ref{sommation} implique
$$g(1/x)g(x)f(x,\frac{x}{\omega} \vert \omega) = x^Pg^2(1)\ .$$
\noindent Or
$$\begin{array}{lll}
f(x,\frac{x}{\omega} \vert \omega) & = & \sum_{l=0}^{N-1} \frac{\omega(x\vert l)}{\omega(\frac{x}{\omega}\vert l)}\ \omega^l = 1 + \frac{(1-x)\omega}{1-x\omega} + \frac{(1-x)\omega^2}{1-x\omega^2} + \ldots + \frac{(1-x)\omega^{N-1}}{1-x\omega^{N-1}} \\  
& & \\
& = & (1-x) \left( \frac{1}{1-x} + \frac{\omega}{1-x\omega} + \ldots + \frac{\omega^{N-1}}{1-x\omega^{N-1}} \right) \\
& & \\
& = & (1-x) \ \frac{d}{dx}\left( - \log (1-x^N) \right) \\ \\
& = & N x^{N-1}\frac{1-x}{1-x^N}\ .
\end{array}$$
\noindent Donc
$$\frac{[x]}{h(1/x)} = N^{-1}\frac{1-x^N}{1-x}x^{-P}\frac{g(1)}{g(1/x)} = \frac{g(x)}{g(1)} =  h(x)x^P\ .$$
\hfill $\Box$

\medskip

\noindent {\bf T\'etraèdres d\'ecor\'es.} Consid\'erons $\mr^3$ muni d'une orientation quelconque, et soit $\Delta$ un t\'etraèdre plong\'e dans $\mr^3$. Un \emph{branchement} $b$ de $\Delta$ est la donn\'ee d'une orientation de ses ar\^etes, telle que sur chaque face deux ar\^etes seulement ont une orientation cons\'ecutive. Un branchement d\'etermine un ordre sur les sommets de $\Delta$, et une orientation des faces. Notons $v_0$, $v_1$, $v_2$ et $v_3$ les sommets, $\mathcal{E}$ les ar\^etes orient\'ees de $\Delta$, $\mathcal{F}$ ses faces, et $f_j$ la face oppos\'ee à $v_j$. Soient $e_0$, $e_1$, $e_2 \in \mathcal{E}$ les ar\^etes $e_0 = [v_0,v_1]$, $e_1=[v_1,v_2]$ et $e_2=[v_0,v_2] =-[v_2,v_0]$, et $e_0'$, $e_1'$ et $e_2'$ les ar\^etes oppos\'ees. Le trièdre $(e_0',e_1',e_2')$ d\'etermine une orientation $*=\pm$ de $\Delta$, \emph{positive} si elle est coh\'erente avec celle induite par l'orientation de $\mr^3$, et n\'egative sinon.

\noindent Fixons un branchement $b$ sur $\Delta$. Rappelons que $B$ est le sous-groupe de Borel de $SL(2,\mc)$ des matrices triangulaires sup\'erieures (cf. (\ref{param})). Consid\'erons l'ensemble des $1$-cocycles pleins $z \in Z^1(\Delta;B)$ sur $\Delta$ à valeurs dans $B$, o\`u ``plein'' signifie qu'\emph{aucune} valeur de $z$ n'est une matrice diagonale. On utilise la $b$-orientation des ar\^etes de $\Delta$ pour d\'efinir $z$, et on pose $z(-e)=z(e)^{-1}$. Notons $z(e) = (t(e),x(e))$ pour
$$z(e) = \left( \begin{array}{ll} 
                t(e) & x(e) \\
                0     & t(e)^{-1}
                \end{array} \right)\ .$$
\noindent Fixons une d\'etermination de la racine $N$-ième commune à toutes les composantes matricielles des matrices $z(e)$. Consid\'erons les repr\'esentations standard $r_N(e)$ de $\mathcal{W}_N$, $e \in \mathcal{E}$, telles que $(a_{r_N(e)}, y_{r_N(e)}) = (t(e)^{1/N},x(e)^{1/N})$ (rappelons que $z \in Z^1(\Delta;B)$ est plein). Alors la relation (\ref{multiplication}) traduit exactement la propri\'et\'e de cocycle de $z$.

\noindent Un \emph{\'etat} de $\Delta$ est une application $\alpha: \mathcal{F} \rightarrow \mz/N\mz$; notons $\alpha_{j} = \alpha(f_j)$. A tout 4-uplet $*(\Delta,b,z,\alpha)$, on peut associer 
\begin{eqnarray}
\Xi(*(\Delta,b,r_N,\alpha)) = \left\lbrace 
\begin{array}{l} \label{assoc6j}
R(r_{N}(e_0),r_{N}(e_1),r_{N}(e_0'))_{\alpha_{3},\alpha_{1}}^{\alpha_{2},\alpha_{0}}\  \ \rm{si} \ *=+1\\ \\
\bar{R}(r_{N}(e_0),r_{N}(e_1),r_{N}(e_0'))_{\alpha_{2},\alpha_{0}}^{\alpha_{3},\alpha_{1}} \  \ \rm{si} \ *=-1\ .
\end{array} \right.
\end{eqnarray} 
\noindent Les $6j$-symboles cycliques peuvent donc \^etre vus comme des op\'erateurs sur les t\'etraèdres \emph{d\'ecor\'es} $(\Delta,b,z)$. Il n'est pas encore clair comment d\'efinir des invariants quantiques de $3$-vari\'et\'es à partir des $6$-symboles cycliques, en utilisant des triangulations munies de d\'ecorations g\'en\'eralisant $b$ et $z$. Le problème est le suivant : l'\'equation du pentagone, consid\'er\'ee comme identit\'e entre op\'erateurs d\'efinis sur un polyèdre d\'ecor\'e, n'est v\'erifi\'ee que pour \emph{certains} branchement du polyèdre\footnote{Ce problème très subtil fait l'objet d'un travail en cours, en collaboration avec R. Benedetti et F. Costantino.}(dont celui correspondant à (\ref{P1})) \cite[Ch. 2]{BB1}. 

\smallskip

\noindent On construit donc les c-$6j$-symboles. Ce sont des symm\'etrisations \emph{partielles} des $6j$-symboles, au sens o\`u il existe une action projective non triviale de $SL(2,\mz)$ sur les c-$6j$-symboles (cf. Remarques \ref{matrices1}). Notons que tous les $6j$-symboles utilis\'es pour construire des invariants quantiques ``à la Turaev-Viro'' ont des symm\'etries t\'etra\'edrales, ou sont symm\'etris\'es de fa\c{c}on \'el\'ementaire \cite{TV,T}, \cite{BW,GK}.   

\medskip

\noindent {\bf Construction des c-$6j$-symboles.} Soit $(\Delta,b,z)$ comme ci-dessus. Consid\'erons son $6j$-symbole (\ref{assoc6j}); notons $r_{N}(e_0)=\rho$, $r_{N}(e_1)=\mu$, et $r_{N}(e_0')=\nu$. En permutant les sommets $v_0$ et $v_1$, on obtient 
$$\bar{R}(\bar{\rho},\rho\mu,\nu)_{\alpha_2,\alpha_1}^{\alpha_3,\alpha_0} = \frac{[\frac{x_{\rho\mu \nu}x_{\mu}}{x_{\rho\mu}x_{\mu \nu}}]}{h(\frac{x_{\rho\mu \nu}x_\mu}{x_{\rho\mu}x_{\mu \nu}})}\ \omega^{-\alpha_3\alpha_1}\frac{\delta(\alpha_2 + \alpha_1 - \alpha_0)}{\omega(\frac{x_{\rho\mu}x_{\mu \nu}}{\omega},x_{\bar{\rho}}x_\nu,x_{\mu}x_{\rho\mu \nu}\vert \alpha_2,\alpha_3)}\ .$$

\pagebreak

\noindent On d\'eduit des lemmes \ref{inversion} et \ref{lem1} que
$$\begin{array}{lll}
\bar{R}(\bar{\rho},\rho\mu,\nu)_{\alpha_2,\alpha_1}^{\alpha_3,\alpha_0} & = & \left(\frac{x_{\rho\mu}x_{\mu \nu}}{x_{\rho\mu \nu}x_{\mu}}\right)^P h_{\rho,\mu,\nu} \ \omega^{-\alpha_3 \alpha_0 + \frac{(\alpha_3 - \alpha_2)}{2}} \ \times \\
 &  & \hspace{2cm} \omega(x_{\rho\mu \nu}x_\mu,x_\rho x_\nu,x_{\rho\mu}x_{\mu \nu}\vert \alpha_3,\alpha_2) \ \delta(\alpha_2 + \alpha_1 - \alpha_0) \\  \\
 & = & \left(\frac{x_{\rho\mu}x_{\mu \nu}}{x_{\rho\mu \nu}x_{\mu}}\right)^P R(\rho,\mu,\nu)_{-\alpha_3,\alpha_1}^{-\alpha_2,\alpha_0}\ \omega^{\frac{(\alpha_3 - \alpha_2)}{2} + \frac{(\alpha_2)^2}{2} - \frac{(\alpha_3)^2}{2}}.
\end{array}$$
\noindent Soient $T = \{ T_{m,n} \},\ T^{-1} = \{ T^{m,n} \},\ S = \{ S_{m,n} \}$ et $S^{-1} = \{ S^{m,n} \}$ les matrices de taille $N \times N$ définies par :
$$\begin{array}{ll}
T_{m,n} = \zeta^{-1}\omega^{\frac{m^2}{2}}\delta(m + n),  & \ S_{m,n} = N^{-\frac{1}{2}}\omega^{mn}, \\
T^{m,n} = \zeta\omega^{-\frac{m^2}{2}}\delta(m + n),  & \ S^{m,n} = N^{-\frac{1}{2}}\omega^{-mn}\ ,
\end{array}$$
\noindent o\`u $\zeta \in \mathbb{C}^*$. On trouve finalement
$$\bar{R}(\bar{\rho},\rho\mu,\nu)_{\alpha_2,\alpha_1}^{\alpha_3,\alpha_0} =  \left(\frac{x_{\rho\mu}x_{\mu \nu}}{x_{\rho\mu \nu}x_{\mu}}\right)^P \omega^{\frac{(\alpha_3 - \alpha_2)}{2}} \sum_{\alpha',\gamma'=0}^{N-1} R(\rho,\mu,\nu)_{\alpha',\alpha_1}^{\gamma',\alpha_0} \ T_{\gamma',\alpha_2} T^{\alpha ',\alpha_3}\ .$$
\noindent Le même type de calculs pour les transpositions $(v_2,v_3)$ et $(v_1,v_2)$ donnent respectivement
$$\begin{array}{l}
\bar{R}(\rho,\mu \nu,\bar{\nu})_{\alpha_3,\alpha_0}^{\alpha_2,\alpha_1} =  \left(\frac{x_{\rho\mu \nu}x_{\mu}}{x_{\rho\mu}x_{\mu \nu}}\right)^P \omega^{\frac{(\alpha_2 - \alpha_3)}{2}} \sum_{\beta',\delta'=0}^{N-1} R(\rho,\mu,\nu)_{\alpha_3,\beta '}^{\alpha_2,\delta '} \ S_{\alpha_0,\delta '} S^{\alpha_1,\beta '}, \\
\bar{R}(\rho\mu,\bar{\mu},\mu \nu)_{\alpha_1,\alpha_0}^{\alpha_3,\alpha_2} =  \left(\frac{x_{\rho}x_{\nu}}{x_{\rho\mu}x_{\mu \nu}}\right)^P \sum_{\alpha',\delta'=0}^{N-1} R(\rho,\mu,\nu)_{\alpha',\alpha_1}^{\alpha_2,\delta '} \ T_{\alpha_0,\delta'} S^{\alpha_3,\alpha'}
\end{array}$$
\noindent pour une certaine racine de l'unité $\zeta$.  

\begin{defi} \label{c6jnouv}
Soient $a,c \in \mathbb{Z}_N$ et $(\rho,\mu,\nu)$ une suite régulière de représentations standard de $\mathcal{W}_N$. Les \emph{c-$6j$-symboles} et leurs inverses sont définis par 
$$\begin{array}{l}
R(\rho,\mu,\nu\vert a,c)_{\alpha,\beta}^{\gamma,\delta} = (x_{\rho\mu}x_{\mu \nu})^P\ \omega^{c(\gamma - \alpha) - \frac{ac}{2}}\ R(\rho,\mu,\nu)_{\alpha,\beta - a}^{\gamma - a ,\delta}\ ,\\ 
  \\
\bar{R}(\rho,\mu,\nu\vert a,c)_{\gamma,\delta}^{\alpha,\beta} = (x_{\rho\mu}x_{\mu \nu})^P\ \omega^{c(\gamma - \alpha) + \frac{ac}{2}}\ \bar{R}(\rho,\mu,\nu)_{\gamma+a,\delta}^{\alpha,\beta+a}\ .
\end{array}$$
\end{defi}

\begin{prop} \label{symmetry1} Les c-$6j$-symboles vérifient les relations de symmetrie suivantes :
$$\begin{array}{l}
\sum_{\alpha ',\gamma '=0}^{N-1} R(\rho,\mu,\nu\vert a,c)_{\alpha ',\beta}^{\gamma ',\delta} \ T_{\gamma,\gamma '} T^{\alpha,\alpha'} =\omega^{a/4}\ \bar{R}(\bar{\rho},\rho\mu,\nu\vert a,b)_{\gamma,\beta}^{\alpha,\delta}\ ,\\

 \\

\sum_{\alpha ',\delta '=0}^{N-1} R(\rho,\mu,\nu\vert a,c)_{\alpha ',\beta}^{\gamma,\delta'} \ T_{\delta,\delta '} S^{\alpha,\alpha'} =\omega^{-c/4}\ \bar{R}(\rho\mu,\bar{\mu},\mu \nu\vert b,c)_{\beta,\delta}^{\alpha,\gamma}\ ,\\

\\

\sum_{\beta ',\delta '=0}^{N-1} R(\rho,\mu,\nu\vert a,c)_{\alpha,\beta '}^{\gamma,\delta'} \ S_{\delta,\delta '} S^{\beta,\beta '} =\omega^{a/4}\ \bar{R}(\rho,\mu \nu,\bar{\nu}\vert a,b)_{\alpha,\delta}^{\gamma,\beta}\ ,
\end{array}$$
\noindent o\`u $b = P+1 -a -c \in \mathbb{Z}_N$, et on prend $\zeta = \omega^{\frac{9}{8}}\ (-1)^P \ \vert g(1) \vert/g(1)$. 
\end{prop}
\noindent Une \emph{charge int\'egrale} de $\Delta$ est une application $c: \mathcal{E} \rightarrow \mz$ telle que $c(e)=c(-e)$, et si $e$, $e'$ et $e''$ sont dans une m\^eme face de $\Delta$, alors $c(e) + c(e') + c(e'') = 1$. Posons $c_N(e) = c(e)/2$ mod($N$). Comme dans (\ref{assoc6j}), pour tout 5-uplet $(\Delta,b,z,c,\alpha)$, on pose $\Xi(*(\Delta,b,r_N,c_N,\alpha))$ \'egal à
\begin{eqnarray} \label{assocc6j}
\left\lbrace 
\begin{array}{l}
R(r_{N,i}(e_0),r_{N,i}(e_1),r_{N,i}(e_0') \vert c_N(e_0),c_N(e_1))_{\alpha_{i3},\alpha_{i1}}^{\alpha_{i2},
\alpha_{i0}}\  \ \rm{if} \ *=+1\\ \\
\bar{R}(r_{N,i}(e_0),r_{N,i}(e_1),r_{N,i}(e_0') \vert c_N(e_0),c_N(e_1))_{\alpha_{i2},\alpha_{i0}}^
{\alpha_{i3},\alpha_{i1}} \  \ \rm{if} \ *=-1\ .
\end{array} \right.
\end{eqnarray}
\begin{remarks} \label{matrices1} {\rm 1) Les matrices $T$ et $T^{-1}$ (resp. $S$ et $S^{-1}$) sont inverses l'une de l'autre, et on v\'erifie facilement que $S^4 = id$ et $S^2 = \zeta'(ST)^3$, avec $\vert \zeta' \vert = 1$. Par conséquent, les matrices $S$ et $T$ définissent une représentation projective de dimension $N$ de $SL(2,\mathbb{Z})$, qui admet une pr\'esentation de la forme $\langle s,t\vert s^4=1,s^2=(st)^3\rangle$.

\smallskip

\noindent 2) Notons $\widehat{\mathcal{W}}_N$ l'extension \emph{affine} de $\mathcal{W}_N$. Il existe une action projective bien connue de $SL(2,\mz)$ sur $\widehat{sl}(2,\mc)$ \cite[\S 4]{Wa}. Rappelons que $\widehat{\mathcal{W}}_N$ peut être réalisée comme sous-algèbre de Hopf de $U_{\omega}(\widehat{sl}(2,\mc))$. Les c-$6j$-symboles sont-ils (à $(x_{\rho\mu}x_{\mu \nu})^P$ près) des repr\'esentations cycliques du ddH de $\widehat{\mathcal{W}}_N$ ?  Les valeurs de $c$ dans (\ref{assocc6j}) fixeraient ainsi le \emph{niveau} des représentations irréductibles et cycliques de $\widehat{\mathcal{W}}_N$. }
\end{remarks}

\smallskip

\noindent {\it Démonstration de la proposition \ref{symmetry1}.} On obtient la première et la troisième relation de symm\'etrie en reprenant les calculs qui pr\'ecèdent la d\'efinition \ref{c6jnouv}. La seule diff\'erence est qu'on utilise la relation $b = P+1 -a -c$ mod($N$). La seconde relation de symmétrie est plus difficile à vérifier. Dans (\ref{assocc6j}), notons $r_{N}(e_0)=\rho$, $r_{N}(e_1)=\mu$, $r_{N}(e_0')=\nu$, et $c_N(e_0)=a$, $c_N(e_1)=c$ et $c_N(e_2)=b$. La permutation $(v_1,v_2)$ change $R(\rho,\mu,\nu\vert a,c)_{\alpha,\beta}^{\gamma,\delta}$ en 
$$\bar{R}(\rho\mu,\bar{\mu},\mu \nu\vert b,c)_{\beta,\delta}^{\alpha,\gamma} = \frac{(x_\rho x_\nu)^P[\frac{x_{\rho\mu \nu}x_{\bar{\mu}}}{x_\rho x_\nu}]}{h(\frac{x_\rho x_\nu}{x_{\rho\mu \nu}x_{\bar{\mu}}})} \ \frac{\delta(\beta + \delta - \gamma)\ \omega^{c(\beta - \alpha) + bc/2 -\alpha \delta}}{\omega(\frac{x_{\rho\mu \nu}x_{\bar{\mu}}}{\omega},x_{\rho\mu}x_{\mu \nu},x_\rho x_\nu \vert \beta + b,\alpha)}\ .$$
\noindent Par un calcul direct, on vérifie que
$$\forall x \in \mathbb{C}^*,\ x \ne \omega^j \ (j=1,\ldots,\ N),\ \ \frac{[x^{-1}]}{h(x)} = N^{-1}\frac{g(1)}{g(x)\frac{1-x}{r(x)}} \ x^{-P}r(x)^{N-1}\ ,$$
\noindent o\`u $r(x) = (1-x^N)^{1/N}$ (cf. Convention \ref{conv1}). Comme 
$$r(x_\rho x_\nu/x_{\rho\mu \nu}x_{\bar{\mu}}) = - x_{\rho\mu}x_{\mu \nu}/x_{\rho\mu \nu}x_{\bar{\mu}}\ ,$$ 
\noindent on déduit de cette formule, du lemme \ref{gdilocyc} et de la définition de $[ \ ]$ et de $h$ que 
$$\frac{(x_\rho x_\nu)^P[\frac{x_{\rho\mu \nu}x_{\bar{\mu}}}{x_\rho x_\nu}]}{h(\frac{x_\rho x_\nu}{x_{\rho\mu \nu}x_{\bar{\mu}}})} =  N^{-1}\frac{g(1)\omega^{-1}}{g(\frac{x_\rho x_\nu\omega}{x_{\rho\mu \nu}x_{\bar{\mu}}})} \ \left(-\frac{x_{\rho\mu}x_{\mu \nu}}{x_{\rho\mu \nu}x_{\bar{\mu}}}\right)^{2P} (x_{\rho\mu \nu}x_{\bar{\mu}})^P\ .$$
\noindent D'autre part, (\ref{ompetit}) donne imm\'ediatement
$$\omega(x,y,z \vert n) = \left(\frac{y}{z}\right)^n \  \omega \left(\frac{x}{z} \vert n\right)\ .$$ 
\noindent On peut donc écrire :
$$\begin{array}{l}
\bar{R}(\rho\mu,\bar{\mu},\mu \nu\vert b,c)_{\beta,\delta}^{\alpha,\gamma} = 
N^{-1}\omega^{-1}g(1)\ \frac{(x_{\rho\mu}x_{\mu \nu})^{2P}}{(x_{\rho\mu \nu}x_{\bar{\mu}})^P} \left(\frac{x_\rho x_\nu}{x_{\rho\mu}x_{\mu \nu}}\right)^{\beta + b -\alpha}\hspace{3cm}  \\ \\ 
\hspace{6.5cm} \times \ \frac{\delta(\beta + \delta - \gamma)\ \omega^{c(\beta - \alpha) + \frac{bc}{2} -\alpha \delta - \frac{\alpha^2}{2}}}{g(\frac{x_\rho x_\nu\omega}{x_{\rho\mu \nu}x_{\bar{\mu}}})\ \omega(\frac{x_{\rho\mu \nu}x_{\bar{\mu}}}{x_\rho x_\nu\omega} \vert \beta + b - \alpha)}\ .
\end{array}$$
\noindent Considérons la fonction rationnelle $f(x,y \vert z)$ définie dans (\ref{sommf}). Par le lemme \ref{sommation} on a
$$f(0,y \vert z\omega^m) = (y\omega)^{2P} \frac{g(1)}{g(y\omega)g(\omega /z\omega^m)}\ .$$
\noindent Appliquons le lemme \ref{gdilocyc}  à $g(\omega /z\omega^m)$ et $g(\omega /z)$. On trouve
$$f(0,y \vert z\omega^m) = f(0,y \vert z) \prod_{l=1}^m \frac{1-z\omega^{l-1}}{-yz\omega^l} = (y\omega)^{2P} \ \ \frac{g(1)\ (-yz)^{-m}\ \omega^{-m(m+1)/2}}{g(y\omega)g(\omega /z)\ \omega(\frac{z}{\omega} \vert m)}\ .$$
\noindent Fixons 
$$z = \frac{x_{\rho\mu \nu}x_{\bar{\mu}}}{x_\rho x_\nu},\quad y= \frac{x_{\rho\mu}x_{\mu \nu}}{x_{\rho\mu \nu}x_\mu\omega},\quad m= \beta + b -\alpha\ .$$
\noindent Les variables $y$ et $z$ sont toutes les deux dans le domaine de définition de $f(0,y \vert z)$ et $f(0,y \vert z\omega^m)$. Par substitutions successives, on peut transformer $\bar{R}(\rho\mu,\bar{\mu},\mu \nu\vert b,c)_{\beta,\delta}^{\alpha,\gamma}$ en
$$\begin{array}{l}
(x_{\rho\mu \nu}x_\mu)^PN^{-1}\omega^{-1}(-1)^Pg(\frac{x_{\rho\mu}x_{\mu \nu}}{x_{\rho\mu \nu}x_\mu}) \ \delta(\beta + \delta - \gamma)\  \times \\ \hspace{2cm} \omega^{c(\beta - \alpha) + \frac{bc}{2} -\alpha \delta - \frac{\alpha^2}{2} + \frac{(\beta + b - \alpha)(\beta + b - \alpha-1)}{2}}\ f(0,\frac{x_{\rho\mu}x_{\mu \nu}}{x_{\rho\mu \nu}x_\mu\omega} \vert \frac{x_{\rho\mu \nu}x_{\bar{\mu}}\omega^{\beta + b -\alpha}}{x_\rho x_\nu})\ .
\end{array}$$ 
\noindent Maintenant, on a la chaîne de transformations : 
$$\begin{array}{l}
f(0,\frac{x_{\rho\mu}x_{\mu \nu}}{x_{\rho\mu \nu}x_\mu\omega} \vert \frac{x_{\rho\mu \nu}x_{\bar{\mu}}\omega^{\beta + b -\alpha}}{x_\rho x_\nu})  \\ \\
\hspace{2cm} = \sum_{\sigma =0}^{N-1} \frac{\left( \frac{x_{\rho\mu \nu}x_{\bar{\mu}}\omega^{\beta + b -\alpha}}{x_\rho x_\nu} \right)^{\sigma}}{\omega\left(\frac{x_{\rho\mu}x_{\mu \nu}}{x_{\rho\mu \nu}x_\mu\omega} \vert \sigma \right)} =  \sum_{\sigma =0}^{N-1} \frac{\omega^{\sigma (b+\beta - \alpha - 1/2)}\ (-1)^{\sigma}}{\omega (x_{\rho\mu}x_{\mu \nu},\omega^{1/2}x_\rho x_\nu,\omega x_{\rho\mu \nu}x_\mu \vert \sigma)} \\ \\
\hspace{2cm} =\sum_{\sigma =0}^{N-1} \omega^{\sigma (b+\beta - \alpha - 1/2)}\ \omega^{\frac{\sigma^2}{2}}\ \omega (x_{\rho\mu \nu}x_\mu,x_\rho x_\nu,x_{\rho\mu}x_{\mu \nu}\vert -\sigma) \\ \\
\hspace{2cm} =\sum_{\sigma = 0}^{N-1} \omega^{-\frac{\sigma}{2} + \frac{(\sigma + b + \beta - \alpha)^2}{2} - \frac{(b + \beta - \alpha)^2}{2}} \omega(x_{\rho\mu \nu}x_\mu,x_\rho x_\nu,x_{\rho\mu}x_{\mu \nu} \vert -\sigma )\ .
\end{array}$$
\noindent Dans l'avant-dernière égalité on a utilisé le lemme \ref{inversion}. En considérant le changement de variable $-\sigma = \gamma - a - \alpha'$ et en simplifiant (un peu) les puissances de $\omega$ à l'aide de $\delta(\beta + \delta - \gamma)$, on trouve alors facilement
$$\begin{array}{l}
\bar{R}(\rho\mu,\bar{\mu},\mu \nu\vert b,c)_{\beta,\delta}^{\alpha,\gamma} = \\ \\
\hspace{0.2cm} (x_{\rho\mu}x_{\mu \nu})^P h(\frac{x_{\rho\mu}x_{\mu \nu}}{x_{\rho\mu \nu}x_\mu})\ N^{-1}\omega^{-1}(-1)^Pg(1)\ \omega^{c(\beta - \alpha) + \frac{bc}{2} -\alpha \delta - \frac{\alpha^2}{2} - \frac{(\beta + b - \alpha)}{2}} \ \times \\ \\
\hspace{0.2cm} \delta(\beta + \delta - \gamma)\ \sum_{\sigma = 0}^{N-1} \omega^{\frac{(\gamma - a - \alpha')}{2} + \frac{(-\delta + 1/2 - c + \alpha' - \alpha)^2}{2}}{2} \omega(x_{\rho\mu \nu}x_\mu,x_\rho x_\nu,x_{\rho\mu}x_{\mu \nu} \vert \gamma - a - \alpha ')\ .
\end{array}$$
\noindent On peut encore simplifier les puissances de $\omega$ en utilisant la relation $b = P+1 - a - c$ mod($N$). Ainsi, la contribution totale de $\omega$ est : 
$$\omega^{\frac{c}{4} +  c(\gamma - \alpha') - \frac{ac}{2} - \frac{1}{8} - \alpha \alpha' -  \delta \alpha' + \frac{(\alpha')^2}{2} + \frac{\delta^2}{2}}\ .$$
\noindent D'o\`u la formule 
$$\begin{array}{l}
\bar{R}(\rho\mu,\bar{\mu},\mu \nu\vert b,c)_{\beta,\delta}^{\alpha,\gamma} = \omega^{\frac{c}{4}}(x_{\rho\mu}x_{\mu \nu})^P h(\frac{x_{\rho\mu}x_{\mu \nu}}{x_{\rho\mu \nu}x_\mu}) \ \left(\frac{\omega^{-1-\frac{1}{8}}(-1)^Pg(1)}{N} \right) \ \times \hspace{3.7cm}\\ \\ 
\hspace{1cm} \sum_{\alpha ',\delta ' =0}^{N-1} \omega^{c(\gamma - \alpha ') - \frac{ac}{2} + \alpha ' \delta ' }\ \omega(x_{\rho\mu \nu}x_\mu,x_\rho x_\nu,x_{\rho\mu}x_{\mu \nu} \vert \gamma - a,\alpha ')  \\
\hspace{7cm} \delta(\gamma + \delta ' - \beta)\ \omega^{\frac{\delta^2}{2}}\ \delta(\delta + \delta ')\ \omega^{-\alpha \alpha '}\ .
\end{array}$$
\noindent On a vu que $\vert g(1) \vert = N^{1/2}$. En considérant les matrices $S$ et $T$ et en prenant $\zeta = \omega^{\frac{9}{8}}(-1)^P \vert g(1) \vert/g(1)$, on obtient finalement la seconde relation de symm\'etrie.\hfill $\Box$

\medskip

\begin{prop} \label{unitarite1} Soit $(\rho,\mu,\nu)$ une suite régulière de représentations standard de $\mathcal{W}_N$. On a
$$\bar{R}(\rho^*,\mu^*,\nu^*\vert a,c)_{\gamma,\delta}^{\alpha,\beta} = \left( R(\rho,\mu,\nu\vert a,c)_{-\alpha,-\beta}^{-\gamma,-\delta} \right)^*,\nonumber$$
\noindent o\`u $\rho^*$ est la repr\'esentation complexe conjugu\'ee (D\'efinition {\rm \ref{repstandard}}).
\end{prop}

\noindent {\it Démonstration.} Par le lemme \ref{lem1}, pour tout $x \in \mc^*$ on a 
$$[x^*]^* = h(1/x^*)^*\ \frac{g(x^*)^*}{g(1)^*}\ .$$
\noindent Or
$$g(x^*)^* = \prod_{j=1}^{N-1} (1-x\omega^{-j})^{j/N} = (-x)^P\omega^{-\frac{(N-1)(2N-1)}{6}}\ g(1/x)\ .
$$
\noindent Donc $g(1)/g(1)^* = (-1)^P \omega^{\frac{(N-1)(2N-1)}{6}}$ et
$$\frac{\mathbb{[}(\frac{x_{\rho\mu \nu}x_{\mu}}{x_{\rho\mu}x_{\mu \nu}})^*\mathbb{]}^*}{h((\frac{x_{\rho\mu}x_{\mu \nu}}{x_{\rho\mu \nu}x_{\mu}})^*)^*} = h(\frac{x_{\rho\mu}x_{\mu \nu}}{x_{\rho\mu \nu}x_{\mu}}) \ \frac{g(1)(-1)^P\omega^{- (N-1)(2N-1)/6}}{g(1)^*} = h(\frac{x_{\rho\mu}x_{\mu \nu}}{x_{\rho\mu \nu}x_{\mu}})\ .$$
\noindent Alors
$$\begin{array}{l}
\left(\bar{R}(\rho^*,\mu^*,\nu^*\vert a,c)_{\gamma,\delta}^{\alpha,\beta}\right)^* = (x_{\rho\mu}x_{\mu \nu})^P h(\frac{x_{\rho\mu}x_{\mu \nu}}{x_{\rho\mu \nu}x_{\mu}}) \ \delta(\gamma + \delta - \beta) \ \times \hspace{2.5cm}\\ \\
\hspace{3.5cm} \omega^{c(\alpha - \gamma) - \frac{ac}{2} + \alpha \delta + \frac{\alpha^2}{2}}\ \frac{1}{\left(\omega(\frac{x_{\rho\mu \nu}^*x_{\mu}^*}{\omega},x_{\rho}^*x_{\nu}^*,{x_{\rho\mu}^*x_{\mu \nu}^* \vert \gamma - \alpha)} \right)^*}\ .
\end{array}$$
\noindent Maintenant, on vérifie de façon analogue au lemme \ref{inversion} que
$$\frac{1}{\left(\omega(\frac{x_{\rho\mu \nu}^*x_{\mu}^*}{\omega},x_{\rho}^*x_{\nu}^*,{x_{\rho\mu}^*x_{\mu \nu}^* \vert \gamma - \alpha)} \right)^*} = \omega(x_{\rho\mu \nu}x_{\mu},x_{\rho}x_{\nu},x_{\rho\mu}x_{\mu \nu} \vert \alpha - \gamma)\ ,$$
d'o\`u la conclusion.\hfill $\Box$

\bigskip

\noindent Rappelons que $Y=\omega^{1/2}XZ$ est la matrice de composantes $Y_{m,n} =\omega^{1/2 + n}\delta_{m,n+1}$ (cf. \S \ref{symboles}). Pour tout $k \in \mathbb{N}$, on a 
$$(Y^k)_{i,j} = \omega^{\frac{k^2}{2} + kj}\ \delta_{i,j+k},\quad (Y^{-k})_{i,j} = \omega^{\frac{-k^2}{2} - ki}\ \delta_{i+k,j}\ .$$

\begin{lem}\label{matrix}
Soit $(\rho,\mu,\nu)$ une suite régulière de représentations standard de $\mathcal{W}_N$. On a :
$$ R(\rho,\mu,\nu\vert a,c)_{\alpha,\beta}^{\gamma,\delta} = (x_{\rho\mu}x_{\mu \nu})^P \omega^{\frac{ac}{2}} \left( Y_1^{-a}Z_1^{-c} R(\rho,\mu,\nu) Z_1^cZ_2^{-a} \right)_{\alpha,\beta}^{\gamma,\delta}\ ,$$
\noindent o\`u $Y_1=Y \otimes id$, etc.
\end{lem}

\noindent {\it Démonstration.} On a
$$\begin{array}{l}
R(\rho,\mu,\nu\vert a,c)_{\alpha,\beta}^{\gamma,\delta}=  (x_{\rho\mu}x_{\mu \nu})^P\ \omega^{c(\gamma - \alpha) - \frac{ac}{2}} \ R(\rho,\mu,\nu)_{\alpha,\beta - a}^{\gamma - a ,\delta} \hspace{3.3cm} \\ \\
\hspace{1cm} = (x_{\rho\mu}x_{\mu \nu})^P\ h_{\rho,\mu,\nu}\ \omega^{c(\gamma - \alpha) - \frac{ac}{2}+\alpha \delta} \ \omega(x_{\rho\mu \nu}x_\mu,x_\rho x_\nu,x_{\rho\mu}x_{\mu \nu}\vert \gamma-a,\alpha)\ \times \\ 
\hspace{10cm} \delta(\gamma + \delta - \beta)
\end{array}$$
$$\begin{array}{l} 
\hspace{1cm} = (x_{\rho\mu}x_{\mu \nu})^P\ h_{\rho,\mu,\nu}\ \omega^{c(\gamma - \alpha) - \frac{ac}{2} +\alpha \delta} \ \omega(x_{\rho\mu \nu}x_\mu,x_\rho x_\nu,x_{\rho\mu}x_{\mu \nu}\vert \gamma,\alpha + a)\ \times \\ 
\hspace{8.1cm}  \omega^{\frac{\alpha^2}{2} - \frac{(\alpha + a)^2}{2}} \ \delta(\gamma + \delta - \beta)\\ \\
\hspace{1cm} = \omega^{c(\gamma - \alpha) +\alpha \delta - a(\frac{c+a}{2} + \alpha) - \gamma '\delta} \ \omega^{\gamma '\delta} \ \omega(x_{\rho\mu \nu}x_\mu,x_\rho x_\nu,x_{\rho\mu}x_{\mu \nu}\vert \gamma,\gamma ') \ \times \\ \\
\hspace{3.7cm} (x_{\rho\mu}x_{\mu \nu})^P\ h_{\rho,\mu,\nu}\ \delta(\gamma ' - \alpha - a) \ \delta(\gamma + \delta - \beta) \ \delta(\delta ' - \beta)\\ \\
\hspace{1cm} = (x_{\rho\mu}x_{\mu \nu})^P \ \omega^{-c\alpha - a(\frac{c+a}{2} + \alpha)} \  \omega^{c\gamma - a \delta} \delta(\gamma ' - \alpha - a)\ \times \\ \\
\hspace{3.3cm} \left(h_{\rho,\mu,\nu}\ \omega^{\gamma '\delta}\ \omega(x_{\rho\mu \nu}x_\mu,x_\rho x_\nu,x_{\rho\mu}x_{\mu \nu}\vert \gamma,\gamma ')\ \delta(\gamma + \delta - \delta ') \right) \\ \\
\hspace{1cm} = (x_{\rho\mu}x_{\mu \nu})^P \ \omega^{-c\gamma ' + \frac{ac}{2} - \frac{a^2}{2} - a\alpha} \ R(\rho,\mu,\nu)_{\gamma ',\delta '}^{\gamma,\delta} \ \omega^{c\gamma - a\delta} \ \delta(\gamma ' - \alpha - a) \\ \\
\hspace{1cm} = (x_{\rho\mu}x_{\mu \nu})^P\ \omega^{\frac{ac}{2}} \left( Y_1^{-a}Z_1^{-c} R(\rho,\mu,\nu) Z_1^cZ_2^{-a} \right)_{\alpha,\beta}^{\gamma,\delta}.
\end{array}$$ \hfill $\Box$

\begin{remark} {\rm  Le corollaire \ref{6j1} montre que pour tout triplet r\'egulier $(\rho,\mu,\nu) \in \mathcal{C}^3$, il existe des repr\'esentations normales $\pi_{(\rho,\mu,\nu)}$ et $\pi_{(\rho,\mu,\nu)}^*$ de $\bar{\mathcal{H}}_{\omega^{-1}}$ telles que
$$\begin{array}{l}
R(\rho,\mu,\nu \vert a,c) = (z_{\rho \mu}z_{\mu \nu})^P \ \omega^{\frac{ac}{2}}\ \times \hspace{5cm} \\ \\
\hspace{1.7cm} \pi_{(\rho,\mu,\nu)} \otimes \pi_{(\rho,\mu,\nu)}^* \biggl( ((\bar{E}')^{-a}(E')^c \otimes id)\ R_\omega \ ((E')^{-c} \otimes (E')^{a})\biggr) \\ \\
\bar{R}(\rho,\mu,\nu \vert a,c)= (z_{\rho \mu}z_{\mu \nu})^{\frac{N-1}{2N}} \ \omega^{\frac{-ac}{2}}\ \times \hspace{5cm} \\ \\
\hspace{1.7cm} \pi_{(\rho,\mu,\nu)} \otimes \pi_{(\rho,\mu,\nu)}^* \biggl( ((E')^{-c} \otimes (E')^{a})\ R_\omega^{-1} \ ((E')^c(\bar{E}')^{-a} \otimes id)\biggr).
\end{array}$$
\noindent Les repr\'esentations $\pi_{(\rho,\mu,\nu)}$ et $\pi_{(\rho,\mu,\nu)}^*$ sont obtenues comme dans la proposition \ref{CGdilocyc}.} 
\end{remark}

\subsection{L'équation du pentagone étendue} \label{EPE}

\noindent Dans cette section, on prouve une relation d'orthogonalit\'e et l'extension de l'équation du pentagone pour les c-$6j$-symboles. On renvoie à \cite{BB1} pour l'interpr\'etation g\'eom\'etrique de ces relations.

\medskip

\noindent Le lemme suivant est une conséquence immédiate du corollaire \ref{6j1}, mais peut aussi être obtenu à partir de la proposition \ref{CGdilocyc} et des relations de commutation de l'élément canonique $R_{\omega} \in {\rm End}\left( Q(\widehat{\bar{\mathcal{H}}_{\omega^{-1}}})^{\otimes 2} \right)$. 

\begin{lem} \label{commutation}
On a les relations de commutation suivantes :
$$\begin{array}{rcl}
(C_1) \quad  \quad R_{12}(\rho,\mu,\nu)\ Z_1Y_2 & = & Z_1Y_2\ R_{12}(\rho,\mu,\nu), \label{C1} \\
(C_2) \quad \quad \quad  R_{12}(\rho,\mu,\nu)\ Y_1 & = & Y_1Y_2\ R_{12}(\rho,\mu,\nu), \label{C2} \\
(C_3) \quad \quad  R_{12}(\rho,\mu,\nu)\ Z_1Z_2 & = & Z_2\ R_{12}(\rho,\mu,\nu)\ . \label{C3}
\end{array}$$
\end{lem}

\begin{prop} \label{EP1} Soit $(\rho,\mu,\nu,\upsilon)$ une suite régulière de représentations standard de $\mathcal{W}_N$. On a la \emph{relation du pentagone étendue (EP)} :
$$\begin{array}{l}
R_{12}^4(\rho,\mu,\nu \vert i,m-k)\ R_{13}^0(\rho,\mu \nu,\upsilon \vert j,l+m)\ R_{23}^0(\mu,\nu,\upsilon \vert k,l-i) = \hspace{3cm} \\ \\
\hspace{5cm} x_{\mu \nu}^{2P}\ R_{23}^1(\rho\mu,\nu,\upsilon \vert j+k,l)\ R_{12}^3(\rho,\mu,\nu\upsilon \vert i+j,m)\ .
\end{array}$$
\end{prop}

\noindent {\it Démonstration.} Avec le lemme \ref{matrix}, on voit que le membre de gauche s'écrit :
$$\begin{array}{l}
\omega^{\frac{(m-k)i}{2} + \frac{(l+m)j}{2} + \frac{(l-i)k}{2}} \ (x_{\rho\mu}x_{\mu \nu})^{P}(x_{\rho\mu \nu}x_{\mu \nu\upsilon})^{P}(x_{\mu \nu}x_{\nu\upsilon})^{P} \ \left( Y_1^{-i}Z_1^{k-m}R_{12}(\rho,\mu,\nu)\right.\\  \\ \hspace{0.5cm} \times\  \left. Z_1^{m-k}Z_2^{-i}Y_1^{-j}Z_1^{-l-m}R_{13}(\rho,\mu \nu,\upsilon)Z_1^{l+m}Z_3^{-j} Y_2^{-k}Z_2^{i-l}R_{23}(\mu,\nu,\upsilon)Z_2^{l-i}Z_3^{-k} \right)  .
\end{array}$$
\noindent Considérons d'abord le terme entre parenthèses. Comme $ZY = \omega YZ$, en éliminant $Z_2^{-i}$ et $Z_2^{i}$ et en réordonnant les autres, on trouve
$$\begin{array}{l}
\omega^{-j(m-k) +ik} \ \times \hspace{12cm} \\
Y_1^{-i}Z_1^{k-m}R_{12}(\rho,\mu,\nu)Y_1^{-j}Z_1^{-k-l}Y_2^{-k}R_{13}(\rho,\mu \nu,\upsilon)Z_1^{l+m}Z_2^{-l}Z_3^{-j}R_{23}(\mu,\nu,\upsilon)Z_2^{l-i}Z_3^{-k}\end{array}$$
\noindent Dans la chaîne de transformations qui suit, on indique qu'on utilise la relation $C_i,\ i=1,2,3$, du Lemme \ref{commutation} en écrivant $\stackrel{C_i}{=}$. La dernière expression est \'egale à
$$\begin{array}{l}
\stackrel{C_2}{=} \omega^{-j(m-k) +ik -j(k-m)}\ Y_1^{-i-j}Z_1^{k-m}Y_2^{-j}R_{12}(\rho,\mu,\nu)Z_1^{-k-l}Y_2^{-k}R_{13}(\rho,\mu \nu,\upsilon) \ \times \\
\hspace{7cm} Z_1^{l+m}Z_2^{-l}Z_3^{-j}R_{23}(\mu,\nu,\upsilon)Z_2^{l-i}Z_3^{-k}  \\ \\
\stackrel{C_1,C_3}{=} \omega^{ik}\ Y_1^{-i-j}Z_1^{-m}Y_2^{-j-k}Z_2^{-l}R_{12}(\rho,\mu,\nu)Z_2^lR_{13}(\rho,\mu \nu,\upsilon) \ \times \hspace{4cm} \\
\hspace{7cm} Z_1^{l+m}Z_2^{-l}Z_3^{-j}R_{23}(\mu,\nu,\upsilon)Z_2^{l-i}Z_3^{-k} \\ \\
= \omega^{ik}\ Y_1^{-i-j}Y_2^{-j-k}Z_1^{-m}Z_2^{-l}R_{12}(\rho,\mu,\nu)R_{13}(\rho,\mu \nu,\upsilon) \ \times \hspace{4cm} \\
\hspace{7cm} Z_1^{l+m}Z_3^{-j}R_{23}(\mu,\nu,\upsilon)Z_2^{l-i}Z_3^{-k}  
\end{array}$$
$$\begin{array}{l}
\stackrel{C_3}{=} \omega^{ik}\ Y_1^{-i-j}Y_2^{-j-k}Z_1^{-m}Z_2^{-l}Z_3^{l+m}R_{12}(\rho,\mu,\nu)R_{13}(\rho,\mu \nu,\upsilon)\ \times \hspace{4cm} \\
\hspace{7cm} Z_3^{-j-l-m}R_{23}(\mu,\nu,\upsilon)Z_2^{l-i}Z_3^{-k}  \\ \\
\stackrel{C_3}{=} \omega^{ik}\ Y_1^{-i-j}Y_2^{-j-k}Z_1^{-m}Z_2^{-l}Z_3^{l+m}R_{12}(\rho,\mu,\nu)R_{13}(\rho,\mu \nu,\upsilon) \ \times \hspace{4cm} \\
\hspace{7cm} R_{23}(\mu,\nu,\upsilon)Z_2^{-j-m-i}Z_3^{-j-l-m-k}.
\end{array}$$
\noindent Par la relation du pentagone (\ref{P1}), ceci vaut 
$$\begin{array}{l}
\omega^{ik}  \ Y_1^{-i-j}Y_2^{-j-k}Z_1^{-m}Z_2^{-l}Z_3^{l+m} R_{23}(\rho\mu,\nu,\upsilon)R_{12}(\rho,\mu,\nu\upsilon)Z_2^{-j-m-i}Z_3^{-j-l-m-k}\\ \\
 = \omega^{ik}  \ Y_1^{-i-j}Y_2^{-j-k}Z_1^{-m}Z_2^{-l}Z_3^{l+m} R_{23}(\rho\mu,\nu,\upsilon)Z_3^{-l-m}R_{12}(\rho,\mu,\nu\upsilon)Z_2^{-j-m-i}Z_3^{-j-k}\\ \\
 \stackrel{C_3}{=}\omega^{ik}\ Y_1^{-i-j}Y_2^{-j-k}Z_1^{-m}Z_2^{-l}R_{23}(\rho\mu,\nu,\upsilon)Z_2^{l+m}R_{12}(\rho,\mu,\nu\upsilon)Z_2^{-j-m-i}Z_3^{-j-k}\\ \\
 \stackrel{C_3}{=}\omega^{ik} \  Y_1^{-i-j}Y_2^{-j-k}Z_1^{-m}Z_2^{-l}R_{23}(\rho\mu,\nu,\upsilon)Z_2^{l}Z_3^{-j-k}R_{12}(\rho,\mu,\nu\upsilon)Z_1^{m}Z_2^{-i-j}\end{array}$$
$$\begin{array}{l}
 = \omega^{ik} \ Y_2^{-j-k}Z_2^{-l}R_{23}(\rho\mu,\nu,\upsilon)Z_2^{l}Z_3^{-j-k} Y_1^{-i-j}Z_1^{-m}R_{12}(\rho,\mu,\nu\upsilon)Z_1^{m}Z_2^{-i-j} \\ \\
 = \omega^{ik}\ \omega^{\frac{-(i+j)m}{2}}\ \omega^{\frac{-(j+k)l}{2}}\ (x_{\rho\mu \nu}x_{\nu\upsilon}x_{\rho\mu}x_{\mu \nu\upsilon})^{-P}\ \times \hspace{4cm} \\
\hspace{5cm}  R_{23}(\rho\mu,\nu,\upsilon \vert j+k,l)\ R_{12}(\rho,\mu,\nu\upsilon \vert i+j,m)\ .
\end{array}$$
\noindent Or 
$$\omega^{\frac{(m-k)i}{2} + \frac{(l+m)j}{2} + \frac{(l-i)k}{2} + ik} = \omega^{\frac{(i+j)m}{2} + \frac{(j+k)l}{2}}\ ,$$ 
\noindent d'o\`u le r\'esultat. \hfill $\Box$

\begin{prop} \label{orth1} Soit $(\rho,\mu,\nu)$ une suite régulière de représentations standard de $\mathcal{W}_N$. On a la relation \emph{d'orthogonalité}:
$$R(\rho,\mu,\nu\vert a,c)\bar{R}(\rho,\mu,\nu \vert -a,-c) = (x_{\rho\mu}x_{\mu \nu})^{2P} id \otimes id\ .$$
\end{prop}

\noindent {\it Démonstration.} On a
$$\begin{array}{l}
x_{\rho\mu}^{2P}R_{12}^4(\rho,\mu,\nu \vert i,m-k)R_{13}^2(\rho,\mu \nu,\upsilon \vert j,l+m)R_{23}^0(\mu,\nu,\upsilon \vert k,l-i) = \\ \\
R_{12}^4(\rho,\mu,\nu \vert i,m-k)\bar{R}_{12}^4(\rho,\mu,\nu \vert -i,k-m)R_{23}^1(\rho\mu,\nu,\upsilon \vert j+k,l)R_{12}^3(\rho,\mu,\nu\upsilon \vert i+j,m).
\end{array}$$
\noindent La relation du pentagone étendue pour le membre de gauche donne
$$\begin{array}{l}
(x_{\rho\mu}x_{\mu \nu})^{2P} R_{23}^1(\rho\mu,\nu,\upsilon \vert j+k,l)R_{12}^3(\rho,\mu,\nu\upsilon \vert i+j,m) = \\
\hspace{3cm} R_{12}^4(\rho,\mu,\nu \vert i,m-k)\bar{R}_{12}^4(\rho,\mu,\nu \vert -i,k-m) \\
\hspace{5.5cm}R_{23}^1(\rho\mu,\nu,\upsilon \vert j+k,l)R_{12}^3(\rho,\mu,\nu\upsilon \vert i+j,m) \ .
\end{array}$$ 
\noindent Comme les c-$6j$-symboles sont des opérateurs inversibles, ceci achève la démonstration. \hfill $\Box$

\end{document}